\documentclass[11pt]{article}

\usepackage[utf8]{inputenc}
\usepackage[a4paper,margin=1in]{geometry}
\usepackage{amsmath,amssymb,amsthm,mathtools,bm}
\usepackage[expansion=false]{microtype}
\usepackage[round,authoryear,sort]{natbib}
\makeatletter
\newcommand\refnumfmt[1]{[#1]}
\renewcommand\NAT@biblabel[1]{\refnumfmt{#1}}
\renewcommand\NAT@bibsetup[1]{%
  \settowidth\labelwidth{\refnumfmt{#1}}%
  \setlength{\leftmargin}{\labelwidth}\addtolength{\leftmargin}{\labelsep}%
  \setlength{\itemindent}{\z@}%
  \setlength{\itemsep}{\bibsep}\setlength{\parsep}{\z@}}
\makeatother
\usepackage{booktabs,tabularx,longtable,array}
\usepackage{enumitem}
\usepackage{xcolor}
\usepackage{tikz}
\usetikzlibrary{arrows.meta,positioning,fit,backgrounds,calc,decorations.pathreplacing,intersections}
\definecolor{figblue}{RGB}{46,94,150}
\definecolor{figorange}{RGB}{205,108,38}
\usepackage[most]{tcolorbox}
\usepackage[colorlinks=true,linkcolor=blue!45!black,citecolor=green!30!black,urlcolor=blue!55!black]{hyperref}
\usepackage[capitalise,noabbrev]{cleveref}

\setlist{itemsep=2pt,topsep=4pt}
\newcolumntype{L}{>{\raggedright\arraybackslash}X}
\newcolumntype{P}[1]{>{\raggedright\arraybackslash}p{#1}}
\allowdisplaybreaks
\theoremstyle{plain}
\newtheorem{theorem}{Theorem}[section]
\newtheorem{proposition}[theorem]{Proposition}
\newtheorem{lemma}[theorem]{Lemma}

\theoremstyle{definition}
\newtheorem{definition}[theorem]{Definition}

\newtheorem{example}[theorem]{Example}
\theoremstyle{remark}
\newtheorem{remark}[theorem]{Remark}
\crefname{assumption}{Assumption}{Assumptions}

\newenvironment{keybox}[1][]{\par\medskip\noindent\rule{\linewidth}{0.3pt}\par\nopagebreak\smallskip\noindent{\scshape #1.}\enspace\ignorespaces}{\par\nopagebreak\smallskip\noindent\rule{\linewidth}{0.3pt}\par\medskip}
\makeatletter\newenvironment{openprob}[1][]{\par\medskip\noindent\rule{\linewidth}{0.3pt}\par\nopagebreak\smallskip\noindent{\scshape #1}\par\nopagebreak\@beginparpenalty=\@M}{\par\nopagebreak\smallskip\noindent\rule{\linewidth}{0.3pt}\par\medskip}\makeatother
\theoremstyle{definition}
\newtheorem{problem}{Problem}[part]
\crefname{problem}{Problem}{Problems}
\makeatletter
\renewcommand\paragraph{\@startsection{paragraph}{4}{\z@}{1.25ex \@plus .5ex \@minus .2ex}{-0.6em}{\normalfont\scshape}}
\makeatother

\newcommand{\R}{\mathbb{R}}

\newcommand{\E}{\mathbb{E}}
\renewcommand{\P}{\mathbb{P}}
\newcommand{\Q}{\mathbb{Q}}
\newcommand{\cP}{\mathcal{P}}
\newcommand{\cH}{\mathcal{H}}
\newcommand{\cF}{\mathcal{F}}
\newcommand{\cL}{\mathcal{L}}
\newcommand{\cN}{\mathcal{N}}

\newcommand{\M}{\mathsf{M}}
\newcommand{\dd}{\mathrm{d}}
\newcommand{\Id}{\mathrm{I}}
\newcommand{\one}{\mathbf{1}}
\newcommand{\law}{\operatorname{Law}}
\newcommand{\KL}{\operatorname{KL}}
\newcommand{\TV}{\operatorname{TV}}
\newcommand{\W}{\mathsf{W}}
\newcommand{\tr}{\operatorname{tr}}
\newcommand{\divg}{\nabla\!\cdot}

\newcommand{\Lip}{\operatorname{Lip}}

\newcommand{\pf}{_{\#}}
\newcommand{\D}{\mathcal{D}}
\newcommand{\Sk}{\delta}
\newcommand{\Dom}{\operatorname{Dom}}
\newcommand{\gammaM}{\gamma}
\newcommand{\snr}{\operatorname{SNR}}

\newcommand{\argmin}{\operatorname*{arg\,min}}

\newcommand{\inner}[2]{\langle #1,#2\rangle}
\newcommand{\norm}[1]{\lVert #1\rVert}
\newcommand{\abs}[1]{\lvert #1\rvert}
\newcommand{\Rref}{\mathsf{R}}

\title{Geometry and Stochastic Analysis of\\ Flow and Diffusion Models: A Survey}
\author{%
\begin{tabular}[t]{c}Ehsan Mirafzali\\[2pt]\small\href{mailto:smirafza@ucsc.edu}{\texttt{smirafza@ucsc.edu}}\end{tabular}\qquad
\begin{tabular}[t]{c}Razvan Marinescu\\[2pt]\small\href{mailto:ramarine@ucsc.edu}{\texttt{ramarine@ucsc.edu}}\end{tabular}\\[6pt]
\small Department of Computer Science and Engineering, University of California Santa Cruz}
\date{}

\hypersetup{pdftitle={Geometry and Stochastic Analysis of Flow and Diffusion Models: A Survey},pdfauthor={Ehsan Mirafzali and Razvan Marinescu}}

\begin{document}
\maketitle

\begin{abstract}
Flow-based and diffusion-based models have become the principal means of learning and sampling high-dimensional probability distributions, from images and video to proteins, molecules and physical fields, and increasingly language. Each carries a simple law to the law of the data along a curve of probability measures, and behind the variety of their names lie three geometries, those of the state space on which the data live, of the space of probability measures through which the curve passes, and of the space of paths traced by the generating process. A single analytic idea runs through the subject. Every score identity used to train a diffusion model is an integration by parts, on the noise, on the data or on Wiener space, and Malliavin calculus extends this to nonlinear and infinite-dimensional diffusions. We develop these foundations with full definitions and follow them through normalising flows, flow matching and few-step maps, diffusion models and their convergence theory, conditioning, Schr\"odinger bridges and stochastic control, and their extension to manifolds, Lie groups, discrete and function spaces and jump processes. We close with statistical physics, renormalisation and, more speculatively, holography, in which the radial direction of anti-de Sitter space can play the role of flow time. A table of the main methods and a collection of open problems complete the survey.
\end{abstract}

\tableofcontents
\newpage

\section{Introduction}\label{sec:intro}

A generative model learns to produce new samples that resemble a given collection of data. The models surveyed here do so by transport. They move probability mass from a simple law to the law of the data, and the mathematics of such motion is the subject of this survey.

\subsection{Generative modelling as transport of measures}

Let $\mu_{\rm data}$ be a probability measure on a state space $\M$, known only through samples $x^{(1)},\dots,x^{(n)}$. A generative model is a procedure that turns samples of a simple reference measure $\mu_{\rm ref}$, usually a standard Gaussian and sometimes a uniform or heat-kernel measure, into samples of an approximation of $\mu_{\rm data}$. Every model in this survey does so by moving probability mass along a curve $(\rho_t)_{t\in[0,1]}$ in the space $\cP(\M)$ of probability measures, or, for flows built from finitely many layers, along a finite sequence of such laws. It specifies the curve with $\rho_0=\mu_{\rm ref}$ and $\rho_1\approx\mu_{\rm data}$ or the reverse, depending on convention, together with a mechanism that realises the curve by transporting individual points. Two mechanisms dominate.

One is deterministic transport. A time-dependent vector field $v_t$ on $\M$ generates a flow $\phi_t$, with $\partial_t\phi_t=v_t\circ\phi_t$, and the curve $\rho_t=(\phi_t)\pf\rho_0$ solves the continuity equation $\partial_t\rho_t+\divg(\rho_t v_t)=0$. Normalising flows, continuous normalising flows, flow matching, rectified flow, probability-flow ODEs and flow-map models are of this kind. The other is stochastic transport. A Markov process $(X_t)$ with time-dependent generator $\cL_t$ realises the curve through the Kolmogorov forward equation $\partial_t\rho_t=\cL_t^{*}\rho_t$. Score-based diffusion models, Schr\"odinger bridges, stochastic interpolants, diffusion samplers, and discrete and jump-process models are of this kind.

The two mechanisms are closely tied. Under mild regularity conditions, every diffusion has a deterministic flow, the probability-flow ODE, with the same one-time marginals, and every deterministic flow can be dressed with noise without changing its marginals, provided the drift is corrected by a score term. The art lies in choosing the curve $(\rho_t)$, the regression by which the velocity or generator is learned from samples with a computable target, and the scheme by which the learned dynamics are integrated accurately and quickly.

\subsection{Three geometries and one calculus}

\begin{figure}[t]
\centering
\begin{tikzpicture}[
  box/.style={draw,rounded corners=3pt,align=center,inner sep=6pt,text width=3.5cm,minimum height=3.65cm,fill=#1},
  box/.default=white, >=Stealth]
\node[box=blue!5] (state) {{\scshape State space}\\ $\M$\\[3pt]\small points $x$; vector fields; Brownian motion and heat kernels on $\M$\\[3pt]\footnotesize Part V};
\node[box=green!5,right=1.9cm of state] (meas) {{\scshape Measure space}\\ $\cP(\M)$\\[3pt]\small laws $\rho_t$; continuity and Fokker--Planck equations; OT and Fisher--Rao geometry\\[3pt]\footnotesize Parts II, IV};
\node[box=orange!6,right=1.9cm of meas] (path) {{\scshape Path space}\\ $C([0,T];\M)$\\[3pt]\small laws of trajectories; Girsanov; time reversal; Doob; control; Malliavin\\[3pt]\footnotesize Parts III, IV};
\draw[->,thick] (state) -- node[above,font=\footnotesize]{law} (meas);
\draw[->,thick] (path) -- node[above,font=\footnotesize]{marginals} (meas);
\node[draw,rounded corners=3pt,fill=gray!8,below=0.6cm of meas,text width=15.2cm,align=center,inner sep=5pt] (found) {\small{\scshape Part I:} transport (\S\ref{sec:transport}), path-space analysis (\S\ref{sec:pathspace}), integration by parts (\S\ref{sec:ibp})};
\node[draw,rounded corners=3pt,fill=violet!6,above=0.6cm of meas,text width=15.2cm,align=center,inner sep=5pt] (phys) {\small{\scshape Parts VI--VII:} statistical physics, renormalisation and holography; a unified view and open problems};
\end{tikzpicture}
\caption{The three geometries of generative dynamics. A model is specified by objects at one or more levels; most results in the literature translate a statement at one level into a statement at another.}
\label{fig:three}
\end{figure}
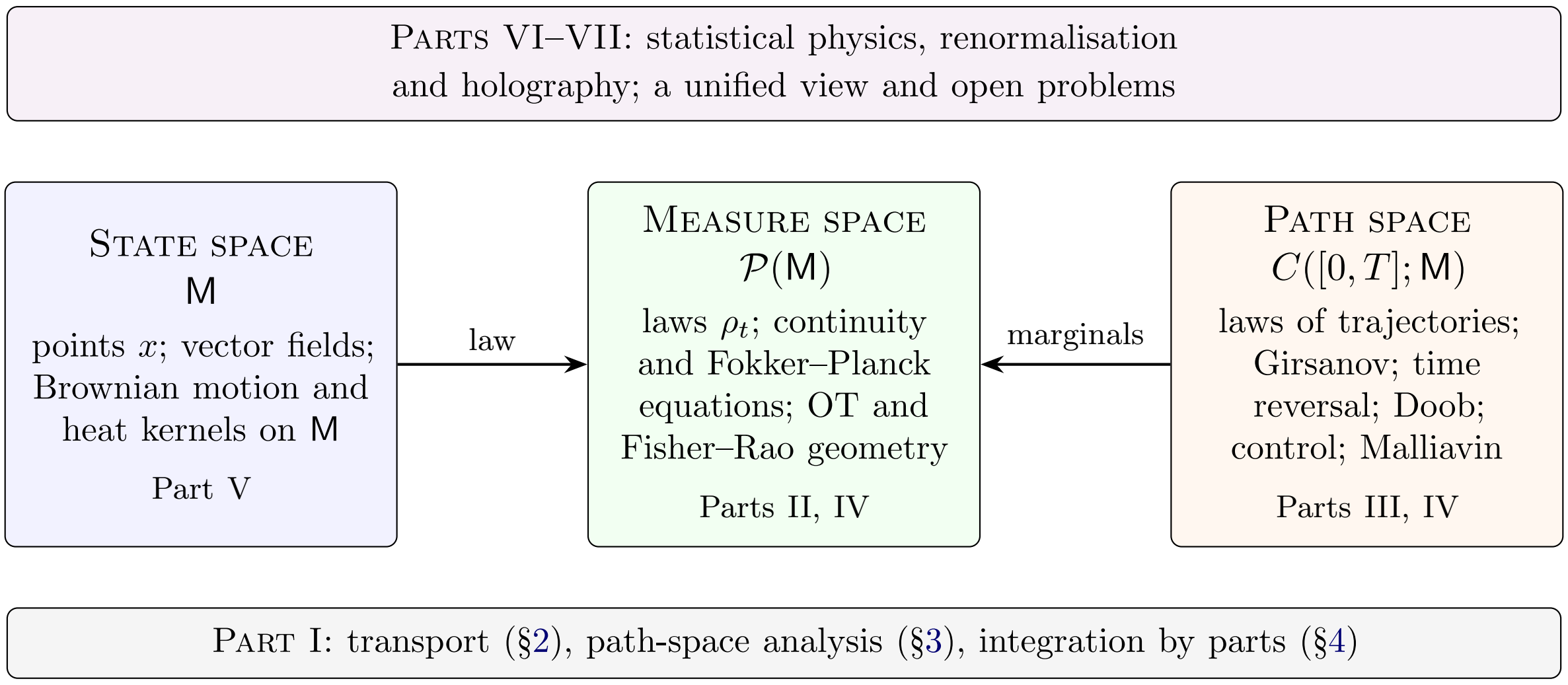

The objects just introduced live at three levels, and each level carries its own geometry and its own tools (\cref{fig:three}). The first level is the state space $\M$, which may be Euclidean space, a Riemannian manifold, a Lie group, a simplex or a finite set, or an infinite-dimensional Hilbert space. Its geometry fixes what a velocity field, a Brownian motion, a heat kernel and a score mean. The second is the measure space $\cP(\M)$, whose points are laws and whose curves have velocities and lengths. Optimal transport equips $\cP(\M)$ with the Wasserstein distance and with a formal Riemannian structure, Otto's calculus. Information geometry supplies the Fisher--Rao metric. Combining the two gives the Wasserstein--Fisher--Rao geometry, and restricting Wasserstein velocities to a reproducing-kernel Hilbert space gives the Stein geometry. The third is the path space $C([0,T];\M)$, or a Skorokhod space for jump processes, whose points are entire trajectories. On path space, relative entropy decomposes by Girsanov's theorem, and time reversal and conditioning are changes of measure. Schr\"odinger bridges and control problems are problems of entropy minimisation on path space, and Malliavin calculus is its differential calculus.

Most results in the field are translations between these levels. Flow matching turns a desired curve in $\cP(\M)$ into a regression for a velocity field on $\M$. The reverse-time SDE of a diffusion model is a statement about path measures whose drift involves the score of one-time marginals. Convergence bounds for diffusion models bound a distance between terminal laws by a relative entropy between path measures, and Schr\"odinger bridges minimise entropy on path space subject to constraints on two marginals. Making these translations explicit is the organising aim of this survey.

The second organising theme is analytic. Diffusion models never observe the score $\nabla\log\rho_t$; they regress onto a random target whose conditional expectation is the score. Each such identity is an integration by parts, whether on the Gaussian noise (Tweedie's formula and denoising score matching), on the data distribution (target score matching), or on the Wiener space that drives a general SDE (Bismut-type and Malliavin formulas). Unlike the first two, the last needs neither a closed-form transition kernel nor the data score, so it extends to nonlinear, state-dependent and infinite-dimensional forward processes. It also shows that the classical denoising target is the special case of a deterministic covering vector field on Wiener space (\cref{ex:linear}). \Cref{sec:ibp} develops this viewpoint, \cref{sec:score-estimation} uses it to compare score estimators, and \cref{sec:infinite} carries it to Hilbert spaces and SPDEs.

\subsection{Scope and organisation}

The survey covers generative models that transport measures by flows or by Markov processes in continuous or discrete time. It treats normalising flows viewed as structured diffeomorphisms, continuous normalising flows, flow matching and its variants, stochastic interpolants, score-based and denoising diffusion models, discrete-state and jump-process diffusions, Schr\"odinger bridges and bridge matching, flow maps and consistency models, diffusion samplers for unnormalised densities, conditioning and fine-tuning, and their theory. Variational autoencoders and generative adversarial networks are outside scope; latent-space models appear only where the latent dynamics are themselves flows or diffusions. Network architectures are discussed only when they encode geometry, chiefly equivariance. Applications appear as illustrations rather than as a catalogue. Several surveys already catalogue applications \citep{yang2023diffusion,cao2024survey,croitoru2023diffusion}, and tutorials and overviews exist for individual frameworks and for some of the links between them \citep{luo2022understanding,tang2024score,chen2024overview,chan2024tutorial,nakkiran2024step,lipman2024flow,holderrieth2025introduction,lai2025principles,tang2026foundations,vayer2026notes}. Normalising flows are reviewed by \citet{papamakarios2021normalizing} and \citet{kobyzev2020normalizing}, and deep generative models more broadly by \citet{bondtaylor2022deep} and \citet{ruthotto2021introduction}.

Part~I sets out the three geometries, developing those of the space of measures and of path space together with the integration-by-parts calculus; particular state spaces are taken up in Part~V. Part~II treats deterministic flows, from triangular maps to flow matching, one-step flow maps and gradient flows of discrepancies. Part~III treats diffusion models through time reversal, score estimation, samplers, convergence theory and the geometry of what is learned. Part~IV treats conditioning, Schr\"odinger bridges and stochastic optimal control, all as changes of measure on path space. Part~V leaves Euclidean space and Brownian noise for manifolds, Lie groups, discrete and statistical-manifold state spaces, Hilbert spaces, and L\'evy and jump drivers. Part~VI turns to physics, with phase transitions and memorisation in high dimension, nonequilibrium thermodynamics, stochastic quantisation, renormalisation as optimal transport, and, more speculatively, holographic generative flows. Part~VII gathers the main methods in a single table, collects the open problems and concludes. The appendices give a dictionary of parameterisations and a table of forward kernels on standard spaces.

\subsection{How to read this survey}\label{sec:reading}

The survey can be read in order, but it need not be. Part~I is the foundation for everything that follows, and \cref{tab:reading} suggests routes through the rest for four kinds of reader. The Malliavin material (\cref{sec:malliavin} and the subsection after it) is the most technical part of the survey and can be skipped on a first reading, since later sections use only its conclusions, which \cref{tab:ibp} collects. Parts~II to VI each end with open problems, which \cref{sec:outlook} gathers, and the table of methods in \cref{sec:zoo} doubles as an index.

\begin{table}[t]
\centering\small
\begin{tabularx}{\textwidth}{@{}P{3.4cm}LP{3.5cm}@{}}
\toprule
Reader & Suggested route & Can wait \\
\midrule
New to the area & \Cref{sec:transport}, the Gaussian part of \cref{sec:ibp}, then \cref{sec:simfree,sec:flowmaps,sec:diffusion,sec:samplers,sec:conditioning}; \cref{app:dictionary} as a reference & Malliavin calculus, \cref{sec:infinite,sec:levy}, Part~VI \\
\addlinespace[2pt]
Statistician or learning theorist & Part~I, then \cref{sec:score-estimation,sec:theory,sec:generalisation,sec:sb}, and the theory of flow matching in \cref{sec:simfree} & architectures, Part~VI \\
\addlinespace[2pt]
Probabilist or analyst & Part~I, then \cref{sec:score-estimation,sec:sb,sec:control,sec:infinite,sec:levy} & samplers, applications \\
\addlinespace[2pt]
Physicist & \Cref{sec:transport,sec:pathspace,sec:diffusion,sec:control}, then Part~VI & convergence theory \\
\bottomrule
\end{tabularx}
\caption{Four routes through the survey. Every route starts from Part~I.}
\label{tab:reading}
\end{table}

\subsection{Notation and conventions}\label{sec:notation}

\Cref{tab:notation} fixes notation. Time runs forward from data to noise in diffusion models ($X_0\sim\mu_{\rm data}$, $\law(X_T)\approx\mu_{\rm ref}$), and generation integrates the reverse process; for flows and interpolants we follow the more common convention $x_0\sim\mu_{\rm ref}$, $x_1\sim\mu_{\rm data}$, and we say so whenever the two meet. Densities are with respect to Lebesgue or Riemannian volume unless stated otherwise. We often use the same letter for a measure and its density, and write $p$ for the density of a law $\rho$ where the distinction matters. A few letters carry more than one meaning, which context separates; $T$ is a map or a time horizon, $h$ is a step size, a direction in a Hilbert space or a Doob $h$-function, and $\sigma$ is a diffusion coefficient whereas $\sigma_t$ is a noise scale. Symbols used only locally are defined where they appear.

\begin{table}[t]
\centering\small
\begin{tabularx}{\textwidth}{@{}lL@{}}
\toprule
Symbol & Meaning \\
\midrule
$\M$, $d$ & state space and its dimension; $\R^d$ unless stated \\
$\cP(\M)$, $\cP_2(\M)$ & probability measures; those with finite second moment \\
$\mu_{\rm data}$, $\mu_{\rm ref}$, $\gamma_d$ & data law, reference law, standard Gaussian measure $\cN(0,\Id)$ on $\R^d$ \\
$\hat\mu_n$ & empirical measure of $n$ training samples $x^{(1)},\dots,x^{(n)}$ \\
$T\pf\mu$ & pushforward, $(T\pf\mu)(A)=\mu(T^{-1}(A))$ \\
$\rho_t$, $p_t$ & marginal law at time $t$ and its density \\
$v_t$, $u_t$ & velocity fields; $u_t(\cdot\mid z)$ conditional velocity \\
$s_t=\nabla\log p_t$, $s_\theta$ & score of the marginal at time $t$; a network that approximates it \\
$f$, $g$, $\sigma$, $a=\sigma\sigma^\top$ & drift, scalar diffusion, matrix diffusion, diffusion matrix of an SDE \\
$\alpha_t$, $\sigma_t$, $\lambda_t$ & signal scale, noise scale and log-SNR $\lambda_t=\log(\alpha_t^2/\sigma_t^2)$ of a Gaussian path $X_t=\alpha_tX_0+\sigma_t\epsilon$ \\
$\epsilon$, $\epsilon_\theta$ & standard Gaussian noise of a Gaussian path, independent of $X_0$; a network that predicts it \\
$T$, $N$ & time horizon; number of steps of a sampler on $[0,T]$, with step size $h=T/N$ \\
$W_t$, $\overline{W}_t$ & Brownian motion; a reverse-time Brownian motion \\
$\P$, $\Q$, $\Rref$ & path measures (laws on $C([0,T];\M)$); $\Rref$ a reference process, $\Rref^{xy}$ its bridge from $x$ to $y$ \\
$\KL$, $\TV$, $\W_p$ & relative entropy, total variation, $p$-Wasserstein distance \\
$\cL_t$, $\cL_t^*$ & generator of a Markov process and its formal adjoint \\
$H=L^2([0,T];\R^m)$ & Hilbert space of the isonormal process of $m$-dimensional Brownian motion, isometric to its Cameron--Martin space via $h\mapsto\int_0^\cdot h$ \\
$\D$, $\Sk$, $\gammaM_F$ & Malliavin derivative, Skorokhod integral, Malliavin covariance of $F$ \\
$h(t,x)$ & Doob $h$-function; $a\nabla\log h$ the guidance drift ($g^2\nabla\log h$ for scalar diffusion $g$) \\
\bottomrule
\end{tabularx}
\caption{Notation used throughout. Where a section keeps the notation of a source paper, for instance $\alpha$, $\beta$, $\gamma$ for stochastic interpolants, this is said where it occurs.}
\label{tab:notation}
\end{table}

\part{Foundations: three geometries and one calculus}

\section{Measure transport and the geometry of \texorpdfstring{$\cP(\M)$}{P(M)}}\label{sec:transport}

Generative modelling asks for a way of carrying a simple law to the law of the data. When the carrier is a map, or a family of maps indexed by time, the question acquires a geometry, which we develop in three stages. We begin with a single map and the two standard ways of choosing one. We then let the map depend on time, which turns the space of probability measures into a space with tangent vectors, gradients and distances. We close with two constructions, entropic regularisation and the heat flow, through which this geometry meets the stochastic processes of \cref{sec:pathspace}.

\subsection{Pushforwards and optimal transport}

The simplest carrier is a single invertible map, and the first question is how it changes densities.

\paragraph{Pushforwards and change of variables.}

If $T:\R^d\to\R^d$ is a $C^1$ diffeomorphism and $\rho$ has density $p$, then $T\pf\rho$ has density
\begin{equation}\label{eq:cov}
\begin{aligned}
p_T(y)&=p\big(T^{-1}(y)\big)\,\abs{\det \nabla T^{-1}(y)},\\
\log p_T\big(T(x)\big)&=\log p(x)-\log\abs{\det\nabla T(x)}.
\end{aligned}
\end{equation}
Normalising flows (\cref{sec:nf}) are built on this formula. They parameterise $T$ so that both $T^{-1}$ and $\log\abs{\det\nabla T}$ are cheap to evaluate, and fit $T$ by maximum likelihood. Many maps push one law to another, and two existence results single out which of them to look for. The first is constructive and asymmetric; the second is variational and canonical.

\paragraph{Triangular maps.} Given absolutely continuous $\mu,\nu$ on $\R^d$, there is a map $T$ with $T\pf\mu=\nu$ whose $k$-th component depends only on $(x_1,\dots,x_k)$ and is increasing in $x_k$, and it is unique up to $\mu$-null sets \citep{bogachev2005triangular}; this is the Knothe--Rosenblatt rearrangement. It is assembled from successive conditional cumulative distribution functions, so it depends on the order in which the coordinates are listed, which is the asymmetry referred to above. \citet{rosenblatt1952remarks} used these to map an absolutely continuous law to the uniform law on the unit cube, and $T$ is the map of this kind for $\mu$ followed by the inverse of the one for $\nu$. Its Jacobian is therefore triangular and its determinant is the product of the diagonal entries; this is the structural fact behind autoregressive flows. \citet{bogachev2005triangular} develop the general theory of such maps, including the change-of-variables formula $\rho_\mu=(\rho_\nu\circ T)\prod_k\partial_{x_k}T_k$, in which $\rho_\mu,\rho_\nu$ are the densities of $\mu,\nu$ and $T_k$ is the $k$-th component of $T$. \citet{marzouk2016sampling} survey their use for sampling and inference, and \citet{irons2022triangular} give statistical rates for estimating them.

\paragraph{Optimal transport.}

Optimal transport chooses, among all maps pushing $\mu$ to $\nu$, one of least cost. For quadratic cost the Monge problem is $\inf\{\int\abs{x-T(x)}^2\dd\mu: T\pf\mu=\nu\}$, and its Kantorovich relaxation optimises over couplings rather than maps,
\begin{equation}
\W_2^2(\mu,\nu)=\inf_{\pi\in\Pi(\mu,\nu)}\int\abs{x-y}^2\,\dd\pi(x,y),
\end{equation}
where $\Pi(\mu,\nu)$ is the set of couplings \citep{villani2009optimal,santambrogio2015optimal,peyre2019computational}. When $\mu$ has a density the optimal coupling is carried by a map, and that map has a simple form.

\begin{theorem}[Brenier \citep{brenier1991polar}; see \citep{villani2009optimal,santambrogio2015optimal}]\label{thm:brenier}
If $\mu,\nu\in\cP_2(\R^d)$ and $\mu$ is absolutely continuous, there is a convex function $\varphi$, whose gradient is uniquely determined $\mu$-almost everywhere, such that $T=\nabla\varphi$ pushes $\mu$ to $\nu$; the coupling $(\Id,\nabla\varphi)\pf\mu$ is the unique optimal coupling.
\end{theorem}
\citet{brenier1991polar} proved this for $\mu$ with a density that vanishes outside a bounded, smooth, connected domain and is bounded away from zero and infinity on it, and remarked that the result should extend to unbounded domains under a moment condition. The statement above, for every absolutely continuous $\mu\in\cP_2(\R^d)$, is the form found in \citet{villani2009optimal}, while \citet{santambrogio2015optimal} proves it on compact domains. \citet{mccann1995existence} showed that, without any moment assumption, the monotone map $\nabla\varphi$ still exists and is unique.
When both measures have densities $p$ and $q$, the potential $\varphi$ solves the Monge--Amp\`ere equation
\begin{align*}
p(x)=q\big(\nabla\varphi(x)\big)\det\nabla^2\varphi(x)
\end{align*}
in a weak sense. Brenier maps are the targets of convex-potential flows \citep{huang2021convex} and of optimal flow matching \citep{kornilov2024optimal}. How readily such a map can be learned depends on its regularity, and the cleanest regularity statement is Caffarelli's contraction theorem. If $\nu=e^{-V}$ with $\nabla^2V\succeq\Id$, the Brenier map from $\gamma_d$ to $\nu$ is $1$-Lipschitz \citep{caffarelli2000monotonicity}.

\subsection{Geometries of the space of measures}

A single map records where mass ends up but not how it travels. Letting the map depend on time exposes the path the mass follows, and the space of probability measures then acquires a formal Riemannian structure.

\paragraph{The continuity equation.}

A curve $(\rho_t)$ and a velocity field $(v_t)$ satisfy the continuity equation if
\begin{equation}\label{eq:continuity}
\partial_t\rho_t+\divg(\rho_t v_t)=0\quad\text{in the sense of distributions}.
\end{equation}
If $v$ is locally Lipschitz in space with suitable integrability, the flow $\phi_t$ of $v$ exists and $\rho_t=(\phi_t)\pf\rho_0$; without Lipschitz regularity, provided $\int_0^T\!\int\abs{v_t}^p\,\dd\rho_t\,\dd t<\infty$ for some $p>1$, the superposition principle still represents $\rho_t$ as the time marginals of a measure on integral curves of $v$ \citep{ambrosio2008gradient}. Many velocity fields generate the same curve. Adding to $v_t$ any field $w_t$ with $\divg(\rho_t w_t)=0$ changes the flow but not the marginals. Benamou and Brenier showed that the cheapest choice recovers optimal transport.

\begin{theorem}[Benamou--Brenier \citep{benamou2000computational,ambrosio2008gradient}]\label{thm:bb}
For $\mu,\nu\in\cP_2(\R^d)$,
\begin{align*}
\W_2^2(\mu,\nu)=\inf\Big\{\int_0^1\!\!\int\abs{v_t}^2\,\dd\rho_t\,\dd t\ :\ \partial_t\rho_t+\divg(\rho_tv_t)=0,\ \rho_0=\mu,\ \rho_1=\nu\Big\},
\end{align*}
and minimisers have gradient velocities $v_t=\nabla\psi_t$ with straight particle paths.
\end{theorem}

\citet{benamou2000computational} proved this for bounded, compactly supported densities, using Lagrangian coordinates for smooth $(\rho,v)$; the statement on all of $\cP_2(\R^d)$, with distributional solutions of the continuity equation, belongs to the general theory of \citet{ambrosio2008gradient}.

Among all velocity fields generating a given curve, the one with least kinetic energy at each time is the gradient field $v_t=\nabla\psi_t$, obtained by projecting any admissible $v_t$ onto the closure of gradients in $L^2(\rho_t)$. This projection supplies the tangent spaces on which Otto's calculus is built. The same fact explains why action matching learns gradient fields (\cref{sec:simfree}), and the kinetic energy that gradient fields minimise reappears as a regulariser for continuous normalising flows (\cref{sec:cnf}).

\paragraph{Otto calculus and Wasserstein gradient flows.}

Otto's formal calculus \citep{otto2001geometry} treats $\cP_2(\R^d)$ as a Riemannian manifold whose tangent space at $\rho$ consists of gradient fields $\nabla\psi$ with inner product $\int\inner{\nabla\psi_1}{\nabla\psi_2}\dd\rho$; the induced distance is $\W_2$. The Wasserstein gradient of a functional $\cF$ at $\rho$ is the tangent vector $-\divg(\rho\nabla\frac{\delta\cF}{\delta\rho})$, represented by the velocity field $\nabla\frac{\delta\cF}{\delta\rho}$, where the first variation $\frac{\delta\cF}{\delta\rho}$ is the function, unique up to an additive constant, with $\frac{\dd}{\dd\varepsilon}\cF(\rho+\varepsilon\chi)\big|_{\varepsilon=0}=\int\frac{\delta\cF}{\delta\rho}\,\dd\chi$ for signed perturbations $\chi$ of total mass zero, so the gradient flow $\partial_t\rho_t=\divg(\rho_t\nabla\frac{\delta\cF}{\delta\rho})$ moves particles with velocity $-\nabla\frac{\delta\cF}{\delta\rho}$. For relative entropy $\cF(\rho)=\KL(\rho\,\|\,\pi)$ with $\pi\propto e^{-V}$ this flow is
\begin{equation}\label{eq:fp-gradflow}
\partial_t\rho_t=\divg\!\Big(\rho_t\nabla\log\frac{\rho_t}{\pi}\Big)=\Delta\rho_t+\divg(\rho_t\nabla V),
\end{equation}
the Fokker--Planck equation of the Langevin diffusion $\dd X_t=-\nabla V(X_t)\dd t+\sqrt2\,\dd W_t$. \citet{jordan1998variational} made this rigorous through the minimising-movement (JKO) scheme
\begin{equation}\label{eq:jko}
\rho_{k+1}\in\argmin_{\rho}\Big\{\KL(\rho\,\|\,\pi)+\frac{1}{2\tau}\W_2^2(\rho,\rho_k)\Big\},
\end{equation}
where $\tau>0$ is the time step and $\rho_k$ approximates the solution at time $k\tau$. As $\tau\to0$ the piecewise-constant interpolations converge to the solution of \eqref{eq:fp-gradflow}, weakly in $L^1$ at every time. The authors assume $V\ge0$ smooth with $\abs{\nabla V}\le C(V+1)$ and an initial density with finite second moment and finite free energy. They do not need $\pi$ to be normalisable, since $\KL(\rho\,\|\,\pi)$ can be replaced by the free energy $\int V\rho+\int\rho\log\rho$. The general theory of such flows is developed by \citet{ambrosio2008gradient}, and \citet{santambrogio2017euclidean} gives an expository overview.

Two consequences recur throughout the survey. One is that diffusion is transport by the score. The heat equation $\partial_t\rho=\Delta\rho$ is the continuity equation with velocity $v_t=-\nabla\log\rho_t$, so the process $X_0+\sqrt2\,W_t$, whose law solves $\partial_t\rho=\Delta\rho$, and the deterministic flow $\dot x=-\nabla\log\rho_t(x)$ have the same marginals. This is the origin of the probability-flow ODE. The other is that sampling is optimisation. Langevin Monte Carlo discretises the Wasserstein gradient flow of $\KL(\cdot\,\|\,\pi)$ \citep{wibisono2018sampling}, which is why functional inequalities such as the log-Sobolev and Poincar\'e inequalities govern the convergence of Langevin-type samplers \citep{vempala2019rapid,chewi2024log}. Diffusion models need no such assumption on the data law. Their guarantees combine Girsanov's theorem with the convergence of the forward Ornstein--Uhlenbeck process, which holds for essentially any data law (\cref{sec:theory}).

\paragraph{Fisher--Rao, Wasserstein--Fisher--Rao and Stein geometries.}

The Wasserstein geometry moves mass only by transport. Other geometries on densities allow mass to be created and destroyed, or restrict the admissible velocities, and each has a place in what follows. The Fisher--Rao metric on densities, $\inner{\dot\rho_1}{\dot\rho_2}_{\rm FR}=\int\frac{\dot\rho_1\dot\rho_2}{\rho}$, where the tangent vectors $\dot\rho_1,\dot\rho_2$ at $\rho$ are signed densities of total mass zero, moves mass by creation and annihilation rather than transport; it is the information-geometric metric on statistical manifolds \citep{amari2016information}. Its gradient flow of $\KL(\cdot\,\|\,\pi)$ is the birth--death equation $\partial_t\rho=-\rho\big(\log\frac{\rho}{\pi}-\E_\rho\log\frac{\rho}{\pi}\big)$, which converges to $\pi$ from any everywhere-positive initial density with finite relative entropy; added to the Fokker--Planck equation of Langevin dynamics, it gives, for initial densities close to $\pi$ in relative entropy and bounded below by a positive multiple of $\pi$, an asymptotic convergence rate independent of the potential barriers of $\pi$ \citep{lu2019accelerating}. Unlike the Wasserstein flow, it cannot be simulated by independent particles. The Wasserstein--Fisher--Rao, or Hellinger--Kantorovich, geometry combines both mechanisms \citep{chizat2018interpolating,liero2018optimal}, and particle methods based on it are reviewed in \cref{sec:gradflows}. Restricting Wasserstein gradients to a reproducing-kernel Hilbert space yields the Stein geometry underlying Stein variational gradient descent \citep{liu2017stein,duncan2023geometry}. Finally, the simplex with the Fisher--Rao metric is isometric to a portion of a sphere through $p\mapsto2\sqrt p$, a fact exploited by flow matching on statistical manifolds (\cref{sec:discrete}).

\subsection{Entropic transport and heat-flow maps}

Two further constructions connect this geometry with stochastic processes. Entropic regularisation turns optimal transport into a problem about Brownian motion, and the heat flow produces transport maps of its own.

\paragraph{Entropic optimal transport.}

Adding entropy to the Kantorovich problem,
\begin{equation}\label{eq:eot}
\min_{\pi\in\Pi(\mu,\nu)}\int c\,\dd\pi+\varepsilon\,\KL(\pi\,\|\,\mu\otimes\nu),
\end{equation}
gives a strictly convex problem whose solution has the form
\begin{align*}
\pi(\dd x,\dd y)=e^{f(x)+g(y)-c(x,y)/\varepsilon}\,\mu(\dd x)\,\nu(\dd y),
\end{align*}
computable by Sinkhorn's alternating scaling \citep{cuturi2013sinkhorn,peyre2019computational}. Here $c$ is the cost and $\varepsilon>0$ the strength of the regularisation, and the potentials $f$ and $g$, unrelated to the drift and diffusion coefficient of \cref{tab:notation}, are fixed by the two marginal constraints. For $c(x,y)=\abs{x-y}^2/2$, problem \eqref{eq:eot} is the static form of the Schr\"odinger bridge problem on $[0,1]$ with Brownian reference of variance $\varepsilon$ (\cref{sec:schrodinger-problem}; on $[0,T]$ the cost becomes $\abs{x-y}^2/(2T)$), and Sinkhorn's iteration is the static form of iterative proportional fitting, which diffusion Schr\"odinger bridges implement on path space (\cref{sec:sb}).

\paragraph{Transport maps from heat flow.}\label{sec:heatmaps}

Running the probability-flow ODE of an Ornstein--Uhlenbeck process backwards from its terminal law $p_T$, which approaches $\gamma_d$ as $T\to\infty$, produces a transport map onto $\mu_{\rm data}$; a deterministic diffusion sampler approximates it by starting from $\gamma_d$. Such maps were studied before diffusion models. \citet{kim2012generalization} showed that the map obtained from the reverse heat flow between $\gamma_d$ and $\nu=e^{-V}\gamma_d$ with $V$ convex is a contraction, recovering the conclusion of Caffarelli's theorem with a map built from the heat flow rather than from optimal transport; whether the two maps differ was left open there. They also extended the result to non-Gaussian sources $e^{-U}$, with $U$ convex, block-radially symmetric and satisfying a third-derivative condition, and targets $e^{-(U+V)}$ with $V$ convex and sharing the same symmetry. \citet{mikulincer2023brownian} constructed the Brownian transport map from Wiener space to a target through the F\"ollmer process (\cref{sec:follmer}) and proved Lipschitz bounds with applications to functional inequalities. It is natural to ask whether the heat-flow map is itself the Brenier map. In general it is not. \citet{khrulkov2022understanding} conjectured, in a preprint first circulated in 2022, that the encoder map of the variance-preserving diffusion, the long-time limit of the probability-flow ODE of the Ornstein--Uhlenbeck Fokker--Planck equation, is the optimal transport map from $\mu_{\rm data}$ to $\gamma_d$, and proved this for Gaussian data. \citet{lavenant2022flow} then showed that it fails in general, strengthening an earlier counterexample of \citet{tanana2021comparison}. The probability flow is therefore a canonical but non-optimal transport, a fact behind the design of the straighter, transport-inspired paths of \cref{sec:simfree}.

\section{Stochastic analysis on path space}\label{sec:pathspace}

Path space carries the third geometry. Diffusion models, bridges and control problems concern not a single law but the law of a whole trajectory. The natural operations on them, reversing time, conditioning on the future, steering towards a target, are changes of measure on the space of paths. This section assembles the tools in the order in which they will be used, beginning with stochastic differential equations and their path measures, then time reversal and conditioning, and finally three constructions that drive a process into a prescribed law.

\subsection{Diffusions and their path measures}

We fix notation for diffusions and record the two facts used most often later, that the marginals of a diffusion are carried by a deterministic velocity field and that the relative entropy between two diffusions with the same diffusion coefficient is, apart from the relative entropy of the initial laws, an integral of the squared difference of their drifts.

\paragraph{SDEs, generators and Kolmogorov equations.}

Consider the It\^o SDE on $\R^d$
\begin{equation}\label{eq:sde}
\dd X_t=f(t,X_t)\,\dd t+\sigma(t,X_t)\,\dd W_t,\qquad X_0\sim\rho_0,
\end{equation}
with $W$ an $m$-dimensional Brownian motion and $a=\sigma\sigma^\top$. Its generator acts on test functions by $\cL_t\varphi=\inner{f}{\nabla\varphi}+\frac12\tr(a\nabla^2\varphi)$, and the marginal densities solve the Fokker--Planck, or Kolmogorov forward, equation
\begin{equation}\label{eq:fpe}
\partial_tp_t=\cL_t^*p_t=-\divg(fp_t)+\tfrac12\textstyle\sum_{ij}\partial_{ij}(a_{ij}p_t).
\end{equation}
Writing the right-hand side as $-\divg\big(p_t(f-\frac12\divg a-\frac12 a\nabla\log p_t)\big)$, where $(\divg a)_i=\sum_j\partial_ja_{ij}$, shows that \eqref{eq:fpe} is a continuity equation with velocity
\begin{equation}\label{eq:pf-velocity}
v_t(x)=f(t,x)-\tfrac12\divg a(t,x)-\tfrac12a(t,x)\nabla\log p_t(x).
\end{equation}
This is the probability-flow velocity. When $p_t$ is positive and differentiable and $v_t$ is regular enough, for instance locally Lipschitz with the integrability of \cref{sec:transport}, the ODE $\dot x=v_t(x)$ has the same one-time marginals as \eqref{eq:sde}; for a degenerate initial law, such as a point mass, this holds only from positive times on. Standard references for this section are \citet{karatzas1991brownian} and \citet{pavliotis2014stochastic}, with \citet{kunita1990stochastic} for stochastic flows and the variation processes used in \cref{sec:ibp}; for SDEs on manifolds see \citet{elworthy1982stochastic} and \citet{hsu2002stochastic}.

\paragraph{Path measures, Girsanov's theorem and relative entropy.}

Let $\P$ and $\Q$ be the laws on $C([0,T];\R^d)$ of two diffusions with the same invertible diffusion coefficient $\sigma$ and drifts $b^{\P}$, $b^{\Q}$, and write $X$ for the coordinate process on path space and $\P_t$, $\Q_t$ for the time-$t$ marginals. Under conditions guaranteeing that the stochastic exponential is a martingale, or after a localisation argument when $\Q$ is the unique law of its SDE, Girsanov's theorem gives, with $\P\ll\Q$ whenever the right-hand side is finite,
\begin{equation}\label{eq:girsanov-kl}
\KL(\P\,\|\,\Q)=\KL(\P_0\,\|\,\Q_0)+\frac12\,\E_{\P}\!\int_0^T\big\|\sigma^{-1}(b^{\P}-b^{\Q})(t,X_t)\big\|^2\,\dd t.
\end{equation}
Together with the data-processing inequality $\KL(\P_T\,\|\,\Q_T)\le\KL(\P\,\|\,\Q)$, this identity drives the convergence theory of diffusion models (\cref{sec:theory}), the evidence lower bound in continuous time (\cref{sec:likelihood}) and the losses for diffusion samplers (\cref{sec:control}). The other basic tool is the chain rule $\KL(\P\,\|\,\Q)=\KL(\P_{0}\,\|\,\Q_{0})+\E_{\P_0}\KL(\P^{x}\,\|\,\Q^{x})$, where $\P^x$ conditions on $X_0=x$.

\subsection{Time reversal and conditioning}

Generation runs a noising process backwards, and conditioning reweights it by what is to come. Both operations turn one diffusion into another, and in both cases the new drift has an explicit form.

\paragraph{Time reversal.}\label{sec:time-reversal}

\begin{theorem}[Time reversal \citep{anderson1982reverse,haussmann1986time}]\label{thm:reversal}
Let $X$ solve \eqref{eq:sde} on $[0,T]$ with marginal densities $p_t$, and suppose the coefficients and $p_t$ are regular enough, for example under the conditions of \citet{haussmann1986time}. Then $Y_s=X_{T-s}$ is a diffusion solving
\begin{equation}\label{eq:reverse}
\dd Y_s=\Big[-f+\divg a+a\nabla\log p\Big](T-s,Y_s)\,\dd s+\sigma(T-s,Y_s)\,\dd\overline{W}_s
\end{equation}
for a Brownian motion $\overline W$ in the reversed filtration.
\end{theorem}

For additive noise $\sigma(t,x)=g(t)\Id$ the reverse drift is $-f+g^2\nabla\log p_t$, the form used by score-based diffusion models \citep{song2021score}. The theorem has been proved in several ways, among them through Kolmogorov equations, through relative entropy and through Malliavin calculus (another route uses enlargement of filtrations), and the last two recur in this survey. Relative entropy on path space drives the convergence theory of \cref{sec:theory}, and integration by parts on Wiener space is the subject of \cref{sec:ibp}. \citet{anderson1982reverse} argued through Kolmogorov equations, assuming that the forward equations for $X$ and for the pair $(X,\overline W)$ have unique smooth solutions, conditions he notes are hard to translate into requirements on the coefficients. \citet{haussmann1986time} gave a rigorous proof using weak ($H^1$-valued) solutions of the forward and backward Kolmogorov equations for $X$. They assume Lipschitz coefficients of linear growth and an integrability condition on the density; the condition holds, for instance, when $X_0$ has a density in a weighted $L^2$ space and either the diffusion is uniformly elliptic or $\sum_{ij}\partial_{ij}a_{ij}$ is bounded, and it holds for every initial law when $\partial_t+\cL_t$ is hypoelliptic.

\citet{follmer1986time} took the entropic route, using Nelson's stochastic derivatives on Wiener space. The path law must have finite entropy relative to Wiener measure, so that the diffusion coefficient is the identity, while the drift need only have finite energy and may depend on the whole past. Under that condition the forward drift and the drift of the reversed process add up, given the present state, to $\nabla\log p_t$ for almost every $t$, with no regularity assumed on the drift. With a localised entropy condition the argument extends to infinite-dimensional diffusions with countably many interacting coordinates, where the duality equation involves the conditional density of each coordinate given the others and can fail without bounds on the interaction \citep{follmer1986infinite}. \citet{cattiaux2023time} carried it from Wiener measure to entropy relative to a general reference diffusion whose time reversal is known, for instance a reversible one, allowing singular drifts, through an integration-by-parts formula for the carr\'e du champ. The Malliavin route is due to \citet{millet1989integration}, who obtained a condition that is necessary as well as sufficient. For Lipschitz coefficients of linear growth and $X_t$ with a density for every $t>0$, the reversed process is a diffusion exactly when the distributional divergence of $a\,p_t$ is locally integrable on $(0,T]\times\R^d$, a condition implied by those of \citet{haussmann1986time}. Nelson's stochastic mechanics had already identified half the difference between the forward and backward drifts, which he called the osmotic velocity, as $\frac12 a\nabla\log p_t$ for constant $a$ \citep{nelson1967dynamical}. Up to the factor $\frac12a$, this is the score, half a century before diffusion models.

\paragraph{Doob $h$-transforms and bridges.}\label{sec:doob}

Let $h(t,x)>0$ be space--time harmonic for $X$, that is, $\partial_th+\cL_th=0$. Reweighting the path measure by $h(T,X_T)/h(0,X_0)$ produces a new Markov process, with the same initial law $\rho_0$, whose generator is $\cL_t^h\varphi=\cL_t\varphi+\inner{a\nabla\log h}{\nabla\varphi}$; the drift acquires the term $a\nabla\log h$. This is the space--time form of the transform that \citet{doob1957conditional} introduced for Brownian motion on a Green space, with transition density $p(t,\xi,\eta)h(\eta)/h(\xi)$, where $p$ is the transition density of the Brownian motion and $h$ is strictly positive and superharmonic. Reweighting instead by $h(T,X_T)/\int h(0,\cdot)\,\dd\rho_0$ gives the same dynamics started from the tilted law $\propto h(0,\cdot)\rho_0$. Two choices of $h$ are used repeatedly below. The first is $h(t,x)=\E[G(X_T)\mid X_t=x]$ for a nonnegative terminal weight $G$. Started from the tilted law $\propto h(0,\cdot)\rho_0$, the transformed process has terminal law proportional to $G\,\rho_T$; started from a fixed point $x_0$, the case $\rho_0=\delta_{x_0}$, its terminal law is proportional to $G\,p_{T\mid0}(\cdot\mid x_0)$, where $p_{T\mid0}$ is the transition density of $X$. Classifier guidance (\cref{sec:conditioning}) is this construction applied to the reverse-time generative process, whose terminal state is the clean sample, with $G(x)=p(y\mid x)$ for a label $y$. The second is $h(t,x)=p_{T\mid t}(z\mid x)$, the transition density to a fixed point $z$, for which the transformed process is the bridge from $X_0$ to $z$; for Brownian motion the bridge drift is $(z-x)/(T-t)$, the conditional velocity of bridge matching (\cref{sec:sb}). Bridges can be simulated with guided proposals \citep{delyon2006simulation,schauer2017guided} and with learned scores \citep{heng2025simulating}.

\subsection{Control, localisation and the Schr\"odinger problem}

Three constructions on path space produce a process whose terminal law is prescribed. Optimal control does so by steering at a quadratic cost, stochastic localisation by gradually revealing a sample, and the Schr\"odinger problem by fixing both endpoints and asking for the least costly path measure between them.

\paragraph{Stochastic optimal control and the F\"ollmer process.}\label{sec:follmer}

Consider controlled dynamics with law $\P^u$, whose uncontrolled ($u=0$) generator we again write $\cL_t$, and a terminal cost $\Phi$,
\begin{align*}
\dd X_t^u&=\big(b(t,X^u_t)+\sigma u_t\big)\dd t+\sigma\,\dd W_t,\\
J(u)&=\E\Big[\int_0^T\tfrac12\norm{u_t}^2\dd t+\Phi(X^u_T)\Big],
\end{align*}
and let $V(t,x)=\inf_u\E\big[\int_t^T\tfrac12\norm{u_s}^2\dd s+\Phi(X^u_T)\,\big|\,X^u_t=x\big]$ be the least expected cost from state $x$ at time $t$. The value function $V(t,x)$ solves the Hamilton--Jacobi--Bellman equation
\begin{align*}
\partial_tV+\cL_tV-\tfrac12\norm{\sigma^\top\nabla V}^2=0,\qquad V(T,\cdot)=\Phi,
\end{align*}
and the optimal control is $u^*=-\sigma^\top\nabla V$ \citep{fleming2006controlled}. The Hopf--Cole substitution $V=-\log h$ linearises this equation into the backward Kolmogorov equation for $h$, whose solution, by the Feynman--Kac formula, is $h(t,x)=\E[e^{-\Phi(X_T)}\mid X_t=x]$ for the uncontrolled process. Optimal control with quadratic cost is therefore a Doob $h$-transform, and by Girsanov's theorem, for a fixed initial point $X_0=x_0$,
\begin{equation}\label{eq:gibbs}
\min_u J(u)=-\log\E\big[e^{-\Phi(X_T)}\big],\qquad \frac{\dd\P^{u^*}}{\dd\P^{0}}=\frac{e^{-\Phi(X_T)}}{\E e^{-\Phi(X_T)}},
\end{equation}
which is the variational formula of \citet{boue1998variational}, the Gibbs variational principle on Wiener space. They prove it for any bounded measurable functional of the Brownian path, with the infimum over all progressively measurable controls, extend it to functionals bounded above, and transfer it to diffusions with unique strong solutions. A terminal cost that is unbounded above lies outside their theorems. Path-integral control \citep{kappen2005path} and neural solvers built on path-space losses \citep{nusken2021solving} exploit the same structure. A distinguished special case is the F\"ollmer process, which is Brownian motion started at the origin, run on $[0,1]$, and steered to a target $\nu$ with density $f=\dd\nu/\dd\gamma_d$ (here not the drift) by the drift $\nabla\log Q_{1-t}f(x)$, where $Q_sf(x)=\E f(x+W_s)$ is the heat semigroup of Brownian motion, generated by $\frac12\Delta$. It minimises relative entropy with respect to Wiener measure among processes with terminal law $\nu$, and this yields Lehec's representation of entropy and of functional inequalities \citep{lehec2013representation,follmer1988random}. The F\"ollmer drift is the prototype of the diffusion samplers of \cref{sec:control}.

\paragraph{Stochastic localisation.}\label{sec:localisation}

Stochastic localisation, introduced by \citet{eldan2013thin} and abstracted into a general framework by \citet{chen2022localization}, constructs a measure-valued martingale $(\nu_t)$ that starts at $\nu_0=\mu$ and converges to a point mass. In its simplest Gaussian form one observes $Y_t=tX+W_t$ with $X\sim\mu$ and sets $\nu_t=\law(X\mid Y_{[0,t]})$. In its own filtration $(Y_t)$ solves $\dd Y_t=\E[X\mid Y_t]\,\dd t+\dd B_t$ for a Brownian motion $B$, and $Y_t/t\to X$ as $t\to\infty$, so sampling $X$ amounts to integrating this SDE with a learned posterior-mean, or denoiser, drift. \citet{montanari2023sampling} observed that, after a time change and a time-dependent rescaling, this is the reverse process of a diffusion model, with $t$ playing the role of signal-to-noise ratio; \citet{elalaoui2022sampling} had used the same idea to design algorithmic samplers for spin glasses (\cref{sec:statphys}). Stochastic localisation also supplies convergence bounds for diffusion models that are linear in the dimension, up to logarithmic factors, under only a finite second moment (\cref{sec:theory}).

\paragraph{The Schr\"odinger problem.}\label{sec:schrodinger-problem}

In 1931 and 1932 Schr\"odinger asked for the most likely way in which a large cloud of independent Brownian particles, observed at two times, could have reached an unexpected final configuration \citep{schrodinger1931umkehrung,schrodinger1932theorie}. The first paper is a study of the time reversal of diffusions, and in the second he called the analogy between the solution and wave mechanics so striking that he found it hard to believe accidental. Through large deviations his question becomes a minimisation of relative entropy on path space \citep{follmer1988random,leonard2014survey}. Given a reference path measure $\Rref$, say reversible Brownian motion with variance $\varepsilon$, whose initial and reversing measure is Lebesgue measure (an unbounded reference, which relative entropy on path space allows), and marginals $\mu,\nu$, the dynamic Schr\"odinger problem is
\begin{equation}\label{eq:sb}
\min\big\{\KL(\P\,\|\,\Rref):\ \P_0=\mu,\ \P_T=\nu\big\}.
\end{equation}
Write $\P_{0T}$ and $\Rref_{0T}$ for the joint laws of $(X_0,X_T)$ and $\P^{xy}$, $\Rref^{xy}$ for the bridges from $x$ to $y$. By the chain rule for relative entropy, $\KL(\P\,\|\,\Rref)=\KL(\P_{0T}\,\|\,\Rref_{0T})+\E_{\P_{0T}}\KL(\P^{xy}\,\|\,\Rref^{xy})$, so the optimiser is a mixture of reference bridges $\Rref^{xy}$ over the static optimiser $\pi^*$, which solves \eqref{eq:eot} with cost $-\varepsilon\log r_{0T}(x,y)$, where $r_{0T}$ is the Lebesgue density of $\Rref_{0T}$. Under mild conditions on $\Rref$ and on the marginals, the solution has the product form $\dd\P^*/\dd\Rref=\mathfrak f(X_0)\mathfrak g(X_T)$, where $(\mathfrak f,\mathfrak g)$ solve the Schr\"odinger system, the pair of equations expressing that $\P^*$ has marginals $\mu$ and $\nu$. Without such conditions the density need not factorise, for instance when its support is not a product set. The solution is Markov, it is reciprocal (given its values at two times, its path between them is independent of its path outside them), and it is an $h$-transform of $\Rref$ \citep{follmer1988random,leonard2014survey}. As $\varepsilon\to0$ the bridge converges to the displacement interpolation of optimal transport \citep{mikami2004monge,chen2016relation}, and the problem can equivalently be posed as stochastic control with fixed endpoints \citep{daipra1991stochastic,chen2016relation,chen2021stochastic}. Generative modelling instantiates \eqref{eq:sb} in \cref{sec:sb}.

\section{Integration by parts: the source of scores}\label{sec:ibp}

This section develops the main analytic tool of the survey. A score $\nabla\log p$ of an intractable density is not computed directly. It is either written as a conditional expectation of a computable random variable, which a network then regresses onto, or removed from the training objective altogether, as in implicit score matching. Each such representation is an integration by parts (IBP), and the setting of the IBP, whether finite-dimensional Gaussian noise, the data law or Wiener space, determines which forward processes it can handle.

\subsection{Gaussian integration by parts and score matching}

The simplest setting, and the one on which diffusion models were first built, is integration by parts against a finite-dimensional Gaussian variable.

\paragraph{Stein and Tweedie.}

Stein's lemma states that for $Z\sim\cN(0,\Id)$ and $\varphi$ almost differentiable, for instance Lipschitz, with $\E\abs{\nabla\varphi(Z)}<\infty$, $\E[Z\varphi(Z)]=\E[\nabla\varphi(Z)]$ \citep{stein1981estimation}. Applied to a Gaussian corruption, it yields the denoiser--score identity. \citet{robbins1956empirical} published it, crediting the example to Tweedie, for exponential-family kernels $p(x\mid\lambda)=e^{\lambda x}f(x)h(\lambda)$, in which $\lambda$ is drawn from a prior, $f$ is a fixed function and $h(\lambda)$ normalises the kernel; the posterior mean of $\lambda$ given $x$ is then $\frac{\dd}{\dd x}\log\big(p(x)/f(x)\big)$, with $p$ the marginal density of $x$. For Gaussian noise, \citet{miyasawa1961empirical} wrote the least-squares estimate of a normal mean through the gradient of the log-density of the observations, which is the form $X+\nabla\log p(X)$, with $p$ the marginal density, that \citet[\S3]{stein1981estimation} gives as the formal Bayes estimate. \citet{efron2011tweedie} revisits it as an empirical Bayes correction for selection bias.

\begin{proposition}[Tweedie's formula]\label{prop:tweedie}
Let $X_t=\alpha_tX_0+\sigma_t\epsilon$ with $\epsilon\sim\cN(0,\Id)$ independent of $X_0\sim\mu$, and let $p_t$ be the density of $X_t$. Then
\begin{equation}\label{eq:tweedie}
\nabla\log p_t(x)=\frac{\alpha_t\E[X_0\mid X_t=x]-x}{\sigma_t^2}=-\frac{1}{\sigma_t}\,\E[\epsilon\mid X_t=x].
\end{equation}
\end{proposition}
\begin{proof}
Write $p_t(x)=\int\phi_{\sigma_t}(x-\alpha_ty)\,\mu(\dd y)$ with $\phi_\sigma$ the $\cN(0,\sigma^2\Id)$ density, differentiate under the integral using $\nabla\phi_\sigma(z)=-z\phi_\sigma(z)/\sigma^2$, and divide by $p_t(x)$; the ratio is the posterior expectation.
\end{proof}

Tweedie's formula identifies the score with the minimum mean-squared-error denoiser, which is why ``training a denoiser'' and ``learning a score'' are the same problem \citep{raphan2011least}, and why denoisers can serve as priors for inverse problems \citep{kadkhodaie2021solving}.

\paragraph{Score matching.}

The Fisher divergence $\frac12\E_p\norm{s_\theta-\nabla\log p}^2$ cannot be evaluated directly, but integration by parts against $p$ itself removes the unknown score.

\begin{proposition}[Implicit and denoising score matching]\label{prop:sm}
(i) If $p$ and $s_\theta$ are differentiable, $s_\theta$ and $\nabla\log p$ are square-integrable under $p$, the derivatives $\partial_is_{\theta,i}$ are integrable under $p$, and $p\,s_\theta\to0$ at infinity, then $\frac12\E_p\norm{s_\theta-\nabla\log p}^2=\E_p\big[\frac12\norm{s_\theta}^2+\tr\nabla s_\theta\big]+C$ \citep{hyvarinen2005estimation}. (ii) For $X_t=\alpha_tX_0+\sigma_t\epsilon$,
\begin{align*}
\E\big\|s_\theta(X_t)-\nabla\log p_t(X_t)\big\|^2=\E\Big\|s_\theta(X_t)+\frac{\epsilon}{\sigma_t}\Big\|^2+C',
\end{align*}
with $C,C'$ independent of $\theta$ \citep{vincent2011connection}.
\end{proposition}
Part (ii) follows from \cref{prop:tweedie}, since the regression target $-\epsilon/\sigma_t$ has conditional mean $\nabla\log p_t(X_t)$, and the population minimiser of a squared loss is the conditional mean. Part (i) is the trace-form estimator, which sliced score matching approximates with random projections \citep{song2019sliced}. \citet{lyu2009interpretation} interprets score matching as minimising the rate of change of relative entropy under Gaussian perturbation. Every later regression objective, from flow matching to bridge matching and the interpolant losses, follows the same \emph{conditional-to-marginal} pattern, in which an intractable target is replaced by a tractable random variable with the right conditional mean.

\subsection{Malliavin calculus in a nutshell}\label{sec:malliavin}

When the noise passes through a nonlinear SDE, the corrupted variable is no longer a Gaussian perturbation of the data, and integration by parts in finitely many Gaussian variables no longer applies. Its replacement is integration by parts on Wiener space. We recall what is needed, with short proofs; \citet{nualart2006malliavin} and \citet{nualart2018introduction} are the standard references, with \citet{malliavin2006stochastic} for the covering-field language. Let $W$ be an $m$-dimensional Brownian motion on $[0,T]$, defined on a probability space $\Omega$ with the filtration $(\cF_t)$ it generates, and $H=L^2([0,T];\R^m)$. The Wiener integrals $W(h)=\int_0^Th(t)\cdot\dd W_t$ form an isonormal Gaussian process on $H$, that is, a centred Gaussian family $\{W(h):h\in H\}$ with $\E[W(h)W(g)]=\inner hg_H$; isonormal processes are defined in the same way on any real separable Hilbert space. For a smooth cylindrical random variable $F=\varphi\big(W(h_1),\dots,W(h_n)\big)$, with $\varphi$ smooth and of polynomial growth together with its derivatives, the Malliavin derivative is the $H$-valued random variable $\D_tF=\sum_i\partial_i\varphi(\cdots)h_i(t)$.

\begin{lemma}[Integration by parts on Wiener space]\label{lem:ibp-wiener}
For smooth $F$ and $h\in H$, $\E\big[\inner{\D F}{h}_H\big]=\E\big[F\,W(h)\big]$.
\end{lemma}
\begin{proof}
By linearity we may assume $\norm h_H=1$ and, after a Gram--Schmidt change of variables in $\varphi$, that $F=\varphi(Z)$ with $Z=(W(e_1),\dots,W(e_n))$, $e_1=h$ and $(e_i)$ orthonormal. Then $Z\sim\cN(0,\Id_n)$ and $\inner{\D F}{h}_H=\partial_1\varphi(Z)$, so the claim is Stein's identity $\E[\partial_1\varphi(Z)]=\E[Z_1\varphi(Z)]$.
\end{proof}

Applying the lemma to the product $FG$ gives $\E[G\inner{\D F}{h}_H]=\E[FGW(h)]-\E[F\inner{\D G}{h}_H]$, from which closability of $\D$ in $L^2(\Omega)$ follows. Its closure has domain $\mathbb{D}^{1,2}$, the completion of smooth variables under $\norm F_{1,2}^2=\E F^2+\E\norm{\D F}_H^2$. For $\varphi\in C^1_b(\R^d)$ and $F\in(\mathbb{D}^{1,2})^d$ the chain rule $\D\varphi(F)=\sum_i\partial_i\varphi(F)\,\D F^i$ holds. The adjoint of $\D$ is the \emph{Skorokhod integral} or divergence $\Sk$. An $H$-valued $u\in L^2(\Omega;H)$ is in $\Dom\Sk$ if there is a constant $c$ such that $\abs{\E\inner{\D F}{u}_H}\le c\norm{F}_{L^2}$ for all $F\in\mathbb{D}^{1,2}$, and then
\begin{equation}\label{eq:duality}
\E\big[\inner{\D F}{u}_H\big]=\E\big[F\,\Sk(u)\big]\qquad\text{for all }F\in\mathbb{D}^{1,2}.
\end{equation}
For an elementary field $u=\sum_jF_jh_j$ with smooth $F_j$, \cref{lem:ibp-wiener} gives
\begin{equation}\label{eq:delta-elementary}
\Sk(u)=\sum_jF_jW(h_j)-\sum_j\inner{\D F_j}{h_j}_H,
\end{equation}
and more generally $\Sk(Fu)=F\Sk(u)-\inner{\D F}{u}_H$ whenever both sides make sense. If $u$ is adapted to the Brownian filtration and square integrable, $\Sk(u)$ coincides with the It\^o integral $\int_0^Tu_t\cdot\dd W_t$ \citep{nualart2006malliavin}. The correction term in \eqref{eq:delta-elementary} is what distinguishes anticipating from adapted integrands. The Clark--Ocone formula $F=\E F+\int_0^T\E[\D_tF\mid\cF_t]\cdot\dd W_t$ represents a Malliavin-differentiable functional as a stochastic integral, identifying the integrand in It\^o's martingale representation \citep{ocone1984malliavin}. The \emph{Malliavin covariance} of an $\R^d$-valued $F$ is the random matrix $\gammaM_F=\big(\inner{\D F^i}{\D F^j}_H\big)_{ij}$.

To apply the calculus to diffusions one needs the Malliavin derivative of an SDE solution, which is expressed through the first variation process. Let $X$ solve \eqref{eq:sde} from a fixed point with smooth coefficients having bounded derivatives, and let $Y_t=\partial X_t/\partial X_0$ solve the linearised equation $\dd Y_t=\nabla f(t,X_t)Y_t\dd t+\sum_l\nabla\sigma_l(t,X_t)Y_t\dd W^l_t$, $Y_0=\Id$, where $\sigma_l$ is the $l$-th column of $\sigma$ and $W^l$ the $l$-th coordinate of $W$.

\begin{proposition}[Malliavin derivative of an SDE solution]\label{prop:sde-derivative}
For $s\le t$, $\D_sX_t=Y_tY_s^{-1}\sigma(s,X_s)$, and $\D_sX_t=0$ for $s>t$.
\end{proposition}
\begin{proof}[Sketch]
Differentiating the integral equation for $X_t$ in the direction of a perturbation of the noise at time $s$ gives, for $t\ge s$, the linear SDE $\D_sX_t=\sigma(s,X_s)+\int_s^t\nabla f\,\D_sX_r\,\dd r+\sum_l\int_s^t\nabla\sigma_l\,\D_sX_r\,\dd W^l_r$. Its solution is the fundamental matrix $Y_tY_s^{-1}$ of the linearised equation applied to the initial value $\sigma(s,X_s)$. For $s>t$ the variable $X_t$ does not depend on the noise after time $t$.
\end{proof}

By \cref{prop:sde-derivative}, the Malliavin covariance of $X_T$ is
\begin{equation}\label{eq:malliavin-cov}
\gammaM_{X_T}=\int_0^TY_TY_s^{-1}a(s,X_s)(Y_s^{-1})^\top Y_T^\top\,\dd s
\end{equation}
\citep{nualart2006malliavin}. Malliavin's original application was a probabilistic proof of H\"ormander's hypoellipticity theorem \citep{hormander1967hypoelliptic,malliavin1978stochastic,bismut1981martingales}. When the diffusion vector fields $\sigma_1,\dots,\sigma_m$ (the columns of $\sigma$), together with their iterated brackets with all the vector fields including the Stratonovich drift $f-\frac12\sum_k(\nabla\sigma_k)\sigma_k$, span the tangent space at every point (the parabolic form of H\"ormander's condition, in which the drift does not count on its own), $\gammaM_{X_T}$ is invertible with inverse moments of all orders, and $X_T$ has a smooth density even if $a$ is degenerate. \citet{hairer2011malliavin} states this for time-homogeneous smooth coefficients with bounded derivatives. For additive noise the Stratonovich and It\^o drifts coincide. Kinetic diffusion models, whose noise acts only on velocities, are of this kind (\cref{sec:kinetic}). Under H\"ormander's original condition, in which the drift counts as a generator, $X_T$ need not have a density; \citet[\S5]{bismut1981martingales} shows that the time-integrated law $\int_0^TP_t(x,\cdot)\,\dd t$, where $P_t(x,\cdot)$ is the law of $X_t$ started at $x$, then has one.

\subsection{Densities and scores by integration by parts}

With these tools a score can be written as a conditional expectation for any sufficiently regular, nondegenerate functional of Brownian motion. Suppose $F\in(\mathbb{D}^\infty)^d$, where $\mathbb{D}^\infty=\bigcap_{k,p}\mathbb{D}^{k,p}$ consists of variables with Malliavin derivatives of all orders in every $L^p$, and that $F$ is nondegenerate, meaning $\gammaM_F$ is invertible with $\E\abs{\det\gammaM_F}^{-p}<\infty$ for all $p$. The chain rule gives $\D\varphi(F)=\sum_j\partial_j\varphi(F)\D F^j$, hence $\partial_i\varphi(F)=\inner{\D\varphi(F)}{u_i}_H$ with
\begin{equation}\label{eq:covering}
u_i=\sum_{k=1}^d(\gammaM_F^{-1})_{ik}\,\D F^k,\qquad\text{so that}\qquad \inner{\D F^j}{u_i}_H=\delta_{ij}.
\end{equation}
Any $H$-valued $u_i$ with the property on the right is called a \emph{covering vector field} \citep{malliavin2006stochastic}. Duality \eqref{eq:duality} then gives $\E[\partial_i\varphi(F)]=\E[\varphi(F)\Sk(u_i)]$, and comparing with $\E[\partial_i\varphi(F)]=-\int\varphi\,\partial_ip_F$, where $p_F$ is the density of $F$, yields the following representation.

\begin{theorem}[Score as a conditional Skorokhod integral]\label{thm:malliavin-score}
If $F\in(\mathbb{D}^\infty)^d$ is nondegenerate and $u_i\in\Dom\Sk$ is any covering vector field for the $i$-th direction, then $F$ has a smooth density $p_F$ and, for almost every $y$ with $p_F(y)>0$,
\begin{equation}\label{eq:malliavin-score}
\partial_i\log p_F(y)=-\E\big[\Sk(u_i)\,\big|\,F=y\big].
\end{equation}
\end{theorem}

The derivative and the divergence combine into an integration by parts for the law of a random vector, which underlies the score formulas of this section.

\begin{proposition}[Integration by parts for laws]\label{prop:ibp-laws}
Let $F\in(\mathbb{D}^{1,2})^d$ have an almost surely invertible Malliavin covariance $\gammaM_F$, and suppose the covering fields $u_i=\sum_k(\gammaM_F^{-1})_{ik}\D F^k$ belong to $\Dom\Sk$. Then for every $\varphi\in C^1_b(\R^d)$,
\begin{align*}
\E\big[\partial_i\varphi(F)\big]=\E\big[\varphi(F)\,\Sk(u_i)\big].
\end{align*}
Consequently, if $F$ has a $C^1$ positive density $p_F$, then $\partial_i\log p_F(y)=-\E[\Sk(u_i)\mid F=y]$ for almost every $y$, and the same identity holds for any other $u_i'\in\Dom\Sk$ with $\inner{\D F^j}{u'_i}_H=\delta_{ij}$.
\end{proposition}
\begin{proof}
By the chain rule, $\inner{\D\varphi(F)}{\D F^k}_H=\sum_j\partial_j\varphi(F)(\gammaM_F)_{jk}$; multiplying by $(\gammaM_F^{-1})_{ki}$ and summing over $k$ gives $\partial_i\varphi(F)=\inner{\D\varphi(F)}{u_i}_H$, and duality gives the first identity. Writing both sides as integrals against the law of $F$, $\int\partial_i\varphi\,p_F=-\int\varphi\,\partial_ip_F$ and $\E[\varphi(F)\Sk(u_i)]=\int\varphi(y)\,\E[\Sk(u_i)\mid F=y]\,p_F(y)\,\dd y$; since $\varphi$ is arbitrary, $-\partial_ip_F=p_F\,\E[\Sk(u_i)\mid F=\cdot]$. The same computation applies to any $u'_i$ with the covering property, because only $\partial_i\varphi(F)=\inner{\D\varphi(F)}{u'_i}_H$ was used.
\end{proof}

In one dimension the same argument with $\varphi=\one_{(y,\infty)}$ (after approximation) gives the density formula $p_F(y)=\E\big[\one_{\{F>y\}}\,\Sk(\D F/\norm{\D F}_H^2)\big]$ \citep{nualart2006malliavin}.

The right-hand side of \eqref{eq:malliavin-score} does not depend on which covering field is used, because every covering field reproduces the same distributional identity; different choices give different estimators of the same score, with different variances. For SDE solutions this representation was developed into explicit score formulas by \citet{mirafzali2025malliavin} for state-independent noise, and by \citet{mirafzali2025approach} for nonlinear SDEs with state-dependent diffusion. In the latter, the covering field \eqref{eq:covering} with $F=X_T$ is shown to be the unique field of minimal energy $\E\int_0^T\norm{u_t}^2\dd t$; evaluating the field $u_t(x)=x^\top Y_t^{-1}\sigma(t,X_t)$ at the anticipating point $x=Y_T^\top\gammaM_{X_T}^{-1}e_k$, with $e_k$ the $k$-th standard basis vector, gives the covering field $u_k$ for the $k$-th direction. Nualart's substitution formula for Skorokhod integrals of random fields evaluated at random points \citep{nualart2006malliavin} then expresses $\Sk(u_k)$ as an It\^o integral minus a correction involving the first and second variation processes $Y$ and $Y^{(2)}=\partial^2X/\partial X_0^2$. The final formula contains no abstract Malliavin derivatives and can be simulated alongside the SDE.

For a random field $u_t(x)$ adapted for each fixed $x$ and a random point $F$, the substitution formula reads
\begin{equation}\label{eq:substitution}
\Sk\big(u_\cdot(F)\big)=\int_0^Tu_t(x)\cdot\dd W_t\Big|_{x=F}-\sum_j\int_0^T\partial_ju_t(F)\cdot\D_tF^j\,\dd t
\end{equation}
under integrability conditions \citep{nualart2006malliavin}. It goes back to \citet[Prop.~4.12]{nualart1988stochastic}, who proved it within their anticipating calculus for adapted random fields that are continuously differentiable in $x$.

\begin{example}[The linear case]\label{ex:linear}
For $\dd X_t=A X_t\,\dd t+\Sigma\,\dd W_t$ started at $x$, the variation process $Y_t=e^{tA}$ is deterministic, so $\D_sX_T=e^{(T-s)A}\Sigma$, $\gammaM_{X_T}=C_T:=\int_0^Te^{(T-s)A}\Sigma\Sigma^\top e^{(T-s)A^\top}\dd s$ and the covering field $u_i(s)=e_i^\top C_T^{-1}e^{(T-s)A}\Sigma$ is deterministic. Its Skorokhod integral is a Wiener integral, $\Sk(u_i)=e_i^\top C_T^{-1}(X_T-e^{TA}x)$, and \eqref{eq:malliavin-score} gives $\nabla\log p_T(y)=-C_T^{-1}(y-e^{TA}x)$, the score of $\cN(e^{TA}x,C_T)$. For the Ornstein--Uhlenbeck process with random $X_0$, conditioning additionally on $X_0$ and averaging reproduces Tweedie's formula \eqref{eq:tweedie} and the denoising-score-matching target.
\end{example}

\paragraph{Bismut--Elworthy--Li formulas.}

A different integration by parts differentiates with respect to the initial point rather than the terminal point. Let $X^x$ be an elliptic diffusion started at $x$, with semigroup $P_t$, first variation $Y$ and square invertible $\sigma$. For bounded measurable $f$ (here not the drift) and any $v\in\R^d$,
\begin{equation}\label{eq:bel}
\inner{\nabla P_Tf(x)}{v}=\frac1T\,\E\Big[f(X_T^x)\int_0^T\inner{\sigma^{-1}(X_s^x)Y_sv}{\dd W_s}\Big],
\end{equation}
the Bismut--Elworthy--Li formula \citep{bismut1984large,elworthy1994formulae}. \citet{thalmaier1997differentiation} gives local versions and \citet{fournie1999applications} turned such weights into Monte Carlo estimators of sensitivities. For a $d\times m$ coefficient with $a=\sigma\sigma^\top$ invertible, $\sigma^{-1}$ is replaced by the right inverse $\sigma^\top a^{-1}$. Formula \eqref{eq:bel} computes the gradient of an expectation in the initial condition, not the gradient of the log-density in the terminal variable. \citet{mirafzali2025approach} derive instead a Bismut-type formula that targets the score directly, and the substitution formula is what turns its anticipating Skorokhod integral into adapted quantities that can be simulated along the SDE. Related integration-by-parts representations of densities, through the Malliavin--Thalmaier formula and Riesz transforms, give Monte Carlo density estimators \citep{kohatsuhiga2009estimating} and regularity estimates for densities \citep{bally2011riesz}.

\subsection{A family of unbiased score estimators}\label{sec:score-family}

Every identity of the form $\nabla\log p_t(y)=\E[\Xi\mid X_t=y]$ with computable $\Xi$ yields a regression target. \Cref{tab:ibp} collects the ones used in practice and reads each as an integration by parts. Convex combinations of unbiased targets for the same score are again unbiased, so the rows that target $\nabla\log p_t$ span a convex set of estimators within which variance can be optimised.

\begin{table}[t]
\centering\small
\begin{tabularx}{\textwidth}{@{}P{3.1cm}P{4.2cm}L@{}}
\toprule
Identity & Random target $\Xi$ (score $=\E[\Xi\mid X_t]$) & Integration by parts on \dots \\
\midrule
Tweedie / DSM \citep{vincent2011connection} & $-(X_t-\alpha_tX_0)/\sigma_t^2$ & the Gaussian noise $\epsilon$ (Stein's lemma); needs Gaussian transition kernels \\
Target score matching \citep{debortoli2024target} & $\alpha_t^{-1}\nabla\log p_0(X_0)$ & the data law; needs the data score, accurate at low noise \\
Malliavin / covering field \citep{mirafzali2025malliavin,mirafzali2025approach} & $-\Sk(u)$ for a covering field $u$ of $X_t$ & Wiener space; nonlinear drift and state-dependent noise, no closed-form kernel \\
Generalised Tweedie \citep{pidstrigach2025conditioning} & Malliavin-weighted expressions for conditional scores of nonlinear SDEs & Wiener space; robust to the singularities of diffusion bridges \\
Divergence--kernel \citep{ni2025scores,ni2025divergence} & forward covector process $\nu_t$ & the transition kernels step by step, with a free backward-adapted parameter \\
Implicit score matching \citep{hyvarinen2005estimation} & (objective, not a target) $\frac12\norm{s}^2+\tr\nabla s$ & the data density $p$ itself (\cref{prop:sm}(i)) \\
\bottomrule
\end{tabularx}
\caption{Score identities read as integrations by parts. The first five give random targets whose conditional expectation is a score, that of the marginal law except in the generalised Tweedie row, which targets the score of a conditioned process. They differ in what they require and in their variance. In the divergence--kernel row, $\nu_t$ is a covector process simulated forward with $X_t$, and the identity is rigorous for time-discretised dynamics but formal in the continuous-time limit (\cref{sec:score-estimation}). In the last row, $s$ is the model score.}
\label{tab:ibp}
\end{table}

Three remarks will be needed later. For Gaussian forward kernels, the DSM target is the special case of \cref{thm:malliavin-score} with a deterministic covering field (\cref{ex:linear}), so the Malliavin representation strictly generalises the standard objective. The minimal-energy covering field minimises $\E\norm{u}_H^2$, not the variance of $\Sk(u)$, and choosing covering fields to minimise estimator variance is an open problem with direct consequences for training (\cref{sec:score-estimation}). Finally, in infinite dimensions there is no Lebesgue density, and \eqref{eq:malliavin-score} must be read as a statement about logarithmic (Fomin) derivatives along Cameron--Martin directions \citep{bogachev2010differentiable}; \cref{sec:infinite} develops this, following \citet{mirafzali2025infinite}, \citet{mirafzali2026logarithmic} and \citet{greco2025malliavin}. Malliavin calculus is also the natural language of Gaussian approximation by Stein's method \citep{nourdin2012normal}, which suggests further tools for the quantitative analysis of diffusion models.

\begin{keybox}[Summary of Part I]
Generative dynamics live on three levels. The continuity equation turns curves in $\cP(\M)$ into velocity fields, and Girsanov's theorem turns differences of drifts into relative entropies of path measures. Time reversal, conditioning and optimal control are all changes of measure on path space. In each the drift acquires a term $a\nabla\log(\cdot)$, the diffusion matrix times the gradient of a logarithm, either of a marginal density, which gives the score term, or of a space--time harmonic function $h$, which gives the drift of the $h$-transform. These gradients can be represented by integration by parts, and Malliavin calculus is the setting in which that integration by parts works for any SDE with smooth coefficients and a nondegenerate Malliavin covariance.
\end{keybox}

\part{Deterministic flows}

\section{Normalising flows}\label{sec:flows}

The most direct way to carry a Gaussian to the data is to learn the map itself. Normalising flows do this, and they make concrete the observation of \cref{sec:transport} that many maps carry one law to another, so that the choice among them is a design decision. We describe flows built as finite compositions of maps, then their continuous-time limit, and in both cases ask what such flows can represent and at what cost.

\subsection{Finite compositions of diffeomorphisms}\label{sec:nf}

\paragraph{Likelihood and the two directions of relative entropy.}

A normalising flow is a diffeomorphism $T_\theta$ of $\R^d$ chosen from a family in which $T_\theta^{-1}$ and $\log\abs{\det\nabla T_\theta}$ are cheap \citep{tabak2010density,rezende2015variational}. With $\rho_\theta=(T_\theta)\pf\gamma_d$, the change of variables \eqref{eq:cov} gives the exact log-likelihood
\begin{align*}
\log\rho_\theta(x)=\log\gamma_d\big(T_\theta^{-1}(x)\big)+\log\abs{\det\nabla T_\theta^{-1}(x)},
\end{align*}
and maximum likelihood minimises the forward relative entropy $\KL(\mu_{\rm data}\,\|\,\rho_\theta)$. When samples are unavailable but an unnormalised target density $\pi\propto e^{-U}$ is known, the reverse relative entropy $\KL(\rho_\theta\,\|\,\pi)$ can be minimised from model samples alone. This ``training by energy'' is the main principle of Boltzmann generators for equilibrium distributions of molecular systems \citep{noe2019boltzmann}, which combine it with maximum likelihood on a few available configurations because the reverse direction alone tends to concentrate on the most stable metastable state. The two directions have the familiar mass-covering and mode-seeking behaviours, and they reappear throughout Part~IV.

\paragraph{Architectures as constrained diffeomorphisms.}

Three structural devices make the log-determinant tractable. Coupling layers split $x=(x_A,x_B)$ and set $y_A=x_A$, $y_B=\tau(x_B;c(x_A))$ with $\tau$ an invertible elementwise map, so that the Jacobian is block-triangular; additive couplings \citep{dinh2014nice}, affine couplings \citep{dinh2016density,kingma2018glow} and monotone rational-quadratic splines \citep{durkan2019neural} are the standard choices. Autoregressive layers set $y_k=\tau(x_k;c_k(x_{<k}))$ with $x_{<k}=(x_1,\dots,x_{k-1})$, realising the Knothe--Rosenblatt structure of \cref{sec:transport} \citep{kingma2016improved,papamakarios2017masked,huang2018neural}; here sampling and density evaluation trade off, since one direction is parallel and the other sequential. In both designs the conditioners $c$ and $c_k$ are unconstrained networks that output the parameters of $\tau$. Residual layers $T(x)=x+g(x)$ with $\Lip(g)<1$ are invertible by fixed-point iteration, and $\log\det(\Id+\nabla g)$ is estimated by a truncated power series with stochastic trace estimates \citep{behrmann2019invertible}, made unbiased by randomised, or ``Russian roulette'', truncation \citep{chen2019residual}. A fourth design takes the geometry of optimal transport as its starting point. Convex potential flows parameterise $T=\nabla\varphi$ with an input-convex network, so that each layer is a Brenier map by construction \citep{huang2021convex}. The first three families are reviewed by \citet{papamakarios2021normalizing} and \citet{kobyzev2020normalizing}.

\paragraph{Expressivity and its limits.}

In principle these families are highly expressive. Coupling networks whose layers include affine couplings and invertible linear maps are $L^p$-universal, and hence distributionally universal, approximators of $C^2$ diffeomorphisms \citep{teshima2020coupling}; neural autoregressive flows are universal density approximators, in the sense of convergence in distribution between laws with continuous, everywhere-positive densities \citep{huang2018neural}; and invertible networks more broadly enjoy universal approximation properties \citep{ishikawa2023universal}. Universality does not imply efficiency, however. A diffeomorphism cannot change topology, so pushing a Gaussian onto a measure with disconnected support or a hole forces the Lipschitz constant of $T$ or $T^{-1}$ to blow up; continuously indexed flows relax bijectivity to sidestep this \citep{cornish2020relaxing}. Lower bounds on depth and conditioning \citep{koehler2021representational} and limitations of simple flow families \citep{kong2020expressive} make the same point quantitatively. On the statistical side, for smooth densities bounded below on compact convex supports, maximum-likelihood estimators of the Knothe--Rosenblatt map are consistent, with relative-entropy rates $n^{-1/2}$ when $d<2s$ and $n^{-s/d}$ when $d>2s$ for smoothness $s$; when the source is strongly log-concave the expected squared Sobolev-norm error of the estimated map obeys the same rates, and without smoothness the rates can be arbitrarily slow \citep{irons2022triangular}.

\paragraph{Stochastic normalising flows.}

\citet{wu2020stochastic} interleave deterministic flow layers with stochastic MCMC or Langevin steps and train by relative entropy between the forward and backward path measures of the resulting Markov chain. This is a discrete-time precursor of diffusion samplers, closely related to annealed importance sampling and nonequilibrium transport (\cref{sec:control}). It is also, within normalising flows, an early instance of the path-space viewpoint, in which likelihoods of individual samples are replaced by Radon--Nikodym derivatives between path measures, and the balance between deterministic transport and stochastic mixing becomes a design parameter.

\subsection{Continuous normalising flows}\label{sec:cnf}

Letting the number of layers grow without bound turns a composition of maps into the flow of a vector field. The log-determinant then becomes an integral of a divergence, and the architectural constraints of the discrete case disappear.

\paragraph{Instantaneous change of variables.}

A continuous normalising flow (CNF) parameterises the velocity $v_\theta(t,x)$ of an ODE $\dot x_t=v_\theta(t,x_t)$ and sets $T_\theta=\phi_1$ \citep{chen2018neural}. Along a trajectory the log-density obeys
\begin{equation}\label{eq:inst-cov}
\frac{\dd}{\dd t}\log\rho_t(x_t)=-\divg v_\theta(t,x_t),
\end{equation}
which is the continuity equation \eqref{eq:continuity} written along characteristics, $\partial_t\rho+\inner{v}{\nabla\rho}=-\rho\,\divg v$. Equivalently, since $\rho_t(\phi_t(x))\det\nabla\phi_t(x)=\rho_0(x)$, the Jacobian obeys $\partial_t\log\det\nabla\phi_t(x)=\divg v_\theta(t,\phi_t(x))$. Nothing here depends on coordinates. The Lie derivative of the volume form along a vector field $v$ is $\divg v$ times the volume form, and the time derivative of a pull-back along a time-dependent flow is the pull-back of the Lie derivative \citep[proof of Prop.~16.33 and Prop.~22.14]{lee2013introduction}, which is why \eqref{eq:inst-cov} carries over to manifolds (\cref{sec:riemannian}). Gradients with respect to $\theta$ are obtained by the adjoint method, which integrates a backward ODE and avoids storing the forward trajectory, and the divergence is estimated by Hutchinson's identity $\divg v=\E_\varepsilon[\varepsilon^\top\nabla v\,\varepsilon]$ for $\varepsilon$ with zero mean and identity covariance, as in FFJORD \citep{grathwohl2019ffjord}. \citet{hutchinson1989stochastic} proposed independent $\pm1$ entries, the minimum-variance choice among independent zero-mean, unit-variance entries for symmetric matrices, in place of the random test vectors of Girard's earlier estimator \citep{girard1989fast}. The same adjoint ideas extend to neural SDEs \citep{li2020scalable,kidger2022neural}, which also arise as diffusion limits of deep latent Gaussian models \citep{tzen2019neural}.

\paragraph{What ODE flows can and cannot represent.}

The flow of a Lipschitz vector field is a homeomorphism, so CNFs inherit the topological obstruction described above. Augmenting the state with extra dimensions removes it \citep{dupont2019augmented}, although a density model built in this way must then marginalise over the added dimensions. Approximation results for neural ODEs and invertible residual networks characterise which maps are reachable; in $\R^d$ not every homeomorphism can be approximated, whereas every homeomorphism of $\R^d$ can be approximated by such a model acting on $\R^{2d}$ \citep{zhang2020approximation}.

\paragraph{Regularisation by transport.}

Maximum likelihood does not determine the velocity, for infinitely many fields generate the same terminal law. Selecting the field of least kinetic energy is both a regulariser and, by the Benamou--Brenier theorem (\cref{thm:bb}), a route to optimal transport. \citet{finlay2020train} penalise kinetic energy and the Frobenius norm of the Jacobian, and obtain straighter trajectories that are cheaper to integrate. OT-Flow restricts to potential velocities $v=-\nabla\Phi$ and adds a penalty enforcing the Hamilton--Jacobi equation satisfied by the optimal potential \citep{onken2021otflow}. Earlier, \citet{yang2020potential}, in work first circulated in 2019, had imposed $L^2$ optimal-transport regularity on potential generators, and the Monge--Amp\`ere flow of \citet{zhang2018monge} integrates the gradient field of a learned scalar potential, a continuous-time analogue of a Brenier map. These are the first appearances of a theme developed in full in \cref{sec:simfree}, that the choice of path in $\cP(\M)$ determines how hard the velocity is to learn and to integrate.

\paragraph{Hamiltonian and symplectic flows.}

If the state is doubled to $(q,p)$ and the velocity is Hamiltonian, the flow preserves the symplectic form \citep[Prop.~22.17]{lee2013introduction}, hence its top exterior power, which is the phase-space volume. This is Liouville's theorem, and it makes the log-determinant vanish. Hamiltonian generative networks \citep{toth2020hamiltonian} and equivariant Hamiltonian flows \citep{rezende2019equivariant} exploit this, and the same augmentation by momentum underlies kinetic diffusion models (\cref{sec:kinetic}).

\paragraph{Statistical and analytic theory.}

\citet{marzouk2024distribution} give a general nonparametric convergence analysis for distributions learned by likelihood-trained ODE models. They measure model complexity by the $C^1$ metric entropy of the velocity class and obtain rates for $C^k$-smooth target densities on the unit cube, bounded below, both for $C^k$ velocity fields and for neural networks. These rates, in squared Hellinger distance, are minimax-optimal only in a qualified sense. For neural networks the rate equals, up to a logarithm, the minimax rate for $(k-1)$-smooth functions on the $(d+1)$-dimensional space--time domain, and for $C^k$ velocity classes, when $k>d/2+3/2$, the exponent comes within an arbitrarily small margin of that rate. Both are one derivative and one dimension short of the optimal rate for $C^k$ densities on $[0,1]^d$. For flows built from Gaussian interpolation, the paths used by flow matching and diffusion models, \citet{gao2024gaussian} prove well-posedness and Lipschitz regularity of the flow under conditions on the data law, and \citet{gao2024convergence} analyse the end-to-end error of CNFs trained by flow matching. Together these results explain when the probability-flow map is a well-behaved transport, and they complement the heat-flow results of \cref{sec:heatmaps}.

\section{Simulation-free flows: matching objectives}\label{sec:simfree}

Training a CNF by maximum likelihood requires simulating the ODE and its divergence during training. The decisive simplification, made in 2022 and 2023, was to fix the curve $(\rho_t)$ in advance, as a mixture of simple conditional curves, and to learn its velocity by regression. Several groups arrived at this independently, in flow matching \citep{lipman2023flow}, rectified flow \citep{liu2023flow,liu2022rectified} and stochastic interpolants \citep{albergo2023building}.

\subsection{Flow matching and rectified flow}

The common principle is a statement about mixtures of curves. If each component of a mixture is transported by a known velocity, the mixture is transported by the conditional average of these velocities.

\paragraph{The conditional-to-marginal principle.}

\begin{proposition}[Marginalisation]\label{prop:marginalisation}
Let $Z\sim q$ be a latent variable and, for each $z$, let $(\rho_t(\cdot\mid z))$ be a curve generated by a velocity $u_t(\cdot\mid z)$, so that $\partial_t\rho_t(\cdot\mid z)+\divg(\rho_t(\cdot\mid z)u_t(\cdot\mid z))=0$. Then the mixture $\rho_t=\int\rho_t(\cdot\mid z)\,q(\dd z)$ is generated by the conditional-expectation velocity
\begin{equation}\label{eq:marginal-velocity}
u_t(x)=\E\big[u_t(X_t\mid Z)\,\big|\,X_t=x\big],\qquad (X_t,Z)\sim\rho_t(\dd x\mid z)\,q(\dd z),
\end{equation}
and for any $v_\theta$, $\E\norm{v_\theta(t,X_t)-u_t(X_t\mid Z)}^2=\E\norm{v_\theta(t,X_t)-u_t(X_t)}^2+C_t$ with $C_t$ independent of $\theta$.
\end{proposition}
\begin{proof}
Integrating the conditional continuity equations against $q$ gives $\partial_t\rho_t+\divg\int\rho_t(\cdot\mid z)u_t(\cdot\mid z)q(\dd z)=0$, and the flux $\int\rho_t(\cdot\mid z)u_t(\cdot\mid z)q(\dd z)$ equals $\rho_tu_t$ by definition of conditional expectation. The loss identity is the Pythagorean decomposition of a squared loss around a conditional mean.
\end{proof}

\begin{figure}[t]
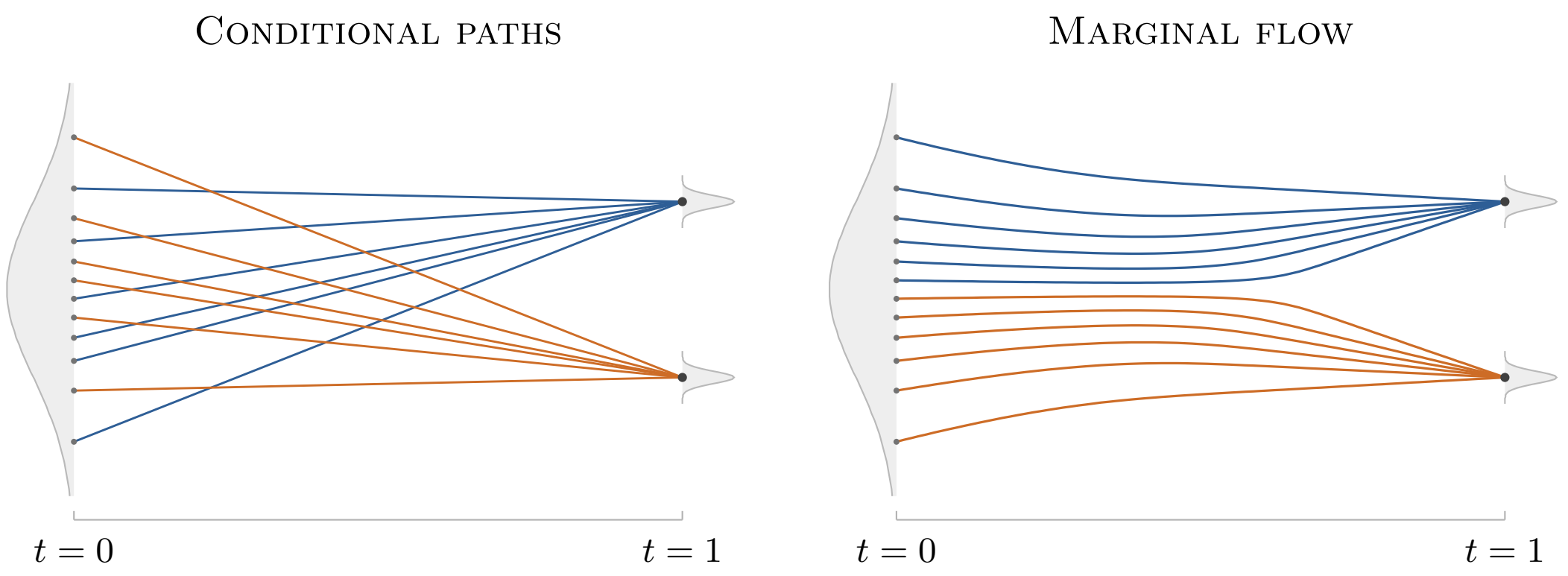

\centering

\caption{The conditional-to-marginal principle of \cref{prop:marginalisation} in one dimension, for a standard Gaussian source, data uniform on $\{-1,1\}$ and the straight-line path $x_t=(1-t)x_0+tx_1$. In the left panel, conditional paths for an independent coupling cross one another. In the right panel, the flow of the marginal velocity $u_t(x)=\E[x_1-x_0\mid x_t=x]=\big(\tanh(tx/(1-t)^2)-x\big)/(1-t)$ carries the same marginals along paths that never cross. Colours mark the endpoint.}
\label{fig:marginalisation}
\end{figure}

\Cref{fig:marginalisation} illustrates the proposition in one dimension. The same statement with the score $\nabla\log\rho_t(\cdot\mid z)$ in place of $u_t(\cdot\mid z)$ is denoising score matching (\cref{prop:sm}); with Markov generators in place of velocities it is generator matching (\cref{sec:generator-matching}). The proposition reduces generative modelling to two choices, the conditional curves and the latent coupling, followed by a regression.

\paragraph{Gaussian conditional paths.}

\citet{lipman2023flow} take $Z=x_1\sim\mu_{\rm data}$ and Gaussian conditional paths $\rho_t(\cdot\mid x_1)=\cN(m_t(x_1),s_t^2\Id)$, whose mean $m_t$ and standard deviation $s_t$ (not a score) satisfy $(m_0,s_0)=(0,1)$ and $(m_1,s_1)=(x_1,s_{\min})$ for a small final width $s_{\min}>0$. The conditional flow $x\mapsto s_tx+m_t$ has velocity
\begin{equation}\label{eq:cfm-velocity}
u_t(x\mid x_1)=\frac{\dot s_t}{s_t}\big(x-m_t(x_1)\big)+\dot m_t(x_1),
\end{equation}
and conditional flow matching regresses $v_\theta(t,X_t)$ on $u_t(X_t\mid x_1)$ with $X_t=s_tX_0+m_t(x_1)$ and $X_0\sim\gamma_d$. The diffusion paths, variance-preserving or variance-exploding, are particular choices of $(m_t,s_t)$ that meet the condition at $t=0$ only approximately. The variance-preserving path has $m_0=\bar\alpha\,x_1$ with small $\bar\alpha>0$ and $s_0=(1-\bar\alpha^2)^{1/2}$, and the variance-exploding path has $m_t\equiv x_1$ and $s_0=\sigma_{\max}\gg1$, so that its source is approximately $\cN(0,\sigma_{\max}^2\Id)$ rather than $\gamma_d$. The so-called optimal-transport path $m_t=tx_1$, $s_t=1-(1-s_{\min})t$ makes each conditional flow a straight line traversed at constant speed. The marginal flow is not straight, because conditional paths cross; the discussion of couplings below shows how changing the coupling reduces the crossings (\cref{sec:couplings}). \citet{lipman2024flow} give a guide covering the Riemannian, discrete and generator-matching extensions.

\paragraph{Rectified flow and reflow.}

Rectified flow \citep{liu2023flow} draws a coupling $(x_0,x_1)$, independent at first, and learns the velocity of the straight interpolation $x_t=(1-t)x_0+tx_1$,
\begin{align*}
v(t,x)=\E[x_1-x_0\mid x_t=x].
\end{align*}
By \cref{prop:marginalisation} its flow preserves the marginals of $x_t$. Reflow then replaces the coupling by $(x_0,\phi_1(x_0))$, the coupling induced by the learned flow, and retrains. \citet{liu2023flow} show that rectification preserves the marginals and does not increase any convex transport cost, and that its fixed points are straight couplings, whose interpolations do not cross. Repeated rectification also straightens paths, in the sense that the smallest straightness defect among the first $K$ rectified couplings is $O(1/K)$ \citep{liu2023flow}, but straight couplings need not be optimal for any prescribed cost.

To reach a particular optimum, \citet{liu2022rectified} restricts the velocity to gradient fields for quadratic cost, and, for a general convex cost $c$ of the displacement $x_1-x_0$, to fields $\nabla c^*(\nabla f_t)$, where $c^*$ is the convex conjugate of $c$ and $f_t$ a learned scalar potential. The resulting variant decreases that one cost monotonically. \citet{liu2022rectified} states that its fixed points are exactly the $c$-optimal couplings, although \citet{hertrich2025relation} give counterexamples showing that this equivalence needs stronger assumptions, such as full support of $X_t$ for $t\in(0,1)$. The smallest value of a surrogate optimality gap over the first $k$ rectifications is $O(1/k)$, while a direct bound on the cost gap is left open. The theory of straightness has since been made quantitative. \citet{bansal2025wasserstein} bound the Wasserstein error of the sampler through a piecewise straightness parameter and give sufficient conditions under which a single rectification yields a straight, or even Monge-optimal, coupling, so that the next rectified flow is exactly straight. From a standard Gaussian source they verify these conditions for Gaussian targets, for mixtures of two isotropic Gaussians with a common variance, and for mixtures of more Gaussians with a common covariance $\sigma^2\Id$ whose means lie within $2\sigma$ of one another, with a similar condition for other common covariances. For general Gaussian mixtures this is left as a conjecture. In practice, \citet{liu2024instaflow} scaled one reflow step followed by distillation to one-step text-to-image generation, with reflow making the teacher's noise--image coupling easier to distil. \citet{lee2024improving} improved the training of rectified flows, and \citet{esser2024scaling} trained rectified-flow transformers with reweighted time sampling at scale. Optimal flow matching takes a different route to the same end. It restricts to vector fields generated by convex potentials, so that a single flow-matching run, rather than repeated rectification, recovers the straight optimal-transport displacement for quadratic cost \citep{kornilov2024optimal}.

\subsection{Stochastic interpolants and couplings}

Flow matching fixes a Gaussian source and a conditional path. Stochastic interpolants free both choices and add a latent noise. Once the source is arbitrary, the coupling between source and data can also be chosen to shape the learned flow.

\paragraph{Stochastic interpolants.}\label{sec:interpolants}

\citet{albergo2023building} learn the velocity of deterministic interpolants $x_t=I_t(x_0,x_1)$ between two densities, and \citet{albergo2023stochastic} generalise to
\begin{equation}\label{eq:interpolant}
x_t=\alpha(t)\,x_0+\beta(t)\,x_1+\gamma(t)\,z,\qquad z\sim\cN(0,\Id)\ \text{independent of }(x_0,x_1),
\end{equation}
with $\alpha(0)=\beta(1)=1$ and $\alpha(1)=\beta(0)=\gamma(0)=\gamma(1)=0$. We keep the notation of \citet{albergo2023stochastic}, in which $\alpha$ multiplies the source $x_0$ and $\beta$ the data $x_1$; by contrast, $\alpha_t$ in \cref{tab:notation} multiplies the data. The scalar $\gamma(t)$ is unrelated to the Gaussian measure $\gamma_d$. The density $\rho_t$ of $x_t$ satisfies a transport equation with velocity $b(t,x)=\E[\dot\alpha x_0+\dot\beta x_1+\dot\gamma z\mid x_t=x]$, and when $\gamma(t)>0$ its score is $s(t,x)=-\gamma(t)^{-1}\E[z\mid x_t=x]$, an identity of Tweedie type. Both are learned by quadratic regression, and together they define a family of generative dynamics indexed by a noise level $\varepsilon(t)$.

\begin{proposition}[Marginal-preserving family]\label{prop:interpolant-family}
Let $\varepsilon:[0,1]\to[0,\infty)$ be continuous, and suppose that the forward SDE $\dd X_t=\big[b(t,X_t)+\varepsilon(t)s(t,X_t)\big]\dd t+\sqrt{2\varepsilon(t)}\,\dd W_t$ with $X_0\sim\rho_0$ is well posed and that its Fokker--Planck equation has a unique solution with initial law $\rho_0$. Then the SDE has marginals $\rho_t$. So does the backward SDE with drift $b-\varepsilon s$, integrated backwards in time from $X_1\sim\rho_1$, under the analogous conditions on its Fokker--Planck equation with terminal law $\rho_1$.
\end{proposition}
\begin{proof}
The Fokker--Planck equation of the SDE reads $\partial_t\rho=-\divg\big(\rho(b+\varepsilon s)\big)+\varepsilon\Delta\rho$, and $\varepsilon\Delta\rho=\divg(\varepsilon\rho s)$ cancels the extra drift, leaving the transport equation satisfied by $\rho_t$. Hence $\rho_t$ solves the Fokker--Planck equation, and uniqueness identifies it with the law of $X_t$. For the backward SDE the Fokker--Planck equation, solved backwards from $\rho_1$, reads $\partial_t\rho=-\divg\big(\rho(b-\varepsilon s)\big)-\varepsilon\Delta\rho$, and the same cancellation leaves the transport equation.
\end{proof}
Setting $\varepsilon\equiv0$ gives the ODE, and diffusion models are recovered from the one-sided interpolant, in which $x_0\sim\gamma_d$ is independent of $x_1$ and $\gamma\equiv0$. The noise level $\varepsilon(t)$ can be tuned at sampling time to trade discretisation error against the propagation of estimation error, and the resulting relative entropy is bounded by time-integrated errors in $b$ and $s$ \citep{albergo2023stochastic}. The authors derive this bound from the Fokker--Planck equations. Girsanov's formula \eqref{eq:girsanov-kl} gives the same drift-error bound on path space. Both bounds carry a factor $1/\varepsilon$, so they require $\varepsilon>0$ and do not extend to the deterministic limit. Data-dependent couplings $(x_0,x_1)$ turn the framework into a tool for conditional generation and inverse problems \citep{albergo2024stochastic}, and several marginals can be interpolated simultaneously \citep{albergo2024multimarginal}. For probabilistic forecasting the base can be a point mass at the current state \citep{chen2024probabilistic}. An SDE started there reaches the conditional law of the future state in finite time, and its drift is learned by regression with a loss of bounded variance. Its diffusion can be retuned after training. The choice that minimises the path-space relative entropy due to estimation error gives, under a monotonicity condition on the schedule, a F\"ollmer process (\cref{sec:follmer}).

\paragraph{Couplings.}\label{sec:couplings}

\Cref{prop:marginalisation} allows any coupling $\pi$ of $\mu_{\rm ref}$ and $\mu_{\rm data}$. If $\pi$ is the optimal coupling for quadratic cost and the conditional paths are straight segments, the segments do not cross, being McCann displacement interpolants; the marginal velocity then transports along them, and the flow at $t=1$ is the Brenier map. Exact optimal couplings are unavailable, but solving optimal transport within each minibatch gives straighter flows and lower-variance regression targets \citep{tong2024improving,pooladian2023multisample}. With an entropic coupling of regularisation $\varepsilon$ for the cost $\abs{x-y}^2/2$, and Brownian bridges of variance $\varepsilon$ in place of straight segments, the same construction approximates Schr\"odinger bridges (\cref{sec:sb}). Conditional variants require care with the conditioning variable, and \citet{chemseddine2024conditional} introduce conditional Wasserstein distances that justify optimal-transport flow matching for posterior sampling. A variational reading of flow matching, which approximates the posterior over endpoints rather than the velocity, is given by \citet{eijkelboom2024variational} and is natural for categorical data.

\subsection{Generalisations}

The marginalisation principle does not depend on the conditional curves being Gaussian interpolations, nor on the dynamics being a velocity field. It extends to gradient velocities learned from marginals alone, to general Markov generators, and to curves dictated by physical equations.

\paragraph{Action matching and gradient velocities.}

Among the velocities generating a curve, the one with least kinetic energy is a gradient $\nabla\psi_t$ (\cref{sec:transport}). Action matching learns $\psi_t$ directly from samples of the marginals. A neural potential is fitted by minimising a tractable objective that equals, up to a constant, the time-integrated squared distance between its gradient and this gradient field \citep{neklyudov2023action}. Action matching needs only samples of $\rho_t$ at each time, not a coupling, which makes it a tool for learning dynamics from snapshots. For Gaussian conditional paths $\cN(\alpha_tx_1,s_t^2\Id)$ with an independent coupling, the marginal flow-matching velocity \eqref{eq:marginal-velocity} is itself a gradient, being affine in $x$ and $\nabla\log\rho_t$ by Tweedie's formula. For general couplings and non-Gaussian endpoints it need not be, and this is why restricting to gradient fields changes the fixed points of rectification \citep{liu2022rectified}.

\paragraph{Generator matching.}\label{sec:generator-matching}

\citet{holderrieth2025generator} replace velocities by Markov generators. If conditional processes with generators $\cL_t^{z}$ realise conditional curves $\rho_t(\cdot\mid z)$, then the generator $\cL_t=\E[\cL_t^Z\mid X_t]$, defined through a linear parameterisation of generators by drift, diffusion and jump kernel, realises the mixture; this is the analogue of \cref{prop:marginalisation}. Flows, diffusions, jump processes and superpositions of them are thereby trained by a single regression principle, and generator matching covers the discrete-state models of \cref{sec:discrete} and the L\'evy-driven models of \cref{sec:levy}.

\paragraph{Physics-guided paths.}\label{sec:physics-paths}

The curve $(\rho_t)$ need not come from a Gaussian interpolation at all. Poisson flow generative models treat data as charges in an augmented space and move samples along electric field lines \citep{xu2022poisson}, and PFGM++ interpolates between this construction and diffusion through the dimension of the augmentation \citep{xu2023pfgm}. GenPhys asks which physical PDEs beyond diffusion and electrostatics can be turned into generative models \citep{liu2023genphys}. In holographic generative flows, data are encoded in Fourier space as boundary sources of a field in anti-de Sitter space, and the flow-matching dynamics are staged in a normalised radial phase space with a Klein--Gordon backbone and a learned residual correction \citep{mirafzali2026holographic}; \cref{sec:holography} returns to this construction.

\subsection{Theory of flow matching}

Error analysis for diffusion models compares path measures through their drifts by Girsanov's theorem (\cref{sec:theory}); for deterministic flows this argument applies only in part, because they lack the smoothing effect of noise. \citet{benton2024error} bound the $\W_2$ distance between the target and the output of the learned flow by the $L^2$ velocity error times the exponential of the time-integrated spatial Lipschitz constant of the learned velocity. They control this constant for the true velocity of interpolants between compactly supported distributions, provided a Gaussian smoothing term $\gamma(t)z$ with $\gamma(t)>0$ is present at every time, endpoints included (so the condition $\gamma(0)=\gamma(1)=0$ of \eqref{eq:interpolant} is relaxed), and the intermediate mixtures $\alpha(t)x_0+\beta(t)x_1$ satisfy a regularity condition.

Statistical rates have followed in several forms. Flow matching attains nearly minimax rates in $\W_p$ for $1\le p\le2$ over Besov densities bounded below on a cube (and smoother near its boundary), provided the standard deviation of the Gaussian component scales like $\sqrt{1-t}$ as $t\to1$, separate networks are trained on a dyadic partition of time, and the posterior mean has a uniformly bounded Jacobian; faster variance decay gives a slower bound \citep{fukumizu2025flow}. Through a connection with kernel density estimation it is minimax optimal up to logarithmic factors over suitable classes of densities, with faster rates, on a small tube around the manifold, when the data lie on a low-dimensional single-chart manifold \citep{kunkel2025minimax}, and control of the Lipschitz constant of the velocity gives $\W_1$ rates for classes of unbounded, not necessarily log-concave distributions \citep{kunkel2025distribution}. A deterministic, non-asymptotic bound on relative entropy gives near-minimax total-variation rates comparable to those of diffusion models \citep{su2025flow}, and for CNFs trained by flow matching \citet{gao2024convergence} combine approximation, estimation and discretisation errors.

Adding noise to the interpolation brings back the Girsanov route to bounds in relative entropy. When the interpolation between the two endpoints of a coupling is a Brownian bridge rather than a straight line, the Markovian projection of the resulting path measure has the same marginals and is learned as a drift. For this diffusion flow matching, the Euler--Maruyama scheme has relative-entropy error of order $\varepsilon^2+h\,(d^4+M_8)$ for $L^2$ drift error $\varepsilon$ and step size $h$, without early stopping, where $M_8$ collects eighth moments of the two endpoint laws and of the scores of these laws and of the coupling \citep{gentilonisilveri2024theoretical}. A sharper analysis lowers the dimension dependence of the discretisation term from $d^4$ to $d^3$ and needs only the score of the coupling, or, after early stopping, that of the conditional law of the source point given the data point \citep{gentilonisilveri2026diffusion}. It also gives $\W_2$ bounds with the same dimension dependence, under weak log-concavity of the coupling and an integrability condition on its Hessian, at the price of a constant exponential in these regularity parameters.

\section{Flow maps and few-step generation}\label{sec:flowmaps}

\begin{figure}[t]
\centering
\begin{tikzpicture}[lab/.style={font=\footnotesize}]
\draw[->] (0,0) -- (11.4,0) node[right,lab] {time};
\draw (0,0) -- (0,-0.06) node[below,lab] {$0$, noise};
\draw (11.0,0) -- (11.0,-0.06) node[below,lab] {$1$, data};
\draw[black!40] (2.75,0) -- (2.75,-0.06) node[below,lab,text=black] {$s$};
\draw[black!40] (8.25,0) -- (8.25,-0.06) node[below,lab,text=black] {$t$};
\draw[black!25,dotted] (2.75,0) -- (2.750,1.255);
\draw[black!25,dotted] (8.25,0) -- (8.250,2.272);
\draw[figblue,line width=1.1pt] plot coordinates {(0.000,1.410) (0.040,1.407) (0.080,1.405) (0.123,1.402) (0.163,1.399) (0.205,1.396) (0.245,1.393) (0.288,1.391) (0.328,1.388) (0.368,1.385) (0.411,1.383) (0.451,1.380) (0.493,1.377) (0.533,1.375) (0.576,1.372) (0.616,1.369) (0.659,1.367) (0.699,1.364) (0.739,1.361) (0.782,1.359) (0.822,1.356) (0.864,1.354) (0.904,1.351) (0.947,1.348) (0.987,1.346) (1.030,1.343) (1.070,1.341) (1.110,1.339) (1.152,1.336) (1.192,1.334) (1.235,1.331) (1.275,1.329) (1.318,1.326) (1.358,1.324) (1.400,1.322) (1.440,1.319) (1.480,1.317) (1.523,1.315) (1.563,1.312) (1.606,1.310) (1.646,1.308) (1.689,1.305) (1.729,1.303) (1.771,1.301) (1.811,1.299) (1.851,1.297) (1.894,1.294) (1.934,1.292) (1.977,1.290) (2.017,1.288) (2.059,1.286) (2.099,1.284) (2.142,1.282) (2.182,1.280) (2.222,1.278) (2.265,1.276) (2.305,1.274) (2.347,1.272) (2.387,1.270) (2.430,1.268) (2.470,1.266) (2.513,1.265) (2.553,1.263) (2.593,1.261) (2.635,1.259) (2.676,1.258) (2.718,1.256) (2.758,1.254) (2.801,1.253) (2.841,1.251) (2.884,1.249) (2.924,1.248) (2.964,1.246) (3.006,1.245) (3.046,1.243) (3.089,1.242) (3.129,1.241) (3.172,1.239) (3.212,1.238) (3.254,1.237) (3.294,1.236) (3.334,1.234) (3.377,1.233) (3.417,1.232) (3.460,1.231) (3.500,1.230) (3.542,1.229) (3.582,1.228) (3.625,1.227) (3.665,1.226) (3.705,1.226) (3.748,1.225) (3.788,1.224) (3.831,1.224) (3.871,1.223) (3.913,1.222) (3.953,1.222) (3.996,1.222) (4.036,1.221) (4.076,1.221) (4.119,1.221) (4.159,1.221) (4.201,1.221) (4.241,1.221) (4.284,1.221) (4.324,1.221) (4.367,1.221) (4.407,1.221) (4.447,1.222) (4.489,1.222) (4.529,1.223) (4.572,1.223) (4.612,1.224) (4.655,1.225) (4.695,1.226) (4.737,1.227) (4.777,1.228) (4.818,1.229) (4.860,1.231) (4.900,1.232) (4.943,1.234) (4.983,1.236) (5.026,1.238) (5.066,1.240) (5.108,1.242) (5.148,1.245) (5.188,1.247) (5.231,1.250) (5.271,1.253) (5.314,1.256) (5.354,1.259) (5.396,1.263) (5.436,1.266) (5.479,1.270) (5.519,1.274) (5.559,1.278) (5.602,1.283) (5.642,1.288) (5.684,1.293) (5.724,1.298) (5.767,1.304) (5.807,1.310) (5.850,1.316) (5.890,1.323) (5.930,1.329) (5.973,1.337) (6.013,1.344) (6.055,1.352) (6.095,1.361) (6.138,1.370) (6.178,1.379) (6.221,1.388) (6.261,1.398) (6.301,1.408) (6.343,1.420) (6.383,1.431) (6.426,1.443) (6.466,1.455) (6.509,1.468) (6.549,1.481) (6.591,1.495) (6.631,1.509) (6.671,1.523) (6.714,1.539) (6.754,1.554) (6.797,1.571) (6.837,1.587) (6.879,1.604) (6.919,1.621) (6.962,1.639) (7.002,1.657) (7.042,1.675) (7.085,1.694) (7.125,1.713) (7.168,1.733) (7.208,1.751) (7.250,1.772) (7.290,1.791) (7.333,1.812) (7.373,1.831) (7.413,1.851) (7.456,1.872) (7.496,1.892) (7.538,1.913) (7.578,1.933) (7.621,1.955) (7.661,1.975) (7.704,1.996) (7.744,2.016) (7.784,2.037) (7.826,2.058) (7.866,2.078) (7.909,2.100) (7.949,2.120) (7.992,2.142) (8.032,2.162) (8.075,2.183) (8.115,2.203) (8.155,2.224) (8.197,2.245) (8.237,2.265) (8.280,2.287) (8.320,2.307) (8.363,2.329) (8.403,2.349) (8.445,2.370) (8.485,2.391) (8.525,2.411) (8.568,2.432) (8.608,2.453) (8.651,2.474) (8.691,2.494) (8.733,2.516) (8.773,2.536) (8.816,2.558) (8.856,2.578) (8.896,2.598) (8.939,2.620) (8.979,2.640) (9.021,2.661) (9.061,2.681) (9.104,2.703) (9.144,2.723) (9.187,2.745) (9.227,2.765) (9.267,2.785) (9.310,2.807) (9.350,2.827) (9.392,2.848) (9.432,2.869) (9.475,2.890) (9.515,2.910) (9.558,2.932) (9.598,2.952) (9.638,2.972) (9.680,2.994) (9.720,3.014) (9.763,3.036) (9.803,3.056) (9.846,3.077) (9.886,3.098) (9.928,3.119) (9.968,3.139) (10.008,3.159) (10.051,3.181) (10.091,3.201) (10.134,3.223) (10.174,3.243) (10.217,3.265) (10.257,3.285) (10.299,3.306) (10.339,3.326) (10.379,3.347) (10.422,3.368) (10.462,3.388) (10.505,3.410) (10.545,3.430) (10.587,3.452) (10.627,3.472) (10.670,3.493)};
\draw[black!70,dashed,line width=0.8pt] (2.750,1.255) -- (8.250,2.272);
\draw[figorange,line width=1.2pt,->,>=Stealth] (2.750,1.255) -- (3.960,1.206);
\fill[black] (2.750,1.255) circle (1.7pt) node[below left,lab] {$x_s$};
\fill[black] (8.250,2.272) circle (1.7pt) node[below right,lab] {$x_t$};
\draw[->,>=Stealth,black!80,line width=0.8pt] (2.750,1.495) to[out=60,in=150] node[midway,above,lab] {flow map $X_{s,t}$} (8.140,2.482);
\node[lab,text=figorange,anchor=north west] at (2.882,1.150) {tangent slope $v_s(x_s)$};
\draw[black!50] (7.700,2.170) -- (8.800,1.800);
\node[lab,anchor=west,text=black!80] at (8.800,1.800) {chord slope (average velocity)};
\end{tikzpicture}
\caption{Three objects attached to one trajectory of the marginal velocity of \cref{fig:marginalisation}, started at $x_0=0.25$. The velocity $v_s(x_s)$ is the slope of the tangent at time $s$, the object that flow matching learns. The flow map $X_{s,t}$ sends $x_s$ to $x_t$ along the trajectory (\cref{prop:flowmap}), and the slope of the chord, $(x_t-x_s)/(t-s)$, is the average velocity over $[s,t]$. MeanFlow, whose time runs from data to noise, learns this average velocity as a function of the noisier endpoint, here $x_s$, with the opposite sign because its time runs backwards. Few-step methods learn the last two, which remain exact however large the step.}
\label{fig:flowmap}
\end{figure}
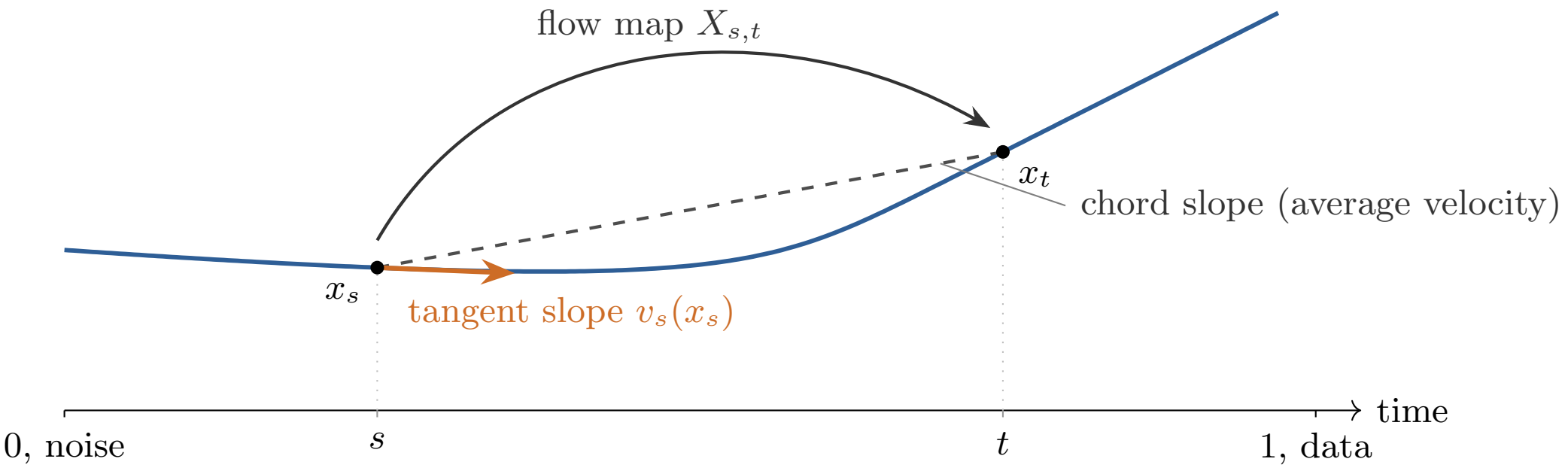

Flow matching learns an Eulerian object, the velocity $v_t(x)$, and generation integrates it numerically. Each integration step costs at least one evaluation of a large network, so the number of steps sets the cost of a sample. Few-step methods instead learn Lagrangian objects, maps that jump along trajectories (\cref{fig:flowmap}). The flow map can be characterised in several equivalent ways, each of which suggests a different training objective. Several of the methods of this section can be classified by which characterisation they enforce.

\paragraph{The flow map and its three characterisations.}

Let $X_{s,t}(x)$ denote the solution at time $t$ of $\dot x=v_t(x)$ started from $x$ at time $s$. For smooth $v$, each $X_{s,t}$ is a diffeomorphism from its open domain onto its image, with inverse $X_{t,s}$, and $X_{u,t}\circ X_{s,u}=X_{s,t}$ wherever both sides are defined \citep[Thm.~9.48]{lee2013introduction}.

\begin{proposition}[Characterisations of the two-time flow map \citep{boffi2024flow,boffi2025build}]\label{prop:flowmap}
For a smooth velocity $v$, the map $X_{s,t}$ is the unique family with $X_{s,s}=\Id$ satisfying any one of
\begin{align*}
&\text{(Lagrangian)} && \partial_tX_{s,t}(x)=v_t\big(X_{s,t}(x)\big),\\
&\text{(Eulerian)} && \partial_sX_{s,t}(x)+\nabla X_{s,t}(x)\,v_s(x)=0,\\
&\text{(semigroup)} && X_{t',t}\circ X_{s,t'}=X_{s,t}\ \text{for all } s,t',t,\ \text{together with}\ \partial_tX_{s,t}\big|_{t=s}=v_s.
\end{align*}
\end{proposition}

Each characterisation yields a training objective, either from scratch or by distillation from a learned velocity. \citet{boffi2024flow} show that the continuous-time limits of consistency distillation and of consistency trajectory models are Eulerian distillation losses, and that progressive distillation enforces the semigroup property. \citet{boffi2025build} turn each characterisation into a self-distillation loss, Eulerian, Lagrangian or progressive, the last built on the semigroup property. They find the Lagrangian family, which avoids spatial derivatives and bootstrapping from small steps, markedly more stable.

\paragraph{Consistency models.}

Consistency models are written in the diffusion convention, with the data near $t=0$. A consistency model learns the one-time map $f_\theta(x,t)\approx X_{t,\delta}(x)$ from any point of a probability-flow trajectory to its endpoint at a small time $\delta>0$ near the data, with the boundary condition $f_\theta(x,\delta)=x$ built into the parameterisation, and enforces $f_\theta(x_t,t)=f_{\theta^-}(x_{t'},t')$ for adjacent points of the same trajectory, where $t'<t$ is the less noisy time and $\theta^-$ is a copy of $\theta$ through which no gradient passes, such as a running average of its past values \citep{song2023consistency}. Adjacent points come either from a pretrained ODE solver, which is distillation, or from sharing the same noise draw, which is consistency training. The idea has been refined in three directions, through better training recipes \citep{song2024improved}, continuous-time formulations stabilised at scale \citep{lu2025simplifying}, and two-time trajectory models that learn $X_{t,s}$ for all pairs $s\le t$, from any time to any earlier time \citep{kim2024consistency}. A consistency model must map every point of a trajectory to its endpoint, so it has to denoise nearly clean inputs as well as generate from noise, and its capacity is shared between accuracy at small times and accuracy at the noise end, which is what one-step generation uses. Truncated consistency models settle this trade-off in favour of the noise end by training only on $[t_{\rm tr},T]$, with the boundary condition at $t_{\rm tr}$ supplied by a first-stage consistency model, since without one a constant map minimises the truncated loss; at similar network size this lowers the one-step FID of improved consistency training from 2.83 to 2.46 on CIFAR-10 and from 4.02 to 2.88 on ImageNet at $64\times64$ \citep{lee2025truncated}. Convergence guarantees given accurate consistency losses \citep{lyu2023convergence} and rates for consistency training \citep{li2024consistency} are available.

\paragraph{Average velocities: shortcut models and MeanFlow.}

A second family parameterises the two-time map by an average velocity evaluated at the later point. Writing $(z_\tau)$ for the trajectory of $v$ through $z_t$, set $u(z_t,r,t)=\frac{1}{t-r}\int_r^tv_\tau(z_\tau)\dd\tau$, the average velocity over $[r,t]$, so that $z_r=X_{t,r}(z_t)=z_t-(t-r)\,u(z_t,r,t)$. This is the convention of MeanFlow, whose time runs from data at $t=0$ to noise at $t=1$, the reverse of the flow convention of \cref{sec:notation}. In MeanFlow's convention the later point $z_t$ is the noisier one, $X_{t,r}$ with $r<t$ moves towards the data, and one-step sampling reads $z_0=z_1-u(z_1,0,1)$ with $z_1=\epsilon$. The following identity is stated in this convention.

\begin{proposition}[MeanFlow identity \citep{geng2025mean}]\label{prop:meanflow}
With $r$ fixed and $z_t$ moving along a trajectory of $v$,
\begin{align*}
u(z_t,r,t)&=v_t(z_t)-(t-r)\,\frac{\dd}{\dd t}u(z_t,r,t),\\
\frac{\dd}{\dd t}u&=\partial_tu+\nabla_zu\;v_t(z_t).
\end{align*}
\end{proposition}
\begin{proof}
Differentiate $(t-r)\,u(z_t,r,t)=\int_r^tv_\tau(z_\tau)\dd\tau$ with respect to $t$ along the trajectory.
\end{proof}
MeanFlow trains $u_\theta$ to satisfy this identity, using the conditional flow-matching velocity in place of $v_t$, both in the target and as the tangent of the Jacobian--vector product that computes the total derivative. The network terms in the target are held fixed by a stop-gradient, so the target is affine in $v_t$ and the substitution leaves the expected update, and hence its fixed points, unchanged. At $r=t$ the loss reduces to flow matching. MeanFlow achieves one-step generation from scratch. Shortcut models enforce the semigroup property instead, on dyadic step sizes \citep{frans2025one}. The structure of the MeanFlow objective has itself been analysed. AlphaFlow decomposes it into a trajectory flow-matching part and a consistency part whose gradients conflict during training, places trajectory flow matching, shortcut models and MeanFlow in one family of losses, and anneals from the first to the last by a curriculum \citep{zhang2025alphaflow}. Improved MeanFlow recasts the objective of \cref{prop:meanflow} as a regression of the instantaneous velocity parameterised through $u_\theta$, so that the target no longer depends on the network, and evaluates the Jacobian--vector product along a predicted marginal velocity, for instance $u_\theta(z_t,t,t)$, rather than the conditional one; with guidance supplied as a conditioning variable, it reaches a one-step FID of 1.72 on ImageNet at $256\times256$ without distillation \citep{geng2025improved}. Without a latent space, pixel MeanFlow has the network predict a clean image, which lies near a low-dimensional manifold, while the loss remains the MeanFlow loss on velocities, related to the prediction by an explicit conversion; with a perceptual loss added, it reaches one-step FIDs of 2.22 and 2.48 on ImageNet at $256\times256$ and $512\times512$ \citep{lu2026one}.

Inductive moment matching works at the level of distributions rather than trajectories, matching laws across pairs of times by a kernel discrepancy, and trains stably from scratch \citep{zhou2025inductive}. \citet{sabour2025align} develop flow-map distillation at scale, and \citet{cao2025highorder} analyse the expressivity and efficiency of second-order extensions with average accelerations.

Other work asks what a one-step generator computes. Decoding the intermediate layers of a one-step model with its own output head recovers, across depth, the denoising that a multistep sampler performs across steps \citep{asiimwe2026depth}. The effect depends on the transport task the model is trained for. MeanFlow, probed on short intervals, shows denoising followed by renoising, while drifting models (\cref{sec:gradflows}) show no such progression. Models that exhibit the effect compress well; a MeanFlow model shrinks by $16.6\times$ in parameters into a single time-conditioned block, at an FID of 4.9 against 4.0 after post-training. Conditioning a one-step map, which leaves no trajectory to guide, is taken up with inverse problems in \cref{sec:conditioning}.

\paragraph{Distillation.}

Distillation compresses a pretrained sampler. The first methods worked at the level of trajectories, by regressing directly onto the endpoints of a teacher's ODE \citep{luhman2021knowledge} and by progressively halving the number of steps \citep{salimans2022progressive}. Later methods match distributions instead. Distribution matching distillation descends an approximate relative entropy between the law of the one-step generator and that of the pretrained diffusion model; its gradient is the difference between a score fitted to generated samples and the pretrained model's score, both evaluated at noisy inputs \citep{yin2024one}. Score identity distillation derives its loss from score identities of the forward process \citep{zhou2024score}, and moment matching distils multistep samplers \citep{salimans2024multistep}. In the language of \cref{sec:ibp}, distribution matching distillation and score identity distillation need the score of the generator's own output law, and they obtain it, or eliminate it, through the score identities of \cref{tab:ibp}.

\paragraph{Stochastic transitions.}

Between a deterministic flow map and a full diffusion lies a chain of a few stochastic steps. Transition matching learns the kernels of a discrete-time Markov chain from noise to data by matching the transitions of a supervising process, which may be any process between the two, even a discontinuous one. In its difference version the supervising process is the linear interpolation, and each step samples the increment $X_1-X_0$ given the current state with a small flow-matching head on top of a large backbone. The mean of a step is the Euler step of flow matching, and the chain converges to flow matching as the number of steps grows; with 16 to 32 transitions, against 128 to 256 flow-matching steps, it improves image quality and prompt adherence over flow matching at fixed architecture, data and hyperparameters, while sampling up to almost seven times faster \citep{shaul2025transition}. The same structure serves distillation. Transition matching distillation compresses large video models into few-step generators whose transitions are conditional flows computed by the last layers of the network, and trains them with distribution matching \citep{nie2026transition}.

Stochastic flow maps instead push fresh noise to a conditional law, one for each noisy state. Meta flow maps amortise, over every noisy state $(t,x)$, the flow maps that carry noise to the posterior $p_{1|t}(\cdot\mid x)$ of the clean endpoint, so that one network call returns an independent posterior sample; these samples give differentiable estimates of the value function $V_t(x)=\log\E[e^{r(X_1)}\mid X_t=x]$, an expectation over $p_{1|t}(\cdot\mid x)$, for alignment to a reward $r$ (\cref{sec:control}), which allow steering without inner simulations and unbiased off-policy fine-tuning \citep{potaptchik2026meta}. This $V_t$ measures reward rather than cost. It is $-V$ for the value function $V$ of \cref{sec:follmer}, with the base process as the uncontrolled one and $\Phi=-r$ as terminal cost. Diamond Maps distil such posterior samplers from stochastic transitions of a pretrained flow, or turn an existing deterministic flow map into a consistent estimator of the value function, so that search, sequential Monte Carlo and guidance become cheap \citep{holderrieth2026diamond}.

For a stochastic differential equation with additive noise, the analogue of the flow map is the It\^o map $(X_s,W_{[s,t]})\mapsto X_t$, which is continuous in the Brownian path. Strong stochastic flow maps learn this map \citep{mccallum2026strong}. Their parameterisation resembles an Euler--Maruyama step whose coefficients depend on the Brownian path. The Brownian path enters the network through a polynomial approximation of Brownian motion that converges in H\"older norm and composes by explicit Chen relations. Unlike the stochastic flow maps above, which target transition distributions, they reproduce trajectories pathwise rather than only in law. They generate images in a few stochastic steps and equilibrium molecular conformations in one or two.

\section{Gradient-flow generative models}\label{sec:gradflows}

The models of the preceding sections prescribe a curve of laws and learn its velocity. A different idea lets the target choose the curve. Particles move along the gradient flow of a discrepancy between their current law and the target, in one of the geometries of \cref{sec:transport}. The choice of discrepancy and of geometry then determines the method, and the four families below correspond to four such choices. A last paragraph turns to models that do without the noise level.

\paragraph{Wasserstein gradient flows and JKO-trained flows.}

In the Wasserstein geometry each step of the flow is a JKO step \eqref{eq:jko}, which can be solved with input-convex networks parameterising the transport map as the gradient of a convex potential \citep{mokrov2021large,alvarezmelis2022optimizing}, or variationally \citep{fan2022variational}. JKO-iFlow trains a normalising flow block by block, each residual block one JKO step of the relative entropy with respect to the Gaussian reference, so that the trained flow carries data to noise and its inverse generates \citep{xu2023normalizing}. Convergence guarantees for JKO-iFlow are obtained through proximal gradient descent in Wasserstein space \citep{cheng2024convergence}. An early generative approach of this kind pushes particles along the negative gradient of the first variation of an $f$-divergence \citep{gao2019deep}.

\paragraph{Kernel and sliced discrepancies.}

When the discrepancy is a maximum mean discrepancy, the gradient flow moves particles along the gradient of the witness function. Because the MMD is not displacement convex, convergence to the target holds only under conditions, and injecting noise into the gradient is the proposed remedy \citep{arbel2019maximum}. Sliced-Wasserstein flows evolve particles by one-dimensional transports along random directions, together with a Brownian term from an entropic regulariser \citep{liutkus2019sliced}, and sliced MMD flows with Riesz kernels scale to images \citep{hertrich2024generative}.

Drifting models carry this idea into training \citep{deng2026generative}. A one-step generator is regressed onto its own samples displaced by a drifting field, held fixed at each step, which attracts samples towards the data and repels them from the generator's own samples through normalised kernels; no flow is integrated at inference. The field is antisymmetric in the two laws, so it vanishes when they coincide. The converse fails for general fields, and the authors give sufficient conditions, in a finite-dimensional model class, under which zero drift forces equality. Gradients of an MMD give a special case, although the normalised kernels do not arise naturally from an MMD. With the field computed in the features of a pretrained masked-autoencoder encoder, fine-tuned for classification, the method reaches a one-step FID of 1.54 on ImageNet at $256\times256$. \citet{li2026longshort} derive the field from the semigroup property of flow maps (\cref{prop:flowmap}). They compose a long flow map with a short terminal map evaluated with the closed-form flow-matching velocity of the data. As the terminal interval shrinks, a first-order treatment of that step yields the attraction alone and a second-order one the full attraction--repulsion field. In their derivation the law of the noise fixes the kernel. Gaussian noise gives a Gaussian kernel, and Laplace noise the Laplacian kernel of the original method. Wasserstein gradient flows give the construction a variational form. W-Flow takes the drifting field to be the velocity of the Wasserstein gradient flow of a Sinkhorn divergence to the data, so that each training step is an explicit Euler step of that flow. The field is the difference between the entropic optimal-transport barycentric projections of a sample onto the data and onto the model's own samples, the latter estimated with an independent batch. This keeps the attraction--repulsion structure of drifting models and removes their diagonal masking. Under suitable assumptions, the finite-sample training dynamics converge to the continuous gradient flow, and with the same pretrained feature encoder as drifting models the method reaches a one-step FID of 1.29 on ImageNet at $256\times256$ \citep{han2026one}.

\paragraph{Stein geometry.}

Stein variational gradient descent moves particles with a kernelised Wasserstein gradient of relative entropy \citep{liu2016stein,liu2017stein}. Its geometry is that of a Riemannian metric induced by the kernel \citep{duncan2023geometry}, and it can equally be read as a kernelised Wasserstein gradient flow of the chi-squared divergence \citep{chewi2020svgd}.

\paragraph{Fisher--Rao and Wasserstein--Fisher--Rao dynamics.}

Adding birth--death moves, which create and destroy mass, to Langevin dynamics implements the Fisher--Rao part of the Wasserstein--\hspace{0pt}Fisher--\hspace{0pt}Rao flow and accelerates convergence for multimodal targets \citep{lu2019accelerating}. For Gaussian location mixtures, Wasserstein--Fisher--Rao gradient descent on the negative log-likelihood, run with a small enough step size from a fully supported initialisation, can only converge to a nonparametric maximum-likelihood estimator (any weak limit is one), while its interacting-particle implementation reaches approximate estimators numerically \citep{yan2024learning}. Which geometry and which energy to choose is itself a question. \citet{chen2023gradient} compare gradient flows for sampling and argue for relative entropy and affine-invariant geometries, and kernel Fisher--Rao flows transport a reference to the target in unit time \citep{maurais2024sampling}. Finally, the generator of the Langevin dynamics can be approximated from samples by diffusion maps, and a spectral approximation of that generator supplies the kernel of a Laplacian-adjusted Wasserstein gradient descent, which gives a particle generative model built directly from samples \citep{li2023diffusion}.

\paragraph{Equilibrium landscapes.}

Diffusion and flow models feed the noise level to the network, yet the network can usually infer it from its input. Without it, the regression target is a posterior mean over noise levels, which in high dimension are nearly determined by the noisy input, and most denoising generative models degrade only mildly, with a bound on the accumulated sampling error; the deterministic sampler of DDIM is the exception, failing in line with a much larger bound, and a noise-unconditional model reaches an FID of 2.23 on CIFAR-10 \citep{sun2025noise}. Equilibrium matching drops time altogether. It regresses a single field onto $(\epsilon-x)\,c(\gamma)$ at $x_\gamma=\gamma x+(1-\gamma)\epsilon$, without access to $\gamma$, where $x\sim\mu_{\rm data}$, $\gamma\in[0,1]$ is drawn at random and the weight $c$ vanishes at $\gamma=1$. The target therefore vanishes at the data, and so, approximately and in high dimension, does the learned field, which is read as the gradient of an implicit energy. The choice $c\equiv1$ recovers the flow-matching objective up to sign. Samples are produced by gradient descent on that energy, with adaptive step sizes, optimisers and stopping, and reach an FID of 1.90 on ImageNet at $256\times256$. The supporting analysis is informal and concerns the exact minimiser, whose local minima are approximately the training samples \citep{wang2025equilibrium}; the generalisation that the authors observe must therefore, as for diffusion models (\cref{sec:generalisation}), rest on imperfect fitting. Energy matching instead learns a scalar potential whose gradient follows optimal-transport paths far from the data and relaxes to a Boltzmann law near it, and the potential doubles as a prior for inverse problems \citep{balcerak2025energy}.

\begin{openprob}[Open problems for Part II]
\begin{problem}
Characterise the straightness of the marginal flow-matching velocity in terms of the geometry of $\mu_{\rm data}$, beyond the Gaussian and Gaussian-mixture targets of \cref{sec:simfree}, and quantify, at a fixed batch size and for absolutely continuous targets, how far minibatch couplings move the induced coupling towards the Brenier map.
\end{problem}
\begin{problem}
Develop a stability theory for flow-map training that explains why some characterisations in \cref{prop:flowmap} train more stably than others, and which parameterisations of $X_{s,t}$ are well conditioned.
\end{problem}
\begin{problem}
For one-step maps trained by consistency or self-distillation losses, give error bounds that depend on intrinsic rather than ambient dimension, for data near curved manifolds, as is now known for the samplers of diffusion models (\cref{sec:theory}).
\end{problem}
\begin{problem}
Identify physical PDEs whose propagators make the velocity field easier to learn, such as the Euclidean Klein--Gordon equation of holographic flows, whose bulk-to-boundary kernel suppresses fine scales first (\cref{sec:holography}), and explain which properties of such kernels are responsible.
\end{problem}
\begin{problem}
Beyond mean-field birth--death sampling and Gaussian location mixtures (\cref{sec:transport,sec:gradflows}), characterise when Fisher--Rao or Wasserstein--Fisher--Rao components provably improve generative modelling of multimodal data known only through samples.
\end{problem}
\end{openprob}

\part{Diffusion models}

\section{Score-based diffusion as time reversal}\label{sec:diffusion}

Diffusion models reach the same goal as flows from the opposite direction. Instead of prescribing a curve of laws and learning its velocity, they prescribe a noising process that destroys the data, and recover the curve by running that process backwards in time. The reverse dynamics require the score of the noised law, and the whole construction rests on the time-reversal theorem of \cref{sec:pathspace} and the integration by parts of \cref{sec:ibp}. This section sets up the forward and reverse dynamics, the likelihood they define and the parameterisations in use, and then the forward processes that depart from the linear Gaussian template.

\subsection{Forward and reverse dynamics}

The construction has two halves, a forward process chosen for its tractability, and a reverse process obtained from it by \cref{thm:reversal}.

\paragraph{Forward processes.}

Diffusion models fix a forward SDE that carries $\mu_{\rm data}$ to an approximately known law and learn to reverse it \citep{sohldickstein2015deep,song2019generative,ho2020denoising,song2021score}. The standard choice is a linear SDE with scalar coefficients,
\begin{equation}\label{eq:forward-linear}
\dd X_t=f(t)X_t\,\dd t+g(t)\,\dd W_t,\qquad X_0\sim\mu_{\rm data},
\end{equation}
whose transition kernel is Gaussian, $X_t\mid X_0\sim\cN(\alpha_tX_0,\sigma_t^2\Id)$ with $\alpha_t=\exp\int_0^tf$ and $\sigma_t^2=\alpha_t^2\int_0^tg(s)^2\alpha_s^{-2}\dd s$. The variance-preserving (VP) SDE takes $f=-\beta/2$, $g=\sqrt\beta$, so that $\alpha_t^2+\sigma_t^2=1$ and $X_t\to\cN(0,\Id)$ exponentially fast; it is the continuous-time limit of DDPM \citep{ho2020denoising}. We assume that the noise schedule satisfies $\beta(t)\ge\beta_{\min}>0$, so that $\alpha_t^2=e^{-\int_0^t\beta}\le e^{-\beta_{\min}t}$. The variance-exploding (VE) SDE takes $f=0$, $g^2=\dd\sigma_t^2/\dd t$; it is the limit of the noise-conditional score networks of \citet{song2019generative,song2020improved}. Up to the rescaling $X_t\mapsto X_t/\alpha_t=X_0+\snr(t)^{-1/2}\epsilon$, everything about such a process is captured by the signal-to-noise ratio $\snr(t)=\alpha_t^2/\sigma_t^2$ or its logarithm $\lambda_t$ \citep{kingma2021variational}.

\paragraph{Reverse dynamics and the probability flow.}

\begin{figure}[t]
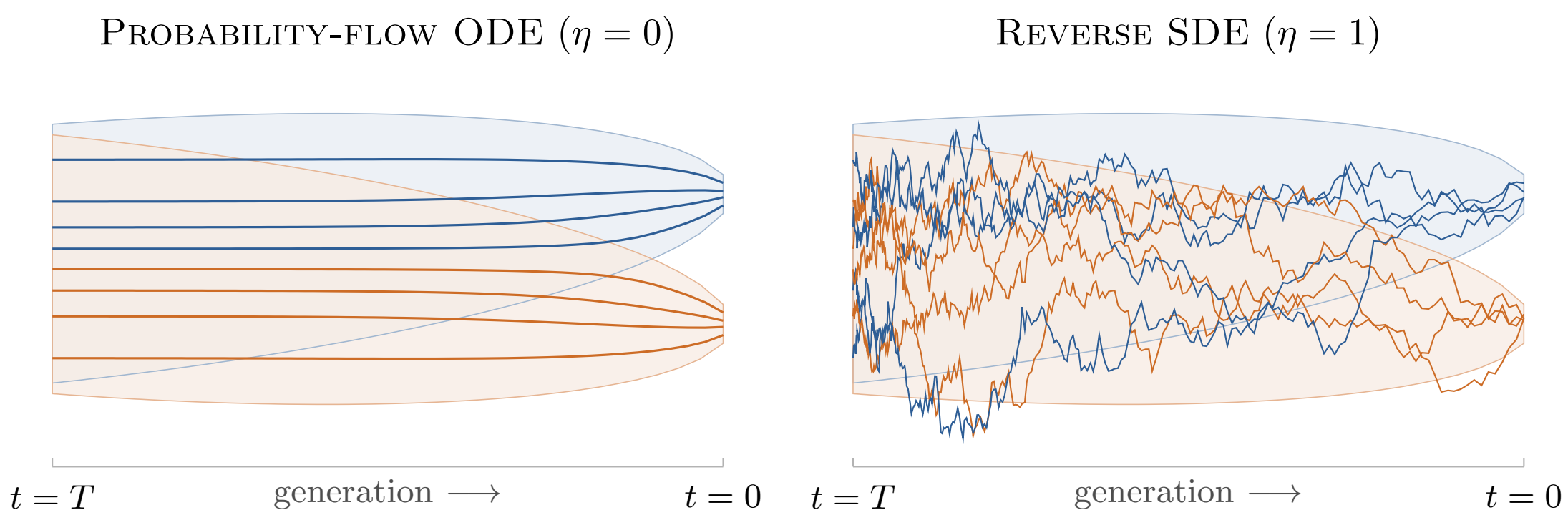

\centering

\caption{Different paths with the same marginals. A one-dimensional variance-preserving diffusion of a two-component Gaussian mixture, run backwards from noise at $t=T$ to data at $t=0$ and plotted against the signal level $\alpha_t$, rescaled to run from $\alpha_T$ at the left edge to $\alpha_0=1$ at the right. The shaded bands mark two standard deviations around the means $\pm\alpha_t$ of the two components of $p_t$, which every member of the $\eta$-family reproduces. In the left panel, the probability-flow ODE ($\eta=0$) follows smooth paths that do not cross. In the right panel, the reverse SDE ($\eta=1$), started from the same points, reaches the same marginals along rough paths; its injected noise contracts accumulated errors, at the price of a larger discretisation error (\cref{sec:samplers}).}
\label{fig:eta-family}
\end{figure}

By \cref{thm:reversal}, the reversal $Y_s=X_{T-s}$ of \eqref{eq:forward-linear} solves
\begin{equation}\label{eq:reverse-linear}
\dd Y_s=\big[-f(T-s)Y_s+g(T-s)^2\,\nabla\log p_{T-s}(Y_s)\big]\dd s+g(T-s)\,\dd\overline W_s,
\end{equation}
and by \eqref{eq:pf-velocity} the probability-flow ODE $\dot x=f(t)x-\frac12g(t)^2\nabla\log p_t(x)$ has the same marginals. More generally, for any $\eta\ge0$, the reverse-time SDE with drift $-f(T-s)y+\frac{1+\eta^2}{2}g(T-s)^2\nabla\log p_{T-s}(y)$ and diffusion $\eta\,g(T-s)$ shares these marginals, interpolating between the ODE at $\eta=0$ and \eqref{eq:reverse-linear} at $\eta=1$. This is the diffusion-model instance of \cref{prop:interpolant-family}, and \cref{fig:eta-family} compares its members at $\eta=0$ and $\eta=1$. Generation replaces $\nabla\log p_t$ by a learned $s_\theta(t,\cdot)$, starts at $\cN(0,\sigma_T^2\Id)$ or $\cN(0,\Id)$, and integrates numerically (\cref{sec:samplers}).

\paragraph{Discrete-time models.}

Before the SDE formulation, diffusion models were built as discrete Markov chains. \citet{sohldickstein2015deep} trained a reverse chain by a variational bound inspired by nonequilibrium thermodynamics. DDPM fixed the forward chain $q(x_k\mid x_{k-1})=\cN(\sqrt{1-\beta_k}\,x_{k-1},\beta_k\Id)$ with a variance schedule $\beta_k\in(0,1)$, parameterised the reverse mean through a noise predictor $\epsilon_\theta$, and observed that a reweighted variational bound reduces to $\E\norm{\epsilon-\epsilon_\theta(x_k,k)}^2$, evaluated at $x_k=\alpha_kx_0+\sigma_k\epsilon$ with $\alpha_k^2=\prod_{j\le k}(1-\beta_j)=1-\sigma_k^2$, for a uniformly drawn step $k$, $x_0\sim\mu_{\rm data}$ and $\epsilon\sim\cN(0,\Id)$ \citep{ho2020denoising}. Learned variances and a cosine schedule improved likelihoods \citep{nichol2021improved}. In parallel, noise-conditional score networks were trained by denoising score matching at many noise levels and sampled by annealed Langevin dynamics \citep{song2019generative,song2020improved}. \citet{song2021score} unified both as discretisations of \eqref{eq:forward-linear}--\eqref{eq:reverse-linear}.

\subsection{Likelihood and parameterisation}

Training regresses a network onto the score, which raises two questions, namely what quantity the regression optimises, and what the network should be asked to output.

\paragraph{Likelihood and the evidence lower bound.}\label{sec:likelihood}

Girsanov's theorem turns the score-matching loss into a likelihood bound.

\begin{proposition}[Relative entropy and score matching \citep{song2021maximum,huang2021variational}]\label{prop:elbo}
Let $\pi_T$ be a law on $\R^d$, such as $\cN(0,\Id)$, and let $q_\theta$ be the law at forward time zero, that is at $s=T$, of the reverse SDE \eqref{eq:reverse-linear} with $\nabla\log p_t$ replaced by $s_\theta$ and started from $Y_0\sim\pi_T$. Then
\begin{equation}\label{eq:likelihood-bound}
\KL(\mu_{\rm data}\,\|\,q_\theta)\le\KL(p_T\,\|\,\pi_T)+\frac12\int_0^Tg(t)^2\,\E\big\|s_\theta(t,X_t)-\nabla\log p_t(X_t)\big\|^2\dd t.
\end{equation}
\end{proposition}
\begin{proof}
Compare the path measures of the true reverse process and the model reverse process. They have the same diffusion coefficient $g$, initial laws $p_T$ and $\pi_T$, and drifts differing by $g^2(\nabla\log p_t-s_\theta)$. Apply \eqref{eq:girsanov-kl} and the data-processing inequality for the time-zero marginal.
\end{proof}
With the likelihood weighting $g(t)^2$, the denoising objective therefore bounds the negative log-likelihood up to $\theta$-independent terms, namely the entropy of $\mu_{\rm data}$, $\KL(p_T\,\|\,\pi_T)$, and minus the constant by which the weighted denoising objective exceeds the weighted score-matching objective. For data with a density, this constant, $\frac12\int_0^Tg(t)^2\big(d/\sigma_t^2-\E\norm{\nabla\log p_t(X_t)}^2\big)\dd t$, and the weighted denoising objective itself diverge at $t=0$, so in practice the time integral is truncated at a small $t_0>0$ \citep{song2021maximum}. Variational diffusion models showed that in continuous time the bound depends on the noise schedule only through the endpoint values of the SNR; they learned the schedule to reduce variance and obtained state-of-the-art likelihoods \citep{kingma2021variational}. \citet{kingma2023understanding} showed further that the common weighted objectives are all weighted integrals of ELBOs over noise levels, and that when the weighting $w(\lambda_t)$ is monotonically non-decreasing in $t$ (non-increasing in $\lambda$), as for $v$-prediction with the cosine schedule and for flow matching with optimal-transport paths (most popular weightings, such as $\epsilon$-prediction with the cosine schedule and EDM, are not), the objective equals the ELBO of data augmented by Gaussian noise. The velocity $v$ targeted by $v$-prediction and the EDM design are described in the next paragraph. Exact likelihoods follow from the probability-flow ODE by \eqref{eq:inst-cov}, and training ODE likelihoods requires matching higher-order score information \citep{lu2022maximum,zheng2023improved}. This concerns score-matching objectives; a continuous normalising flow can instead maximise the exact likelihood directly (\cref{sec:cnf}), at the cost of simulating the ODE during training.

\paragraph{Parameterisations and preconditioning.}

The score, the noise $\epsilon$, the clean data $x_0$ and the velocity $v=\alpha_t\epsilon-\sigma_tx_0$ are affinely related given $x_t$, since $\nabla\log p_t(x_t)=-\E[\epsilon\mid x_t]/\sigma_t$ and $x_t=\alpha_t\E[x_0\mid x_t]+\sigma_t\E[\epsilon\mid x_t]$ (\cref{prop:tweedie}). The velocity is the target of $v$-prediction; on a variance-preserving path, with $\alpha_t=\cos\phi_t$ and $\sigma_t=\sin\phi_t$, it is the rate of change $\dd x_t/\dd\phi_t$ of $x_t=\alpha_tx_0+\sigma_t\epsilon$, which explains the name. Choosing which one the network outputs is a choice of preconditioning and of implicit loss weighting, and \cref{app:dictionary} tabulates the conversions. \citet{karras2022elucidating} separated the design choices, namely the noise schedule, the network preconditioning by input, output and skip scalings, the training noise distribution and loss weighting, and a second-order Heun sampler, and showed that each can be tuned independently. \citet{karras2024analyzing} later revisited training dynamics and normalisation at scale. The geometry of the SNR curve matters for training in other ways too. Weighting by a clipped SNR helps \citep{hang2023efficient}, as does enforcing zero terminal SNR so that the model actually sees pure noise at $t=T$ \citep{lin2024common}. \citet{strasman2024analysis} analyse the effect of the noise schedule on the error bounds of \cref{sec:theory}.

\subsection{Other forward processes}

The linear Gaussian forward process is convenient but not necessary, and it has been relaxed by adding momentum, by allowing nonlinear dynamics, and by replacing noise with deterministic degradation.

\paragraph{Kinetic and augmented forward processes.}\label{sec:kinetic}

Critically damped Langevin diffusion (CLD) augments the state with a velocity $V$ and runs
\begin{align*}
\dd X_t&=M^{-1}V_t\,\beta\,\dd t,\\
\dd V_t&=-\beta X_t\,\dd t-\Gamma\beta M^{-1}V_t\,\dd t+\sqrt{2\Gamma\beta}\,\dd W_t,
\end{align*}
with mass $M>0$, friction coefficient $\Gamma>0$, time scale $\beta>0$ and noise only in the velocity; critical damping means $\Gamma^2=4M$ \citep{dockhorn2022score}. Write $p_t$ for the joint density of $(X_t,V_t)$. The diffusion matrix is degenerate, yet the pair $(X_t,V_t)$ has a smooth positive density, because the drift couples the components and H\"ormander's condition holds (\cref{sec:malliavin}). The time reversal involves only the velocity score $\nabla_v\log p_t=\nabla_v\log p_t(v\mid x)$, the score of the conditional law of the velocity given the position, which the authors argue is easier to learn than the data score. They derive a hybrid score-matching objective and a symmetric splitting sampler of Strang type. Hybrid score matching also has a likelihood interpretation. From the Radon--Nikodym derivative between forward and backward processes whose noise acts only in some coordinates, \citet{blessing2025underdamped} derive the evidence lower bound for such processes and prove that underdamped score matching maximises it, which makes rigorous a connection drawn by the authors of CLD. \citet{pandey2023generative} give a recipe for general augmented linear diffusions. Kinetic processes are where hypoellipticity, the original subject of Malliavin calculus, enters generative modelling directly.

\paragraph{Nonlinear and state-dependent forward processes.}\label{sec:nonlinear-forward}

The Gaussian transition kernel of \eqref{eq:forward-linear} is what makes denoising score matching computable. Learned or nonlinear forward processes give up some of that convenience for flexibility, and the designs proposed so far each keep some tractable structure to avoid computing the score of a fully nonlinear SDE. Diffusion normalising flow learns a forward drift jointly with the reverse model \citep{zhang2021diffusion}. Implicit nonlinear diffusion models run a linear diffusion in the latent space of a normalising flow \citep{kim2022maximum}, and neural flow diffusion models learn the forward process through a time-dependent invertible map \citep{bartosh2024neural}. Multivariate linear diffusions with learned coefficients are explored by \citet{singhal2023where}. The Malliavin formulas of \cref{thm:malliavin-score} remove that constraint in principle. For any forward SDE whose coefficients are smooth with bounded derivatives and whose Malliavin covariance is nondegenerate, the score is a conditional expectation of a quantity simulated together with the process and its first and second variations; \citet{mirafzali2025malliavin} treat state-independent diffusion coefficients and \citet{mirafzali2025approach} the general case. The divergence--kernel method gives another forward-simulable score estimator for multiplicative noise and, by differentiating marginal densities with respect to parameters, allows forward-only training of parametric SDE generative models that fit the data law directly, demonstrated so far on low-dimensional examples (up to 20 dimensions, at most 21 parameters) \citep{ni2025scores,ni2025divergence}. \Cref{sec:score-estimation} compares these estimators; the divergence--kernel formulas are rigorous for time-discretised dynamics and formal in the continuous-time limit.

\paragraph{Deterministic and latent corruptions.}

The forward process can also be a deterministic degradation with little or no noise. Inverse heat dissipation runs the heat equation on the image plane, so that generation proceeds from coarse to fine scales \citep{rissanen2023generative}, and blurring diffusion models combine blurring with noise in frequency space \citep{hoogeboom2023blurring}. Cold diffusion goes further and inverts arbitrary degradations without noise \citep{bansal2023cold}, while soft diffusion develops score matching for linear corruption operators \citep{daras2023soft}. Latent diffusion moves the whole process into the latent space of a pretrained autoencoder \citep{rombach2022high}; the dynamics are unchanged and we do not discuss autoencoders further. These variants connect to the renormalisation-group view of \cref{sec:fieldtheory}, where coarse-graining is the forward process. Non-Gaussian driving noise is treated in \cref{sec:levy}.

\section{Score estimation}\label{sec:score-estimation}

Everything in the reverse dynamics hinges on the score, which is learned by regression onto a random target whose conditional mean is the score. This section asks how good that regression can be. It considers how noisy the standard target is and how it can be improved, which targets exist when the forward process is nonlinear, and what is known about the score as a statistical object.

\subsection{Objectives and their variance}

With weights $w(t)$, the denoising objective is
\begin{align*}
\mathcal J(\theta)=\int_0^Tw(t)\,\E\Big\|s_\theta(t,X_t)+\frac{X_t-\alpha_tX_0}{\sigma_t^2}\Big\|^2\dd t,
\end{align*}
a regression onto a random target whose conditional mean is the score (\cref{prop:sm}). When $\mu_{\rm data}$ has a density, the target has conditional variance $d/\sigma_t^2-\E\norm{\nabla\log p_t(X_t)}^2\approx d/\sigma_t^2$ (summed over coordinates and averaged over $X_t$), which explodes as $t\to0$ and is the main statistical difficulty of score estimation near the data. For data on a $k$-dimensional set the conditional variance is of order $k/\sigma_t^2$, and for the empirical measure it vanishes as $t\to0$. Two remedies change the random target while keeping its conditional mean, exactly or up to a vanishing bias, and a third option avoids the target altogether.

The first remedy averages over the posterior. The stable target field replaces the single-sample target by a self-normalised importance-weighted average of conditional scores over a reference batch of clean data, an estimate of the posterior mean $\E[X_0\mid X_t]$ with much lower variance at intermediate noise levels; the self-normalised average trades a small bias, which vanishes as the reference batch grows, for this variance reduction \citep{xu2023stable}. The second integrates by parts on the data rather than on the noise. Target score matching uses the identity
\begin{align*}
\nabla\log p_t(x)=\alpha_t^{-1}\,\E\big[\nabla\log p_0(X_0)\,\big|\,X_t=x\big],
\end{align*}
which can be used when the data score is known and is accurate exactly where denoising score matching is noisy \citep{debortoli2024target}. Because the denoising and target score-matching targets are both unbiased, any time-dependent convex combination of them is also unbiased (\cref{tab:ibp}). The third option abandons the Gaussian kernel. The trace form of \cref{prop:sm}(i), estimated by random projections, needs no transition density \citep{song2019sliced}. It has no random target, and for data with a smooth density the variance of its integrand stays bounded as $t\to0$. But near a low-dimensional set, or for the empirical measure, $\frac12\norm{s_\theta}^2$ and $\tr\nabla s_\theta$ themselves grow like $\sigma_t^{-2}$, and the variance returns at order $\sigma_t^{-4}$.

\subsection{Scores of nonlinear SDEs via Malliavin calculus}\label{sec:malliavin-scores}

For a general forward SDE \eqref{eq:sde} without closed-form kernels, \cref{thm:malliavin-score} yields an explicit target. We state the result of \citet{mirafzali2025approach} in a simplified form.

\begin{theorem}[Bismut-type score formula for nonlinear SDEs \citep{mirafzali2025approach}]\label{thm:bismut-score}
Let $X$ solve \eqref{eq:sde} on $\R^d$ from a fixed point $X_0=x$, with $C^2$ coefficients having bounded derivatives, let $Y$ and $Y^{(2)}$ be its first and second variation processes, and suppose the Malliavin covariance $\gammaM_{X_T}$ in \eqref{eq:malliavin-cov} is invertible with sufficient inverse moments. Then for each $k$,
\begin{align*}
\partial_{y_k}\log p_T(y)&=-\E\big[\Sk(u_k)\,\big|\,X_T=y\big],\\
u_k(t)&=\sum_j(\gammaM_{X_T}^{-1})_{kj}\,\D_tX_T^j,
\end{align*}
and $\Sk(u_k)$ equals the It\^o integral $\int_0^Tu_t(x')\cdot\dd W_t$ of the adapted field $u_t(x')=x'^\top Y_t^{-1}\sigma(t,X_t)$, computed for fixed $x'$ and then evaluated at $x'=Y_T^\top\gammaM_{X_T}^{-1}e_k$, the point at which $u_t(x')=u_k(t)$, minus an explicit correction built from $Y$, $Y^{-1}$, $Y^{(2)}$, $\partial_x\sigma$ and $\gammaM_{X_T}^{-1}$. For state-independent $\sigma(t)$ the correction simplifies; if in addition the drift is linear, it vanishes and the formula reduces to the Gaussian score. With state-dependent noise the correction survives even for linear drift, as \cref{ex:gbm} shows.
\end{theorem}

\begin{example}[Geometric Brownian motion]\label{ex:gbm}
A simple example with state-dependent noise is geometric Brownian motion, for which the answer is known and the role of each term can be seen. Let $\dd X_t=\mu X_t\dd t+\sigma X_t\dd W_t$ with $X_0=x>0$, so that $X_T=x\exp\big((\mu-\frac{\sigma^2}{2})T+\sigma W_T\big)$. By \cref{prop:sde-derivative}, $\D_sX_T=\sigma X_T$ for $s\le T$, hence $\norm{\D X_T}_H^2=\sigma^2TX_T^2$ and the covering field is the constant-in-time but anticipating process
\begin{align*}
u(s)=\frac{\D_sX_T}{\norm{\D X_T}_H^2}=\frac{1}{\sigma TX_T}=:F_T.
\end{align*}
By \eqref{eq:delta-elementary} with $h=\one_{[0,T]}$, $\Sk(u)=F_TW_T-\inner{\D F_T}{\one_{[0,T]}}_H$. Since $\D_sF_T=-\D_sX_T/(\sigma TX_T^2)=-1/(TX_T)$, the correction equals $-1/X_T$ and
\begin{align*}
\Sk(u)=\frac{W_T}{\sigma TX_T}+\frac{1}{X_T}.
\end{align*}
Given $X_T=y$ the Brownian endpoint is determined, $W_T=\big(\log(y/x)-(\mu-\frac{\sigma^2}{2})T\big)/\sigma$, so \cref{prop:ibp-laws} yields
\begin{align*}
\partial_y\log p_T(y)=-\frac{\log(y/x)-(\mu-\frac{\sigma^2}{2})T}{\sigma^2Ty}-\frac1y,
\end{align*}
which is exactly the score of the log-normal law of $X_T$. The term $-1/y$ comes entirely from the Skorokhod correction. An It\^o-style formula that ignored the anticipation of $u$ would miss it. This is the one-dimensional shadow of the correction term in \cref{thm:bismut-score}. Here the drift and diffusion coefficients are linear in $x$, so the second variation $Y^{(2)}$ vanishes identically, and the correction comes entirely from the state dependence of the noise, $\partial_x(\sigma x)=\sigma\neq0$, through $\D_sY_T=\sigma Y_T$.
\end{example}

The formula of \cref{thm:bismut-score} turns score estimation for nonlinear diffusions into Monte Carlo simulation of an augmented SDE system $(X,Y,Y^{-1},Y^{(2)})$, followed by regression of $-\Sk(u_k)$ on $X_T$; \citet{mirafzali2025malliavin} report experiments in the state-independent case. Three related developments obtain scores by integration by parts in other settings. \citet{pidstrigach2025conditioning} use the same Wiener-space integration by parts to derive a generalised Tweedie formula for conditional scores of nonlinear SDEs, to which \cref{sec:conditioning} returns. \citet{ni2025scores} instead integrates by parts on the transition kernel of each step and derives a divergence--kernel formula, rigorous for time-discretised dynamics and obtained formally in the continuous-time limit, in which the score is the conditional expectation of a forward covector process with a free backward-adapted parameter; that parameter shifts part of the divergence term onto the transition kernel and damps the exponential growth of the pure divergence formula. \citet{ni2025divergence} extends the divergence--kernel formula to parameter derivatives of marginal and stationary densities, which makes gradient-based fitting of multiplicative-noise SDEs possible.

\begin{remark}[Choice of covering field]
All covering fields give unbiased conditional estimators, and so do their convex combinations with those targets of \cref{tab:ibp} that estimate the same score. The minimal-energy field of \cref{thm:bismut-score} is canonical but not obviously the minimal-variance one; the free parameter in the divergence--kernel method plays an analogous role. A theory of variance-optimal covering fields, possibly learned jointly with the score network, is open.
\end{remark}

\subsection{Consistency, the empirical score and sample complexity}

We turn from the estimators to the object being estimated. A score field is far from arbitrary, its exact empirical version is known in closed form, and the rate at which it can be learned from samples has been characterised for some classes of data laws.

\paragraph{Consistency with the Fokker--Planck equation.}

A score field $s_t=\nabla\log p_t$ is constrained by the dynamics that produced it. Since $p_t$ solves \eqref{eq:fpe}, $s_t$ solves a nonlinear PDE obtained by differentiating $\log p_t$. Several works exploit this constraint. Enforcing the score Fokker--Planck equation as a regulariser improves likelihoods \citep{lai2023fp}, and \citet{deveney2023closing} analyse the mismatch between the ODE and SDE samplers when the learned score violates it. Consistent diffusion models enforce a different constraint, the martingale property of the denoiser along reverse trajectories \citep{daras2023consistent}. Self-consistency of the Fokker--Planck equation also yields a solver for that equation. A velocity field is learned so that it equals the drift minus the score of the density it transports, and the residual bounds the Wasserstein error of the resulting solution \citep{shen2022self}. Score-based transport similarly solves Fokker--Planck equations by learning the probability-flow velocity of a given SDE \citep{boffi2023probability}.

\paragraph{The empirical score.}\label{sec:empirical-score}

Minimising the denoising objective over all functions, with the empirical measure $\hat\mu_n=\frac1n\sum_i\delta_{x^{(i)}}$ in place of $\mu_{\rm data}$, yields the score of $\hat\mu_n$ smoothed by the forward kernel.

\begin{proposition}[Optimal empirical denoiser]\label{prop:empirical}
For $\hat\mu_n$ and the Gaussian kernel of \eqref{eq:forward-linear},
\begin{align*}
\E[X_0\mid X_t=x]&=\sum_{i=1}^nw_i(x)\,x^{(i)},\\
w_i(x)&=\frac{\exp\big(-\norm{x-\alpha_tx^{(i)}}^2/(2\sigma_t^2)\big)}{\sum_j\exp\big(-\norm{x-\alpha_tx^{(j)}}^2/(2\sigma_t^2)\big)},
\end{align*}
and the empirical score is $(\alpha_t\sum_iw_ix^{(i)}-x)/\sigma_t^2$.
\end{proposition}

The reverse dynamics driven by this exact empirical score concentrate on the training points as $t\to0$, so a perfectly trained model memorises. That trained networks generalise is therefore a statement about what they fail to fit, which is the subject of \cref{sec:generalisation}. Smoothing the empirical score gives training-free closed-form diffusion models that already produce some novel samples \citep{scarvelis2023closed}.

\paragraph{Statistical complexity of score estimation.}

Score matching is statistically efficient when the data law satisfies good isoperimetric inequalities and can be badly inefficient otherwise \citep{koehler2023statistical}. It can nonetheless be preferable to maximum likelihood on computational grounds. For exponential families of fixed-degree polynomials, zeroth- and first-order oracles for the likelihood are NP-hard, whereas score matching is efficient with polynomial sample complexity \citep{pabbaraju2023provable}. Generalised score matching built on fast-mixing diffusions recovers efficiency \citep{qin2023fit}. These results concern parametric families, in which score matching competes with maximum likelihood; for diffusion models the relevant question is nonparametric. For a sub-Gaussian data law $\rho$ with Lipschitz score $s=\nabla\log\rho$, a regularised Gaussian-kernel estimator $\hat s$ built from $n$ samples by empirical-Bayes smoothing attains the minimax rate $\tilde\Theta(n^{-2/(d+4)})$, up to logarithmic factors, for the squared error $\norm{\hat s-s}^2_{L^2(\rho)}$. This is a curse of dimensionality for score estimation, with consequences for the sample complexity of score-based generative models \citep{wibisono2024optimal}.

\section{Samplers as numerical stochastic analysis}\label{sec:samplers}

Once the score is learned, generation is the numerical solution of a differential equation, and the cost of a sample is the number of network evaluations the solver needs. Diffusion equations have a structure that general-purpose solvers ignore. They have a linear part that can be integrated exactly and a nonlinearity supplied by the network, and they leave a choice, at sampling time, of how much noise to inject. The samplers below exploit this structure, first for the deterministic flow, then for the stochastic one, and finally across time.

\paragraph{Ancestral sampling, Euler--Maruyama and DDIM.}

Ancestral sampling of DDPM is an Euler--Maruyama-type discretisation of \eqref{eq:reverse-linear}. DDIM introduced a family of non-Markovian forward processes with the same marginals and a deterministic sampler which, in suitable variables, is a first-order discretisation of the probability-flow ODE \citep{song2021denoising}; its stochasticity parameter interpolates between deterministic and ancestral sampling, mirroring the $\eta$-family of \cref{sec:diffusion}. General SDE numerics \citep{kloeden1992numerical} apply, but the structure of diffusion equations allows faster schemes.

\paragraph{Exponential integrators.}

In the noise parameterisation the probability-flow ODE is semilinear, $\dot x=f(t)x+\frac{g(t)^2}{2\sigma_t}\epsilon_\theta(x,t)$. Changing variables to the half log-SNR $\bar\lambda_t=\frac12\lambda_t=\log(\alpha_t/\sigma_t)$ and solving the linear part exactly over a step from $s$ to $t<s$ gives
\begin{equation}\label{eq:dpm}
x_t=\frac{\alpha_t}{\alpha_s}x_s-\alpha_t\int_{\bar\lambda_s}^{\bar\lambda_t}e^{-\bar\lambda}\,\hat\epsilon_\theta(x_{\bar\lambda},\bar\lambda)\,\dd\bar\lambda,
\end{equation}
where $\hat\epsilon_\theta(x_{\bar\lambda},\bar\lambda)=\epsilon_\theta(x_{t(\bar\lambda)},t(\bar\lambda))$ is the network evaluated along the solution and $t(\bar\lambda)$ inverts $t\mapsto\bar\lambda_t$, so that only the network term needs approximation \citep{lu2022dpm}. The solvers in use differ in how they approximate that integral. Taylor-expanding $\hat\epsilon_\theta$ in $\bar\lambda$ yields the DPM-Solver family, whose data-prediction and multistep variants are more robust under guidance \citep{lu2022dpmpp}. Polynomial extrapolation of past network evaluations gives DEIS \citep{zhang2023fast}, empirical model statistics improve the expansion \citep{zheng2023dpm}, and predictor--corrector variants unify these schemes \citep{zhao2023unipc}. Pseudo-numerical methods \citep{liu2022pseudo} and the second-order Heun scheme of \citet{karras2022elucidating} complete the deterministic toolbox.

\paragraph{Stochastic solvers, correctors and restarts.}

The same ideas carry over to the reverse SDE. Exponential integrators extend to it \citep{gonzalez2023seeds} and to stochastic multistep schemes \citep{xue2023sa}, and adaptive step sizes control the error of stochastic solvers \citep{jolicoeur2021gotta}. Noise can also be used deliberately to correct errors. Predictor--corrector sampling alternates a reverse-SDE step with Langevin corrector steps targeting the current marginal \citep{song2021score}, and restart sampling alternates large forward noise injections with deterministic backward ODE segments, so that the noise contracts accumulated errors while the ODE keeps the discretisation error small \citep{xu2023restart}.

\paragraph{Parallel-in-time sampling.}

Sequential network evaluations dominate the cost, which suggests parallelising across time. Picard iterations over a sliding window of time steps, refined in parallel until convergence, do so \citep{shih2023parallel}, and with a careful analysis parallel sampling achieves sub-linear time complexity in the dimension \citep{chen2024accelerating}. Randomised-midpoint discretisations give improved sequential and parallel guarantees \citep{gupta2024faster}. The Poisson midpoint method, a variant that emulates many small Euler--Maruyama steps with one large step, gives a provable quadratic speed-up for Langevin Monte Carlo under weak assumptions and, applied to diffusion models, matches 1000-step DDPM sampling with 50--80 network evaluations \citep{kandasamy2024poisson}.

\paragraph{Deterministic or stochastic sampling.}

Deterministic samplers have smaller discretisation error and allow large steps, but errors in the score are transported without damping. Stochastic samplers inject noise that contracts the total-variation distance between the true and perturbed reverse processes, at least for data of bounded support as in the restart analysis above, but they incur a larger discretisation error. The trade-off is explicit in the interpolant family of \cref{prop:interpolant-family}, whose relative-entropy bound weighs score errors by the noise level, and in the theory of \cref{sec:theory}, where SDE samplers enjoy Girsanov-based guarantees under $L^2$ score accuracy alone while ODE samplers need control of the Jacobian of the score.

\section{Convergence and statistical theory}\label{sec:theory}

Why should a sampler driven by an imperfect score produce samples close to the data? For stochastic samplers the theory answers this in a uniform way. Girsanov's theorem converts the question into one about drifts, the error splits into three terms, and the work lies in bounding the hardest of them under assumptions weak enough to cover realistic data. Deterministic samplers need other arguments. This section follows that programme from the first polynomial guarantees to sharp dependence on dimension, deterministic samplers, low-dimensional structure and Wasserstein distances, and ends with end-to-end statistical rates.

\subsection{Error decomposition and polynomial guarantees}

\paragraph{The error decomposition.}

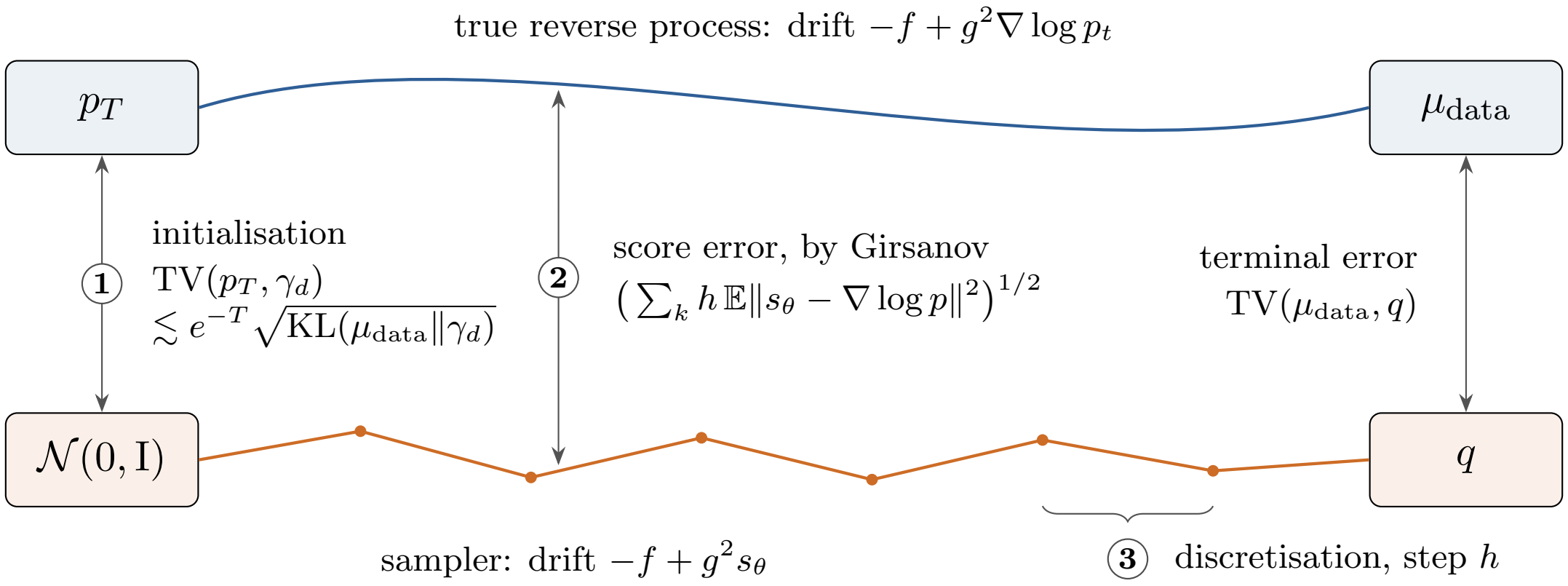
\begin{figure}[t]
\centering
\begin{tikzpicture}[>=Stealth,
  endbox/.style={draw,rounded corners=3pt,inner sep=4pt,minimum width=1.75cm,minimum height=0.85cm,align=center,fill=#1},
  num/.style={circle,draw=black!70,fill=white,inner sep=1.1pt,font=\scriptsize\bfseries},
  lab/.style={font=\footnotesize,align=left}]
\node[endbox=figblue!9] (pT) at (0,1.6) {$p_T$};
\node[endbox=figblue!9] (mu) at (12.4,1.6) {$\mu_{\rm data}$};
\node[endbox=figorange!10] (g) at (0,-1.6) {$\cN(0,\Id)$};
\node[endbox=figorange!10] (q) at (12.4,-1.6) {$q$};
\draw[thick,figblue,name path=truepath] (pT.east) .. controls (3.7,2.45) and (8.4,0.85) .. (mu.west);
\node[lab,anchor=south] at (6.2,2.05) {true reverse process: drift $-f+g^2\nabla\log p_t$};
\draw[thick,figorange,name path=sampler] (g.east) -- (2.35,-1.34) -- (3.9,-1.76) -- (5.45,-1.4) -- (7.0,-1.78) -- (8.55,-1.42) -- (10.1,-1.7) -- (q.west);
\foreach \p in {(2.35,-1.34),(3.9,-1.76),(5.45,-1.4),(7.0,-1.78),(8.55,-1.42),(10.1,-1.7)} \fill[figorange] \p circle (1.5pt);
\draw[<->,gray!65!black] (pT.south) -- (g.north);
\node[num] at (0,0) {1};
\node[lab,anchor=west] at (0.32,0) {initialisation\\[1pt]$\TV(p_T,\gamma_d)$\\[1pt]$\lesssim e^{-T}\sqrt{\KL(\mu_{\rm data}\|\gamma_d)}$};
\path[name path=probe] (4.15,3) -- (4.15,-3);
\path[name intersections={of=truepath and probe,by=A}];
\path[name intersections={of=sampler and probe,by=B}];
\draw[<->,gray!65!black,shorten <=1.5pt,shorten >=1.5pt] (A) -- (B);
\node[num] (n2) at ($(A)!0.5!(B)$) {2};
\node[lab,anchor=west] at ([xshift=5pt]n2.east) {score error, by Girsanov\\[1pt]$\big(\sum_kh\,\E\norm{s_\theta-\nabla\log p}^2\big)^{1/2}$};
\draw[decorate,decoration={brace,mirror,amplitude=4pt},gray!65!black] (8.55,-2.02) -- (10.1,-2.02);
\node[num] (n3) at (9.325,-2.5) {3};
\node[lab,anchor=west] at ([xshift=3pt]n3.east) {discretisation, step $h$};
\node[lab,anchor=center] at (4.3,-2.5) {sampler: drift $-f+g^2s_\theta$};
\draw[<->,gray!65!black] (mu.south) -- (q.north);
\node[lab,anchor=east,align=right] at (12.1,0) {terminal error\\[1pt]$\TV(\mu_{\rm data},q)$};
\end{tikzpicture}
\caption{Where the three terms of \eqref{eq:error-decomp} come from. The sampler (bottom) differs from the true reverse process (top) in three ways. It starts from $\cN(0,\Id)$ rather than from $p_T$ (1), it uses the learned score $s_\theta$ rather than $\nabla\log p_t$ (2), and it moves in steps of size $h$ (3). Girsanov's theorem compares the two path measures, and the data-processing inequality transfers the comparison to the terminal laws $\mu_{\rm data}$ and $q$.}
\label{fig:error}
\end{figure}

Let $q$ be the output law of a discretised reverse SDE with learned score, initialised at $\cN(0,\Id)$ at time $T$ and run with step size $h$ for $N=T/h$ steps. The data-processing inequality and \eqref{eq:girsanov-kl} decompose the error into three terms, traced in \cref{fig:error},
\begin{equation}\label{eq:error-decomp}
\begin{aligned}
\TV(\mu_{\rm data},q)\lesssim{}&\underbrace{\TV(p_T,\gamma_d)}_{\text{initialisation}}+\underbrace{\Big(\sum_kh\,\E\norm{s_\theta-\nabla\log p}^2\Big)^{1/2}}_{\text{score error}}\\
&+\underbrace{\text{discretisation}}_{\text{depends on }h,\,d,\text{ regularity}},
\end{aligned}
\end{equation}
where the first term decays like $e^{-T}$ for the Ornstein--Uhlenbeck process $\dd X_t=-X_t\,\dd t+\sqrt2\,\dd W_t$ (the VP process with $\beta\equiv2$; in general like $e^{-\frac12\int_0^T\beta}$), and the score error sums $h\,\E\norm{s_\theta(t_k,X_{t_k})-\nabla\log p_{t_k}(X_{t_k})}^2$ over the grid times $t_k=kh$, $k=1,\dots,N$, with $X$ the forward process. The difficulty lies in the last term under weak assumptions on $\mu_{\rm data}$, and in the use of Girsanov's theorem when Novikov's condition fails, which \citet{chen2023sampling} circumvent by an approximation argument.

\paragraph{Polynomial guarantees under $L^2$ score accuracy.}

Early guarantees assumed dissipativity or functional inequalities for the data law \citep{block2020generative,lee2022convergence}. \citet{lee2023convergence} and \citet{chen2023sampling} removed such assumptions, and the latter shows, informally, the following.

\begin{theorem}[Convergence of denoising diffusion, informally \citep{chen2023sampling}]\label{thm:ddpm-convergence}
Suppose the scores $\nabla\log p_t$ are $L$-Lipschitz for all $t$, $\mu_{\rm data}$ has finite second moment, and the learned score has $L^2(p_t)$ error at most $\varepsilon_{\rm score}$ at the discretisation times. Then the DDPM-type sampler for the Ornstein--Uhlenbeck forward process $\dd X_t=-X_t\,\dd t+\sqrt2\,\dd W_t$ with step size $h\lesssim1/L$ satisfies
\begin{align*}
\TV(\mu_{\rm data},q)\lesssim{}&e^{-T}\sqrt{\KL(\mu_{\rm data}\,\|\,\gamma_d)}\\
&+\big(L\sqrt{dh}+Lm_2h\big)\sqrt T+\varepsilon_{\rm score}\sqrt T,
\end{align*}
where $m_2^2$ is the second moment of $\mu_{\rm data}$. When $m_2\lesssim d$ and $\KL(\mu_{\rm data}\,\|\,\gamma_d)\le\mathrm{poly}(d)$, $\tilde O(L^2d/\varepsilon^2)$ steps suffice for error $\varepsilon$, provided $\varepsilon_{\rm score}\lesssim\varepsilon/\sqrt T$ with $T\asymp\log(\KL(\mu_{\rm data}\,\|\,\gamma_d)/\varepsilon)$, that is $\varepsilon_{\rm score}\le\tilde O(\varepsilon)$.
\end{theorem}

The Lipschitz assumption on the scores can be avoided by early stopping, that is, by approximating the slightly noised law $p_\delta$ instead of $\mu_{\rm data}$ \citep{chen2023improved}.

\subsection{Refined guarantees}

The polynomial guarantee leaves open how the number of steps should scale with dimension, whether deterministic samplers enjoy similar guarantees, whether low-dimensional structure in the data is exploited automatically, and what can be said in transport distances. Each of these questions has since been addressed.

\paragraph{Dimension dependence.}

\citet{benton2024nearly} use stochastic localisation (\cref{sec:localisation}) to control the expected change of the score along the reverse process. For any law with finite second moment and early stopping at $\delta$, they obtain a number of steps of order $d\log^2(1/\delta)/\varepsilon^2$, up to logarithmic factors, for relative-entropy accuracy $\varepsilon^2$; the dependence on dimension is linear and no smoothness is assumed. \citet{jain2025sharp} improve the dependence on the accuracy to order $d\log^{3/2}(1/\delta)/\varepsilon$ steps, still without smoothness assumptions, by alternating probability-flow ODE steps, which control an error of Wasserstein type, with short forward noising steps that convert it into a bound in relative entropy.

Changing the sampler removes the polynomial dependence on the accuracy. For samplers with Gaussian transitions, DDPM among them, the error of a single step is bounded below by a data-dependent term, governed by the Frobenius norm of the Jacobian of the posterior mean, that no choice of mean can remove (the exponential-integrator step is almost the KL-optimal Gaussian), so the transitions themselves must be corrected. Each backward transition is a Gaussian tilt of the marginal at the next, less noisy, time. First-order rejection sampling corrects the exponential-integrator step by an accept--reject mechanism simulated from score evaluations alone, through randomised estimates of log-density increments, which are integrals of the score. With score error of order $\varepsilon$ in $L^2$, up to logarithms, and only a finite second moment, this reaches accuracy $\varepsilon$ for the data law in $\widetilde O(d_\star\,\mathrm{polylog}(1/\varepsilon))$ steps, where $d_\star\le d$ is an intrinsic dimension defined through covering numbers of the support. Accuracy is measured in a bounded-Lipschitz metric, through a relative-entropy bound for the slightly noised law. The same subroutine gives the first sampler for log-concave targets that uses only gradients of the log-density and depends polylogarithmically on the accuracy. The results are theoretical, and the sampler has not yet been implemented \citep{chen2026high}. Relative-entropy bounds under finite Fisher information are given by \citet{conforti2025kl}. A complementary line obtains a total-variation error of order $d/N$, up to logarithms, plus the score error for $N$ steps under minimal assumptions \citep{li2024od}, together with provably accelerated samplers \citep{li2024accelerating}. These total-variation bounds, and those for DDPM-type and deterministic samplers below, compare the output with the slightly noised law of the first forward step, whose noise variance is polynomially small in $N$. For data concentrated on a Lebesgue-null set the total-variation distance to $\mu_{\rm data}$ itself equals $1$.

For special data laws the dimension can disappear altogether. Dimension-free guarantees hold for mixtures of isotropic Gaussians, and for data close to them in total variation up to an approximation term \citep{li2025dimension}. For data laws supported on a countable set, an information-theoretic analysis free of geometric assumptions bounds the time-discretisation error in relative entropy by order $\mathrm{H}(\mu_{\rm data})^2/N$ for $N$ steps, where $\mathrm{H}(\mu_{\rm data})$ is the Shannon entropy of the data law, with the score-estimation error entering as a separate term \citep{aghapour2026entropy}. The bound holds up to endpoint factors, and the analysis also assumes a finite second moment and sub-exponential information content. Under $L^2$ score accuracy, then, the number of steps needed by stochastic samplers grows at most nearly linearly with dimension, and for structured data laws it need not grow at all.

\paragraph{Deterministic samplers.}

Probability-flow ODEs are harder to analyse, because noise no longer regularises the error. \citet{chen2023probability} obtained the first polynomial guarantees for ODE sampling by interleaving Langevin corrector steps, with better dimension dependence when the corrector is underdamped. \citet{li2024towards} proved rates of order $1/N$ for the ODE sampler, and \citet{li2024sharp} sharpened them to nearly linear dependence on $d$ in total variation, for absolutely continuous data supported in a ball of radius polynomial in $N$, assuming $\ell_2$ accuracy of the score and a bound on the average operator-norm error of its Jacobian. An example in that paper shows that $\ell_2$ score accuracy alone cannot give total-variation guarantees for deterministic samplers, so some additional control, such as Jacobian accuracy, is needed. \citet{cai2025minimax} show that an ODE sampler with a smooth regularised score estimator, which controls both the $L^2$ score error and the mean Jacobian error, attains the minimax rate in total variation, up to logarithms, for sub-Gaussian laws with $\beta$-H\"older densities, $\beta\le2$. The same control lets deterministic samplers adapt to low-dimensional structure. With suitably designed coefficients, the discretised probability-flow ODE reaches total-variation error of order $(k+\log d)\log^3 N/N$ after $N$ steps, plus the score and Jacobian errors, where $k$ is an intrinsic dimension of the support defined through covering numbers \citep{tang2025adaptivity}. Rectified flow with a U-shaped time grid, geometric towards both ends, needs $O(k/\varepsilon)$ steps, up to logarithms, for accuracy $\varepsilon$ in total variation, given accurate velocities, Jacobians, divergences and gradients of divergences; its stochastic counterpart, which stochastic localisation identifies with DDPM, attains the same rate from $L^2$ score accuracy alone \citep{roy2026low}. Deterministic samplers therefore need control of the Jacobian of the score as well as of its values; with that control, their guarantees, including adaptation to intrinsic dimension, are comparable to those of stochastic samplers, and \cref{prob:intrinsic-ode} asks whether intrinsic-dimension guarantees can make do with weaker control.

\paragraph{Low-dimensional structure.}

When $\mu_{\rm data}$ concentrates near a $k$-dimensional set, rates should depend on $k$ rather than on $d$. \citet{debortoli2022convergence} gave the first Wasserstein bounds for manifold-supported data, and \citet{chen2023score} treat linear subspaces. \citet{tang2024adaptivity} show that the distribution estimators induced by score matching, for both Langevin and forward--backward diffusion, converge at rates that depend only on the intrinsic dimension, without knowing or estimating the manifold, the forward--backward version attaining the minimax $\W_1$ rate up to logarithmic factors when the manifold is smooth enough relative to the density; time discretisation is not analysed. For data with a smooth density bounded above and below on a compact smooth manifold of positive reach, \citet{azangulov2024convergence} prove that score estimation and $\W_1$ sampling achieve rates whose exponents depend only on the intrinsic dimension and the smoothness, the ambient dimension entering through logarithmic factors when $\log d=O(\log n)$ for $n$ samples, by showing that the score is governed by points near the denoised preimage; \citet{potaptchik2024linear} show that the number of steps can then be linear in the intrinsic dimension. On the sampling side, DDPM adapts to an unknown intrinsic dimension $k$ provided its update coefficients are chosen suitably, with error $O(k^2/\sqrt N)$ after $N$ steps, up to logarithmic factors \citep{li2024adapting}, later sharpened to $O(k/N)$ in total variation under the same coefficients and only a finite first moment \citep{li2024od}, and an iteration complexity nearly linear in $k$, which is optimal in relative entropy \citep{huang2024denoising}; total-variation versions are in \citet{liang2025lowdim}, including one for the deterministic DDIM sampler that also requires accurate first and second derivatives of the score. Under only moment assumptions, \citet{chakraborty2026generalization} obtain finite-sample Wasserstein rates governed by a Wasserstein-type intrinsic dimension, connecting diffusion models to the sharp rates of empirical optimal transport. Together these results suggest that diffusion models adapt to intrinsic dimension without being told it, both in the samples they need and, with suitable coefficients, in the steps they take.

\paragraph{Wasserstein guarantees.}

Wasserstein bounds obtained by propagating errors along the reverse process require that process to be stable, and this stability comes from regularity of the score or of the data law. \citet{kwon2022score} show that score-based training controls a Wasserstein distance under smoothness. \citet{gao2023wasserstein} give $\W_2$ guarantees for a general class of forward processes under strong log-concavity and smoothness, and \citet{gentilonisilveri2025beyond} relax this to weak log-concavity of the data law together with a one-sided Lipschitz condition on the gradient of its log-density, without regularity assumptions on the score or its estimator. \citet{mimikos2024score} derive robustness bounds through the regularising properties of the reverse dynamics, and their analysis connects to a Wasserstein-proximal view of memorisation \citep{zhang2024wasserstein}. Direct Wasserstein analyses of this kind thus trade the weak assumptions of the total-variation theory for regularity of the data law, such as log-concavity or smoothness. For bounded, manifold-supported or merely moment-bounded data, the Wasserstein rates of the preceding paragraph rely instead on early stopping and on the low-dimensional structure of the data.

\subsection{Statistical optimality}

Combining approximation and estimation errors gives end-to-end rates. \citet{oko2023diffusion} show that diffusion models with neural-network score estimators are nearly minimax optimal in total variation for Besov-smooth densities bounded below on a cube (and smooth near its boundary), and also in $\W_1$ when separate score networks are trained on a geometric partition of the time interval. \citet{zhang2024minimax} remove the density lower bound with a kernel-based score estimator, \citet{cai2025minimax} obtain minimax rates for ODE sampling, and \citet{wibisono2024optimal} pin down the minimax rate of score estimation itself. These results locate the statistical difficulty in estimating the score. Given an estimator with the right rate, sampling adds little, and the generated law attains near-minimax rates of density estimation. Learnability questions for structured families have also been answered in part. Gradient descent on the DDPM objective learns mixtures of two spherical Gaussians from a random initialisation, and of $K$ well-separated spherical Gaussians from a warm start \citep{shah2023learning}. A diffusion-based algorithm learns mixtures of $K$ identity-covariance Gaussians, under a minimum-weight assumption, with quasi-polynomial time and sample complexity \citep{gatmiry2024learning}. \Cref{tab:theory} collects representative results.

\begin{table}[t]
\centering\small
\begin{tabularx}{\textwidth}{@{}P{3.6cm}P{2.3cm}L@{}}
\toprule
Result & Sampler & Assumptions and conclusion \\
\midrule
\citet{chen2023sampling} & SDE & Lipschitz scores, finite second moment, $L^2$ score error; TV error polynomial in $d$, $L$ \\
\citet{chen2023improved} & SDE & minimal smoothness with early stopping; user-friendly bounds \\
\citet{benton2024nearly} & SDE & finite second moment, early stopping; steps linear in $d$ up to logs (KL) \\
\citet{jain2025sharp} & ODE with noising & finite second moment, early stopping; steps of order $d/\varepsilon$ up to logs for KL accuracy $\varepsilon^2$ \\
\citet{chen2026high} & first-order rejection sampling & finite second moment, $L^2$ score error; steps $\widetilde O(d_\star\,\mathrm{polylog}(1/\varepsilon))$ (bounded Lipschitz) \\
\citet{li2024od} & SDE & minimal assumptions, early stopping; TV $\lesssim d/N$ plus score error \\
\citet{chen2023probability} & ODE + corrector & polynomial guarantees for ODE sampling \\
\citet{li2024sharp} & ODE & score and Jacobian accuracy; nearly linear in $d$ \\
\citet{azangulov2024convergence} & SDE & manifold data; rates in intrinsic dimension, ambient dimension only in logarithms \\
\citet{potaptchik2024linear} & SDE & manifold data; steps linear in intrinsic dimension \\
\citet{oko2023diffusion} & SDE & Besov densities; near-minimax TV \\
\citet{cai2025minimax} & ODE & H\"older densities, sub-Gaussian; near-minimax TV \\
\citet{fukumizu2025flow} & flow matching & Besov densities bounded above and below, noise scale $\sim\sqrt{1-t}$ near the data, time-partitioned networks; near-minimax $\W_p$, $1\le p\le2$ \\
\citet{su2025flow} & flow matching & KL bound; near-minimax TV \\
\bottomrule
\end{tabularx}
\caption{A selection of convergence and statistical results, read as constraints on the three terms of \eqref{eq:error-decomp}. Exact rates and constants are in the cited papers.}
\label{tab:theory}
\end{table}

\section{Generalisation, memorisation and the geometry of learned scores}\label{sec:generalisation}

A perfectly trained diffusion model would reproduce its training set. Trained models usually do not, and explaining why is one of the central questions about these models. The explanations proposed so far are of five kinds. They are geometric, through the low-dimensional structure of the data; architectural, through the inductive biases of the denoising network; statistical, through the noise in the regression target; dynamical, through the time scales of training; and representational, through the geometry the model induces on its own features.

\paragraph{Memorisation is the default.}

\Cref{prop:empirical} shows that the exact minimiser of the training objective memorises. Empirically, large models do replicate training images in some regimes \citep{somepalli2023diffusion,carlini2023extracting}, and memorisation appears when the training set is small relative to model capacity \citep{gu2023memorization}. Generalisation is therefore a phenomenon of implicit bias, and several complementary explanations are now available.

\paragraph{Manifolds and intrinsic dimension.}

Near $t=0$ the score points approximately towards the data manifold, along its normal bundle. Its Jacobian has eigenvalues close to $-1/\sigma_t^2$ in the normal directions and of lower order along the manifold (bounded for a flat manifold, of order $\kappa/\sigma_t$ at typical noisy points near a manifold of curvature $\kappa$), so the score itself is of order $1/\sigma_t$ (times $\sqrt{d-k}$ for a $k$-dimensional manifold) at typical noisy points. \citet{pidstrigach2022score} proves that score-based models with suitably accurate scores generate samples on the data manifold and analyses when they reproduce training data. The same geometry can be read off from a trained model, since the singular values of score vectors sampled in a neighbourhood \citep{stanczuk2024your}, or the rate of change of the log-density under Gaussian smoothing computed from the Fokker--Planck equation \citep{kamkari2024geometric}, estimate the local intrinsic dimension; \citet{loaizaganem2024deep} survey generative modelling under the manifold hypothesis. Recent theory shows that geometric structure is learned before densities are. With coarsely learned scores, diffusion models attain a near-parametric rate towards a distribution whose density is comparable to that of $\mu_{\rm data}$ near the manifold, exploiting the regularity of the support rather than of the density, whereas estimating $\mu_{\rm data}$ itself on a $k$-dimensional manifold has minimax rate of order $n^{-1/k}$ in $\W_1$ (for $k\ge3$, with no smoothness assumed of the density) in the number $n$ of samples \citep{shen2026manifold}. Along the reverse dynamics, gaps in the spectrum of the Jacobian of the score reveal a trivial phase, then a manifold-coverage phase that fits the distribution within the manifold, and finally a consolidation phase in which the score becomes orthogonal to the manifold and projects samples onto it; the gaps are computed by random-matrix analysis for linear manifolds with Gaussian internal density, and qualitative signatures of this picture appear in networks trained on images, most clearly a sharp manifold gap on MNIST \citep{ventura2025manifolds}. Memorisation, too, can occur geometrically, by losing manifold dimensions, before it collapses onto points \citep{achilli2024losing}. Smoothing the empirical score, a proxy for the implicit regularisation of training, is itself geometry-adaptive. Because it smooths in the log-density domain, it spreads mass along the data manifold while keeping samples concentrated near it, whereas kernel density estimation, which smooths the density itself, pushes samples off the manifold as the bandwidth grows, and the choice of smoothing selects the manifold along which the model generalises \citep{farghly2025diffusion}. Geometry thus offers one answer to why trained models generalise, since what is learned first, and learned well, is the support of the data rather than its density.

\paragraph{Inductive biases of denoisers.}

Two denoisers trained on disjoint halves of a large data set converge to nearly the same function, and their behaviour is explained by a shrinkage in geometry-adaptive harmonic bases \citep{kadkhodaie2024generalization}. The agreement extends to samples and across designs. Given the same initial noise and a deterministic sampler, diffusion models that differ in framework, architecture and training procedure produce nearly the same images, both in a memorising regime, where they all reproduce the same training points, and in a generalising regime, where they produce new ones; generative adversarial networks and variational autoencoders generally lack this reproducibility \citep{zhang2024emergence}. A linear theory already accounts for much of this agreement. For the optimal linear denoiser fitted to $n$ samples in dimension $d$, deterministic equivalents from random matrix theory show that the sample covariance acts on average as the population covariance at a larger noise variance $\tilde\sigma^2>\sigma^2$, the solution of a Marchenko--Pastur self-consistent equation. The effect is largest at low noise and for $d/n$ large; it over-shrinks low-variance directions and pulls samples towards the mean. The fluctuations across splits factor into anisotropy across eigenmodes, inhomogeneity across inputs and an overall scaling with $n$. Trained U-Nets and transformers follow these predictions qualitatively once they leave the memorising regime \citep{wang2026random}.

Architecture constrains what can be learned. For convolutional networks, locality and translation equivariance prevent learning the empirical score; the minimum-error score subject to these constraints is computable in closed form and, once a single time-dependent locality scale is calibrated, predicts the outputs of trained convolutional diffusion models with high accuracy, producing novel images as mosaics of local training patches \citep{kamb2025analytic}, and patch-based denoisers give a related mechanistic account \citep{niedoba2025towards}. The locality need not be imposed by the architecture, however. The optimal linear denoiser, the Wiener filter of the data covariance, already has the spatial sensitivities of trained networks, local on natural images and nonlocal on aligned faces, and imperceptible edits to the pixel statistics of the training set change the sensitivities a network learns \citep{lukoianov2025locality}. The Wiener filter also predicts the outputs of trained networks better than recent analytical models built on the optimal denoiser, and adding this data-derived locality to the constrained closed-form score above gives the most accurate analytical model so far.

The noise level matters as well. At moderate and high noise, trained scores are close to the score of a Gaussian with the data covariance \citep{wang2024unreasonable}. In the generalising regime, denoisers become increasingly linear and close to the optimal denoiser for a Gaussian with the empirical mean and covariance of the training data, and that Gaussian denoiser reproduces the coarse structure of generated images, although a substantial nonlinear gap remains at intermediate noise levels \citep{li2024understanding}. For mixtures of low-rank Gaussians, with a matching parameterisation of the denoiser, training reduces to subspace clustering \citep{wang2024diffusion}. On this view generalisation reflects the function class and the data statistics it absorbs. An architecture that cannot represent the empirical score represents a smoother, locally consistent field instead, and that field produces new samples.

\paragraph{Noise in the target.}

The denoising target is a noisy estimate of the score, equal to it only on average, and one proposal is that this noise is itself a source of generalisation. \citet{vastola2025generalization} uses a path-integral method to compute the distributions typically learned by expressive linear score models and by infinitely wide networks in the neural-tangent regime, and finds that they resemble the training set with the gaps between training points filled in, an inductive bias due to the covariance of the target, which is largest between training points. \citet{bertrand2025closed} test the idea for flow matching in high dimension and rule it out as a key factor. On image data, the closed-form velocity of the training set aligns with a single conditional target from early in generation (before $t=0.2$ on CIFAR-10, with noise at $t=0$), a transition they relate to the collapse time of \cref{sec:statphys}. Regressing on the closed form instead of the noisy target does not hurt generalisation and can improve it, and generalisation appears where the network fails to fit the closed-form velocity, early in generation. On low-dimensional toy data the noisy regime lasts longer, as in Vastola's examples, so the dimension of the data governs how long the target stays noisy, and with it how much that noise can matter.

\paragraph{Training dynamics.}

\citet{bonnaire2025why} identify two time scales in training, an early time $\tau_{\rm gen}$ after which samples are of high quality and a later time $\tau_{\rm mem}$ after which memorisation emerges. The first is roughly independent of the training-set size $n$, while the second grows linearly in $n$, which opens a widening window of generalisation. They support this with experiments on U-Nets and an analysis of a random-features model in high dimension. This implicit dynamical regularisation complements the architectural and geometric explanations above and connects to the statistical physics of \cref{sec:statphys}.

\paragraph{The geometry of latent and feature spaces.}

Diffusion models also induce geometry on their internal representations. The bottleneck features of the denoising network form a semantic latent space in which edits are approximately linear \citep{kwon2023diffusion}. Pulling back the Euclidean metric of this bottleneck space to the noisy input $x_t$ through the encoder yields a Riemannian structure whose principal directions, the right singular vectors of the encoder Jacobian, give interpretable edits \citep{park2023understanding}. At the level of the whole model, Wasserstein-proximal operators give a variational description of score-based models that explains memorisation and suggests regularisation \citep{zhang2024wasserstein}.

Apart from the contested role of target noise, these explanations complement one another. Geometry describes what can be learned, architecture what can be represented, training dynamics when memorisation sets in and feature geometry how the result is organised. None of them alone yet predicts when a given model will memorise, which is the last of the open problems below.

\begin{openprob}[Open problems for Part III]
\begin{problem}
Design covering fields, and choices of the free parameter of divergence--kernel estimators, that minimise score-estimator variance, and quantify the gain in sample complexity of any variance-optimal target, including the best time-dependent combination of the targets in \cref{tab:ibp}.
\end{problem}
\begin{problem}
Make nonlinear and state-dependent forward processes practical at scale. Malliavin estimators require simulating variation processes and inverting Malliavin covariances, and divergence--kernel estimators require simulating forward covector processes; low-rank or learned surrogates are needed in high dimension.
\end{problem}
\begin{problem}\label{prob:intrinsic-ode}
Prove intrinsic-dimension guarantees for deterministic samplers under conditions weaker than Jacobian accuracy of the score estimate (for instance in a Wasserstein distance, or for estimators that are regular by construction), or show that no such weakening is possible; in total variation, $\ell_2$ score accuracy alone does not suffice \citep{li2024sharp}.
\end{problem}
\begin{problem}
Unify the geometric, architectural and dynamical explanations of generalisation in one model that predicts $\tau_{\rm mem}$ from data geometry and architecture.
\end{problem}
\end{openprob}

\part{Conditioning, bridges and control}

\section{Conditioning as a Doob \texorpdfstring{$h$}{h}-transform}\label{sec:conditioning}

The preceding parts were concerned with learning a generative process. This part is concerned with changing one, by conditioning it on an observation, bridging it to a prescribed endpoint, or steering it towards a target known only through an energy. All three are changes of measure on path space, and the tools of \cref{sec:pathspace} apply directly. We begin with conditioning, where the change of measure is a Doob $h$-transform and every practical method is a way of approximating $h$.

\subsection{Conditional scores and guidance}

\paragraph{Conditional scores and $h$-functions.}

Let $(X_0,Y)$ have joint law $\mu_{\rm data}(\dd x)\,p(y\mid x)$ and run the forward diffusion on $X$. Bayes' rule at each noise level gives
\begin{equation}\label{eq:cond-score}
\begin{aligned}
\nabla_x\log p_t(x\mid y)&=\nabla_x\log p_t(x)+\nabla_x\log p_t(y\mid x),\\
p_t(y\mid x)&=\E\big[p(y\mid X_0)\,\big|\,X_t=x\big].
\end{aligned}
\end{equation}
In the language of \cref{sec:doob}, $h(t,x)=p_t(y\mid x)$ is space--time harmonic for the reverse process, and conditioning on $y$ is the $h$-transform that adds $g^2\nabla\log h$ to the reverse drift, started from $p_T(\cdot\mid y)$, which is close to $p_T$ when the terminal noise level is high. Every conditioning method is a way of approximating $h$ or $\nabla\log h$.

\paragraph{Classifier and classifier-free guidance.}

Classifier guidance trains a noise-conditional classifier $p_\phi(y\mid x_t)$ and adds $w\,g^2\nabla\log p_\phi(y\mid x_t)$ with a guidance scale $w$ \citep{dhariwal2021diffusion}. Classifier-free guidance trains one network on conditional and unconditional inputs, whose outputs $s_\theta(x,y)$ and $s_\theta(x)$ approximate $\nabla\log p_t(x\mid y)$ and $\nabla\log p_t(x)$, and uses $s_{\rm cfg}=(1+w)s_\theta(x,y)-w\,s_\theta(x)$, the unconditional score plus $(1+w)\nabla\log p_t(y\mid x)$ \citep{ho2022classifier}. The two scales are offset by one, since the conditional score corresponds to $w=1$ in classifier guidance and to $w=0$ in classifier-free guidance. From here on $w$ is the classifier-free scale, and $p=p_0$ denotes the data density. For $w>0$ this combination is the score of $p_t(x)p_t(y\mid x)^{1+w}$ at each fixed $t$, but these tilted marginals are not the forward-noised marginals of $p(x)p(y\mid x)^{1+w}$, or of any single law, so the reverse SDE driven by $s_{\rm cfg}$ does not have them as marginals and does not target $p(x)p(y\mid x)^{1+w}$. \citet{bradley2024classifier} show what it does instead. In the SDE limit, classifier-free guidance with the DDPM sampler is equivalent to a predictor--corrector scheme, a DDIM step for $p_t(x\mid y)$ followed by a Langevin step for the law proportional to $p_t(x)^{-2w}p_t(x\mid y)^{1+2w}$, whose exponent differs from the guidance weight. Feynman--Kac correctors aim at the tilted law itself. They weight the guided SDE by the Feynman--Kac terms that a change of score leaves out and resample with sequential Monte Carlo, so that the particles approximate $p_t(x)p_t(y\mid x)^{1+w}$ at every time \citep{skreta2025feynman}.

Analyses in simple settings describe the consequences. Guidance pushes samples towards extremes of the class-conditional support rather than sampling the tilted law \citep{chidambaram2024what}; for Gaussian mixtures it increases classification confidence while reducing diversity \citep{wu2024theoretical}; and in the high-dimensional limit of Gaussian mixtures these distortions vanish and guidance reproduces the target distribution, a property shared by a family of guidance forms that includes nonlinear ones \citep{lehmanpavasovic2025classifier}. This holds when the number of classes grows sub-exponentially with the dimension, a regime in which a dynamical phase transition makes the distortions vanish. With exponentially many classes the distortions persist, and with a continuum of classes they remain in every dimension; standard schedules with positive guidance cannot prevent the variance from shrinking, whereas a schedule with a window of negative guidance improves both class separation and diversity in latent diffusion models \citep{ventura2026emergence}. Practice has moved accordingly. Guidance helps most at intermediate noise levels \citep{kynkaanniemi2024applying}, guiding with a weaker version of the model itself avoids the loss of diversity \citep{karras2024guiding}, and for masked discrete diffusion an exact analysis with one or two tokens favours weak guidance early and strong guidance late, and shows that existing discrete guidance must be renormalised so that it does not speed up unmasking \citep{rojas2025theory}. Guidance is therefore best understood as a controlled departure from the $h$-transform rather than as an exact conditional sampler, and \cref{prob:cfg} asks for a theory of the law it does sample.

\subsection{Posterior sampling}

When the observation is a measurement of the data, the $h$-function is a likelihood averaged over the posterior of the clean signal, and it can be approximated directly, or corrected by sequential Monte Carlo so that the conditioning becomes exact as the number of particles grows.

\paragraph{Inverse problems.}

For $y=\mathcal A(x_0)+\text{noise}$, the $h$-function $p_t(y\mid x)$ is an intractable integral against the posterior of $X_0$ given $X_t$. Diffusion posterior sampling approximates it by evaluating the likelihood at the Tweedie denoiser, $p_t(y\mid x)\approx p(y\mid\hat x_0(x))$, and differentiates through the denoiser \citep{chung2023diffusion}. Other methods replace this point estimate by more structured approximations. Pseudoinverse-guided models use a Gaussian approximation of the posterior and the pseudoinverse of $\mathcal A$ \citep{song2023pseudoinverse}, and denoising diffusion restoration models work in the singular basis of a linear $\mathcal A$ \citep{kawar2022denoising}. A variational formulation regularises the measurement fit with a weighted score-matching term in which denoisers at different noise levels impose structural constraints \citep{mardani2024variational}. Score priors were applied early to medical imaging \citep{song2022solving}. The idea that a denoiser defines a prior usable for linear inverse problems was also developed independently of diffusion models, from Miyasawa's identity \citep{miyasawa1961empirical} and a network trained for blind least-squares denoising \citep{kadkhodaie2021solving}. These methods differ mainly in how they approximate the intractable $h$-function, by a point estimate, a Gaussian or a variational objective, and each trades exactness for speed.

One-step flow maps (\cref{sec:flowmaps}) leave no trajectory to guide, so variational flow maps condition the initial noise instead. An amortised Gaussian adapter $q_\phi(z\mid y)$ is trained jointly with the one-step map $f_\theta$ from noise to data, by a variational bound that matches $q_\phi(z\mid y)\,p(y\mid x)\,p(x)$ to a model in which $x$ and $y$ are Gaussian around $f_\theta(z)$ and $\mathcal A(f_\theta(z))$, with widths $\tau$ and $\sigma$; the data-fit term is replaced by a mean-flow loss on noise--data pairs coupled through the adapter, which keeps $f_\theta$ a flow map. Sampling $z\sim q_\phi(\cdot\mid y)$ and returning $f_\theta(z)$ targets the posterior in the limit $\tau\to0$, which is hard to train, so in practice $\tau$ exceeds $\sigma$. Joint training compensates for the limited expressiveness of a Gaussian adapter. In a linear-Gaussian model with a diagonal adapter, a frozen generator generically leaves the latent posterior correlated across the adapter's axes, which the adapter cannot represent, whereas joint training rotates the latent coordinates so that the principal axes of the posterior align with those of the adapter. On ImageNet at $256\times256$, in the latent space of an autoencoder, a single network evaluation gives inpainting and deblurring samples with better FID and calibration than guidance methods that use 500 evaluations, though a single draw has lower pixel fidelity; replacing the observation term by a differentiable reward fine-tunes a flow map to sample a reward-tilted law in one step \citep{mammadov2026variational}.

Under standard cryptographic assumptions, no posterior sampler can be both fast and even approximately correct in the worst case. If one-way functions exist, there are distributions whose smoothed scores are approximated by small networks, so that unconditional diffusion sampling takes polynomial time, and linear measurements with small Gaussian noise for which approximate posterior sampling takes superpolynomial time \citep{gupta2024diffusion}. If some one-way function takes exponential time to invert, approximate posterior sampling takes time exponential in the number of measurements. Rejection sampling, whose cost is of that order, is then essentially optimal. \citet{daras2024survey} survey the area.

\paragraph{Asymptotically exact conditioning by sequential Monte Carlo.}

Approximations of $h$ bias the samples, and sequential Monte Carlo restores exactness. The twisted diffusion sampler uses approximate $h$-functions as twisting potentials in a particle filter along the reverse process, and is asymptotically exact as the number of particles grows \citep{wu2023practical}. Related schemes were developed for protein motif scaffolding \citep{trippe2023diffusion} and for Bayesian linear inverse problems \citep{cardoso2024monte}, and score-based data assimilation conditions entire trajectories on observations \citep{rozet2023score}. Feynman--Kac steering places these particle methods, the twisted sampler included, in the framework of Feynman--Kac interacting particle systems for the reward-tilted law $p(x_0)\,e^{\lambda r(x_0)}$, where $p$ is the density of the pretrained model, $r$ a reward and $\lambda>0$ the strength of the tilt (not the log-SNR). The potentials are built from intermediate rewards that need not be differentiable, for discrete as well as continuous models, and a 0.8B-parameter text-to-image model steered in this way exceeds a fine-tuned 2.6B-parameter model in prompt fidelity \citep{singhal2025general}. Inference-time search can also be organised as in classical artificial intelligence, with breadth-first and depth-first tree search for global exploration and annealed Langevin moves for local refinement \citep{zhang2025inference}.

\subsection{Learning and composing \texorpdfstring{$h$}{h}-transforms}

Instead of approximating $h$ at sampling time, one can learn it once and reuse it, or combine several models into one.

\paragraph{Learning the $h$-transform.}

DEFT fine-tunes a small network to represent the generalised $h$-transform of a pretrained model \citep{denker2024deft}. Doob's matching estimates $h$ and its gradient jointly by gradient-regularised regression and proves non-asymptotic rates for the guidance and for the generated law in $\W_2$ \citep{chang2026inference}. For bridges, \citet{heng2025simulating} learn the time reversal of the unconditioned diffusion by score matching, which yields a backward representation of the bridge, and in a second step approximate Doob's $h$-transform for its forward representation. These are learned counterparts of the analytic guided proposals of \citet{delyon2006simulation} and \citet{schauer2017guided}. When the process being conditioned is a nonlinear SDE, the $h$-function lacks Gaussian structure. \citet{pidstrigach2025conditioning} use Malliavin integration by parts to obtain a generalised Tweedie formula for the conditional score that remains well behaved near the singularities of bridge conditioning, and use it to add controls to trained diffusion and flow-matching models. Learning $h$ thus amortises the cost of conditioning, which is paid once in training rather than at every sampling step.

\paragraph{Composition.}

Products of experts $p_1(x)p_2(x)$ are natural for compositional generation, but the sum of noisy scores is not the noisy score of the product. Composable diffusion models sum scores anyway \citep{liu2022compositional}, and \citet{du2023reduce} trace the failure to the sampler rather than the models. Annealed MCMC run with the summed scores at each noise level targets a sequence of densities that ends at the product, and an energy parameterisation of the denoiser further allows Metropolis corrections and mixtures. \citet{skreta2025superposition} estimate the log-density of each model along a trajectory with an estimator based on It\^o's formula, without divergence computations, and use it to superpose diffusion models as mixtures or as logical conjunctions. Feynman--Kac correctors treat products of models and temperature-annealed models in the same way as guidance, and can sample a Boltzmann density at temperatures other than the one at which a sampler was trained, without retraining \citep{skreta2025feynman}.

\section{Schr\"odinger bridges and entropic optimal transport}\label{sec:sb}

Conditioning fixes an observation; a bridge fixes an entire endpoint law. The Schr\"odinger problem of \cref{sec:schrodinger-problem} asks for the path measure that connects two prescribed laws while staying as close as possible, in relative entropy, to a reference diffusion. It sits between diffusion models and optimal transport, and the algorithms that solve it inherit ideas from both.

\paragraph{From the Schr\"odinger problem to generative modelling.}

Diffusion models solve a degenerate version of the problem \eqref{eq:sb}, in which the forward process is fixed and its terminal law is only approximately Gaussian. Schr\"odinger bridges instead connect two prescribed laws exactly, in finite time, by the entropically optimal process, which makes them natural for translation between data sets and for shortening the generative path. The laws need only admit some path measure with these marginals and finite relative entropy with respect to the reference. \citet{tang2026foundations} is a thorough recent guide.

\paragraph{Iterative proportional fitting on path space.}

The classical algorithm is iterative proportional fitting (IPF), which alternately enforces the two marginal constraints by solving half-bridge problems. On the static level it is Sinkhorn's algorithm, and in the continuous setting it converges, in total variation and in relative entropy, to the solution of the static Schr\"odinger problem under boundedness conditions on the reference density \citep{ruschendorf1995convergence}. Diffusion Schr\"odinger bridges implement IPF on path space with neural drifts. Each half-bridge is the time reversal of the previous process with corrected initial law, learned by a regression analogous to score matching \citep{debortoli2021diffusion}. Other neural Schr\"odinger-bridge solvers that appeared at about the same time are IPF with Gaussian-process drifts \citep{vargas2021solving}, a likelihood formulation based on forward--backward SDEs \citep{chen2022likelihood}, and a two-stage construction without IPF, which drives a point mass to a Gaussian-smoothed data law and then to the data with drifts built from score and density-ratio estimators \citep{wang2021deep}. Conditional versions handle posterior sampling \citep{shi2022conditional}.

\paragraph{Markovian and reciprocal projections.}

\begin{figure}[t]
\centering
\begin{tikzpicture}[>=Stealth,lab/.style={font=\footnotesize}]
\begin{scope}[xshift=0cm]
\draw[figblue,line width=1pt] (0.200,0.900) -- (5.800,2.900);
\draw[figorange,line width=1pt] (0.200,2.900) -- (5.800,1.300);
\node[lab,text=figblue,anchor=south east] at (5.800,3.000) {terminal law $\nu$};
\node[lab,text=figorange,anchor=north east] at (5.800,1.200) {initial law $\mu$};
\draw[->,black!80,shorten >=0.5pt] (0.750,1.900) -- (1.005,1.187);
\draw[->,black!80,shorten >=0.5pt] (1.005,1.187) -- (1.396,2.558);
\draw[->,black!80,shorten >=0.5pt] (1.396,2.558) -- (1.786,1.466);
\draw[->,black!80,shorten >=0.5pt] (1.786,1.466) -- (2.045,2.373);
\draw[->,black!80,shorten >=0.5pt] (2.045,2.373) -- (2.303,1.651);
\fill[black] (0.750,1.900) circle (1.6pt) node[left,lab] {$\Rref$};
\fill[black] (3.311,2.011) circle (2.1pt);
\node[lab,anchor=west] at (3.431,1.991) {bridge};
\node[font=\small,anchor=south] at (3,3.55) {Iterative proportional fitting};
\node[lab,anchor=north] at (3,0.55) {each step matches one marginal};
\end{scope}
\begin{scope}[xshift=7.4cm]
\draw[figblue,line width=1pt] (0.200,0.900) -- (5.800,2.900);
\draw[figorange,line width=1pt] (0.200,2.900) -- (5.800,1.300);
\node[lab,text=figblue,anchor=south east] at (5.800,3.000) {Markov measures};
\node[lab,text=figorange,anchor=north east] at (5.800,1.200) {reciprocal class of $\Rref$};
\draw[->,black!80,shorten >=0.5pt] (0.648,2.772) -- (1.190,1.254);
\draw[->,black!80,shorten >=0.5pt] (1.190,1.254) -- (1.550,2.514);
\draw[->,black!80,shorten >=0.5pt] (1.550,2.514) -- (1.909,1.510);
\draw[->,black!80,shorten >=0.5pt] (1.909,1.510) -- (2.147,2.344);
\draw[->,black!80,shorten >=0.5pt] (2.147,2.344) -- (2.384,1.680);
\fill[black] (0.648,2.772) circle (1.6pt) node[above,lab] {$\P^0$};
\fill[black] (3.311,2.011) circle (2.1pt);
\node[lab,anchor=west] at (3.431,1.991) {bridge};
\node[font=\small,anchor=south] at (3,3.55) {Iterative Markovian fitting};
\node[lab,anchor=north] at (3,0.55) {both marginals $\mu,\nu$ kept throughout};
\end{scope}
\end{tikzpicture}
\caption{Two ways of reaching the Schr\"odinger bridge by alternating projections, drawn schematically. In the left panel, iterative proportional fitting starts from the reference $\Rref$ and alternates between path measures of the form $\mathfrak f(X_0)\mathfrak g(X_T)\,\Rref$ with the prescribed terminal law and those with the prescribed initial law, so each iterate matches one marginal at a time. In the right panel, iterative Markovian fitting starts from a measure $\P^0$ in the reciprocal class of $\Rref$, whose elements are mixtures of reference bridges. It keeps both marginals at every step and alternates between Markov measures and that class. In both panels the bridge lies where the two sets meet. In the left panel it is the only measure of this product form with both marginals, which is the path measure with both marginals closest to $\Rref$ in relative entropy.}
\label{fig:ipf-imf}
\end{figure}
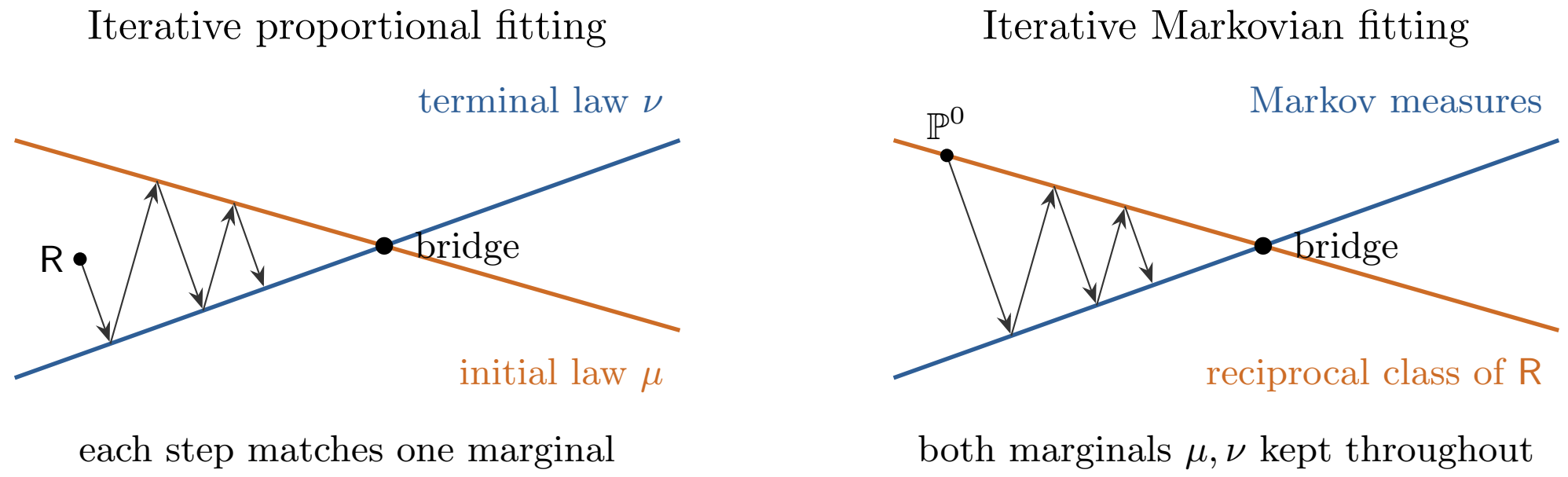

A second family of algorithms alternates two projections. The reciprocal projection of a path measure $\P$ keeps its endpoint coupling $\P_{0T}$ and replaces the bridges by reference bridges, $\int\Rref^{xy}\P_{0T}(\dd x,\dd y)$. The Markovian projection replaces a non-Markov process $\P$ by the Markov diffusion with the same diffusion coefficient whose drift is the conditional expectation, given $X_t$, of the drift of $\P$. It has the same one-time marginals and is the Markov law closest to $\P$ in relative entropy. \Cref{fig:ipf-imf} contrasts the two schemes.

\begin{proposition}[Markovian projection of a bridge mixture]\label{prop:markov-proj}
Let $\Rref$ be Brownian motion with diffusion coefficient $\sigma$ and $\Rref^\pi=\int\Rref^{xy}\,\pi(\dd x,\dd y)$ the mixture of its bridges over a coupling $\pi$ with first marginal $\pi_0$. Under regularity conditions ensuring that the SDE with drift $v$ is well posed on $[0,T)$ and that its Fokker--Planck equation started from $\pi_0$ has a unique solution, the diffusion with drift
\begin{align*}
v(t,z)=\E_{\Rref^\pi}\Big[\frac{X_T-X_t}{T-t}\,\Big|\,X_t=z\Big]
\end{align*}
and diffusion coefficient $\sigma$, started from $\pi_0$, has the same one-time marginals as $\Rref^\pi$.
\end{proposition}
\begin{proof}
Under $\Rref^{xy}$ the drift of the Brownian bridge at $(t,z)$ is $(y-z)/(T-t)$. The marginals of $\Rref^\pi$ therefore solve a Fokker--Planck equation whose drift is the mixture-average of bridge drifts, which is the conditional expectation above; uniqueness for this equation identifies its solution with the marginals of the diffusion. This is \cref{prop:marginalisation} with Markov generators.
\end{proof}
Iterative Markovian fitting alternates the two projections, and its fixed point is the Schr\"odinger bridge. Diffusion Schr\"odinger bridge matching implements it with bridge-matching regressions that avoid the error accumulation of IPF \citep{shi2023diffusion,peluchetti2023diffusion}. Online variants avoid explicit iterations and scale to unpaired image translation \citep{debortoli2024schrodinger}, and a unified view places flow matching, minibatch optimal transport and Schr\"odinger bridge matching in one framework \citep{kim2025unified}. On finite state spaces, discrete-time iterative Markovian fitting converges in relative entropy to the Schr\"odinger bridge for marginals and a Markov reference with full support, which gives categorical bridge matching, applied to unpaired translation of vector-quantised image representations \citep{ksenofontov2025categorical}.

\paragraph{Bridges between paired data.}

When pairs $(x_0,x_1)$ are available, such as a degraded and a clean image or an input and an output, the coupling is given and only the bridges need to be learned. Image-to-image Schr\"odinger bridges exploit a tractable posterior between paired endpoints \citep{liu2023i2sb}, and denoising diffusion bridge models learn the score of a diffusion pinned at a data endpoint through a Doob transform \citep{zhou2024denoising}. Augmented bridge matching adds the initial point to the drift so that the learned process preserves the coupling \citep{debortoli2023augmented}, and aligned bridges handle aligned biological data \citep{somnath2023aligned}.

\paragraph{Simulation-free and closed-form bridges.}

With an entropic optimal-transport coupling between minibatches and Brownian bridges as conditional paths, the score and drift of the bridge mixture can be learned without simulation \citep{tong2024simulation}. In special cases more is possible. Between Gaussians the bridge has a closed form \citep{bunne2023schrodinger}, which yields test cases and initialisations, and parameterising the Schr\"odinger potential by a Gaussian mixture gives a fast light solver \citep{korotin2024light}. Saddle-point formulations give neural entropic optimal transport \citep{gushchin2023entropic}, and benchmark pairs with entropic optimal-transport and Schr\"odinger-bridge solutions known by construction measure how accurately such solvers recover them \citep{gushchin2023building}.

\paragraph{Generalised bridges.}

The reference process and the cost can be enriched. Mean-field and state-cost generalisations correspond to mean-field games and are solved by deep forward--backward schemes \citep{liu2022deep} or by matching algorithms with task-specific costs \citep{liu2024generalized}. When the state cost is a potential, the reference becomes a Feynman--Kac twist of Brownian motion, and iterative Markovian fitting extends exactly \citep{noble2026twisted}. The Markovian projection is then learned with a loss involving the gradient of the potential. When the potential vanishes, this loss reduces to the bridge-matching loss, and in discrete time to the loss of denoising diffusion. In a controlled test the reciprocal projection of twisted Schr\"odinger bridge matching \citep{noble2026twisted} is more accurate than that of the matching algorithm above \citep{liu2024generalized}, and on crowd-navigation and single-cell tasks it attains lower state cost, while the latter matches the terminal marginal more closely. Momentum and multiple marginals give smoother trajectories through several snapshots \citep{chen2023deep}, and reflected reference processes enforce support constraints \citep{haggbom2026reflected}. Trajectory inference, which reconstructs a stochastic process from population snapshots as in single-cell biology, is the statistical counterpart of multi-marginal bridges. Its mathematical theory relates the entropic minimisation to consistent recovery \citep{lavenant2024towards}, and mean-field Langevin dynamics on path space give practical algorithms \citep{chizat2022trajectory}.

\paragraph{Statistical theory.}

Plug-in estimation of Schr\"odinger bridges from the potentials of a static entropic optimal-transport problem between samples, the Sinkhorn bridge, has path-space error $\E[\TV^2]\lesssim n^{-1/2}+C_\tau n^{-1}$, where $\TV$ is taken between the path laws of the estimated and the true bridge on $[0,\tau]$, for a time $\tau$ before the terminal one, the expectation is over the $n$ samples, and $C_\tau$ grows as $\tau$ approaches the terminal time; this holds for fixed regularisation and a target supported in a bounded subset of a smooth $k$-dimensional submanifold. The authors call this the parametric rate. With Euler--Maruyama discretisation, $n\asymp C_\tau\delta^{-2}\vee\delta^{-4}$ samples then give total-variation accuracy $\delta$, and the bound keeps this form when the source is a point mass, which is the F\"ollmer bridge; the constants blow up as the regularisation vanishes or $\tau$ approaches the terminal time, at a rate governed by the intrinsic rather than the ambient dimension of the target \citep{pooladian2024plug}. As the regularisation $\varepsilon$ of the Schr\"odinger problem vanishes, the entropic Brenier map, the conditional mean of the second coordinate of the static bridge given the first, converges to the Brenier map. For smooth, strongly log-concave marginals the first-order correction is $\tfrac{\varepsilon}{2}\nabla\log\mu$, proportional to the score of the source and independent of the target \citep{mulcahy2026diffusion}. The proof approximates the bridge by a noisy McCann interpolation and, under further assumptions, by a mirror Langevin diffusion on the Hessian manifold of the Brenier potential.

Convergence of exact IPF in the static setting was proved under boundedness conditions by \citet{ruschendorf1995convergence}, and quantitative rates, in some regimes exponential, for Sinkhorn's algorithm with continuous marginals have since been obtained \citep{conforti2023quantitative}. For iterative Markovian fitting, quantitative rates are now available as well. With a Brownian or Langevin reference and strongly log-concave marginals, the relative entropy of the static coupling with respect to that of the bridge contracts at every iteration by a factor $L_U/(T\kappa)$, where $L_U$ bounds the regularity of the reference transition densities and $\kappa$ grows with the log-concavity constants of the marginals, so convergence is geometric once the time horizon $T$ is large enough \citep{gentilonisilveri2025exponential}. A smaller but still geometric rate holds for weakly log-concave marginals, and the proof rests on a contraction estimate for the Markovian projection. The bidirectional variant used in practice, which alternates forward and backward fits, in fact interleaves IPF steps with the Markovian and reciprocal projections. This iterative proportional Markovian fitting converges weakly, in discrete or continuous time, for marginals with bounded support, and exponentially for Gaussian marginals in one dimension at any volatility, in continuous time or in discrete time with a single intermediate time, and in higher dimensions for discrete-time fitting at large volatility; exponential convergence for Gaussians in every dimension, at every volatility and for both kinds of fitting, and convergence in general, are conjectured \citep{kholkin2026diffusion}. For the neural versions, in which each projection is only approximated by a learned drift, and for neural IPF, quantitative guarantees in path space, with errors controlled by the training losses, remain open.

\section{Stochastic optimal control for sampling and fine-tuning}\label{sec:control}

The last change of measure in this part is the most general. Steering a diffusion at quadratic cost towards a terminal reward produces, by \eqref{eq:gibbs}, the reference path measure tilted by that reward, so any problem that asks for a tilted law can be posed as a control problem. Two such problems arise in generative modelling. In sampling, the tilt turns a Gaussian into a target known only up to normalisation; in fine-tuning, it turns a pretrained generative model into one that favours a reward. This section treats both, together with the regression methods that have made the control problem tractable.

\subsection{Sampling from unnormalised densities}

Let $\pi\propto\tilde\pi$ be known up to normalisation. Path-space control turns sampling into learning a drift, and the F\"ollmer process of \cref{sec:follmer} is the optimal one.

\begin{proposition}[Control formulation of sampling \citep{tzen2019theoretical,zhang2022path}]\label{prop:control-sampling}
Let $\Rref$ be Brownian motion on $[0,1]$ started at $0$, so that $X_1\sim\gamma_d$ under $\Rref$, let $Z=\int\tilde\pi<\infty$ and $\KL(\pi\,\|\,\gamma_d)<\infty$, and consider $\dd X_t=u_t\,\dd t+\dd W_t$, $X_0=0$, over progressively measurable drifts with $\E\int_0^1\norm{u_t}^2\dd t<\infty$. Then
\begin{align*}
\min_u\ \E\Big[\int_0^1\tfrac12\norm{u_t}^2\dd t+\log\frac{\gamma_d(X_1)}{\tilde\pi(X_1)}\Big]=-\log Z,
\end{align*}
attained by the F\"ollmer drift $u_t^*(x)=\nabla\log Q_{1-t}(\tilde\pi/\gamma_d)(x)$, with $Q$ the heat semigroup of \cref{sec:follmer}, under which $X_1\sim\pi$.
\end{proposition}
This is \eqref{eq:gibbs} with $\Phi=\log(\gamma_d/\tilde\pi)$. The path integral sampler parameterises $u$ by a network and minimises the cost by simulating trajectories \citep{zhang2022path}. F\"ollmer processes were also used for Bayesian posterior sampling, with the drift estimated by Monte Carlo \citep{huang2021schrodinger} or parameterised by a neural network \citep{vargas2023bayesian}. Denoising diffusion samplers use an Ornstein--Uhlenbeck reference and the time-reversal structure of diffusion models \citep{vargas2023denoising}, and \citet{berner2024optimal} recast diffusion-based generative modelling and sampling in a unified control framework.

A central practical difficulty is the variance of the training signal. Relative-entropy gradients have high variance and require on-policy samples, and alternative divergences on path space, in particular the log-variance divergence, give off-policy, lower-variance training \citep{richter2024improved,nusken2021solving}. Losses in the style of generative flow networks with partial trajectories \citep{zhang2024diffusion} and improved off-policy schemes \citep{sendera2024improved} pursue the same goal. The discrete-time and continuous-time versions of these objectives agree in the limit. As the step size $h$ tends to zero, trajectory-level losses of generative-flow-network type converge to divergences between path measures. In the same limit, the local detailed-balance discrepancy vanishes at order $\sqrt h$ exactly when the forward and backward drifts satisfy Nelson's relation (\cref{sec:time-reversal}) for the learned marginal densities, and the vanishing of its mean at order $h$ then gives the Fokker--Planck equation. Training can therefore use coarse, randomly placed time grids of about ten steps while sampling on fine ones, approaching state-of-the-art results at a fraction of the training cost, whereas coarse uniform grids perform worse and are sensitive to the number of sampling steps \citep{berner2026discrete}.

Other methods use the energy more directly. Targets for the score can be estimated by Monte Carlo from the energy, as in iterated denoising energy matching \citep{akhoundsadegh2024iterated}, and particle methods combine diffusion samplers with sequential Monte Carlo \citep{phillips2024particle}. Without any training, the denoiser of stochastic localisation (\cref{sec:localisation}) can be estimated by Markov chain Monte Carlo along the observation process. For a class of non-log-concave targets, the posterior of the sample given the observation becomes strongly log-concave once enough signal has accumulated, while the law of the observation is still strongly log-concave, so the process can be started inside this window \citep{grenioux2024stochastic}. The resulting sampler is competitive with recent methods on multimodal mixtures, Bayesian logistic regression and a high-dimensional field system. \citet{blessing2024beyond} benchmark many of these samplers, among them the path integral, time-reversed and denoising diffusion samplers, against Monte Carlo and variational baselines at scale.

Underdamped versions of these samplers learn forward and backward controls that act, like the noise, only on a velocity variable. They need splitting integrators, since under Euler--Maruyama steps the backward transition kernel degenerates; with hyperparameters such as the time horizon and the prior learned end to end, the underdamped diffusion bridge sampler outperforms its overdamped counterpart and competing methods on standard benchmarks, often with as few as eight steps, a gain the authors describe as empirical \citep{blessing2025underdamped}. \citet{chen2023sampling} give evidence, through a lower bound for their discretisation analysis, that the better step complexity of underdamped Langevin dynamics does not carry over to score-based models built on the critically damped Langevin diffusion, whose learned drift depends on the velocity as well as the position.

\paragraph{Nonequilibrium transport and Jarzynski reweighting.}

A parallel tradition corrects an imperfect transport by importance weights instead of learning a perfect one. Annealed importance sampling moves particles through a sequence of intermediate densities and corrects with importance weights \citep{neal2001annealed}, and Jarzynski's equality (\cref{sec:statphys}) is its continuous-time form. Learned flows inside sequential Monte Carlo reduce the variance of these weights \citep{arbel2021annealed,matthews2022continual}. Flow annealed importance sampling bootstraps a flow with an $\alpha$-divergence objective estimated by annealed importance sampling \citep{midgley2023flow}, and the annealing path itself can be learned \citep{mate2023learning}. Controlled Monte Carlo diffusions learn forward and backward kernels of an annealed process jointly, interpreting the objective through path-space identities of Crooks type \citep{vargas2024transport}, and sequential controlled Langevin diffusions combine this with resampling \citep{chen2025sequential}. NETS adds a learned drift to annealed Langevin dynamics, chosen to make the Jarzynski weights as uniform as possible, and keeps the estimator unbiased through the weights \citep{albergo2025nets}. Samplers that solve the continuity equation of a prescribed annealing path \citep{tian2024liouville,sun2024dynamical} are the deterministic counterparts.

\subsection{Control by regression: SOC matching and adjoint matching}

Minimising a control cost by simulating trajectories and backpropagating through them can be expensive and unstable. The methods of this subsection replace it by regression, in the spirit of the matching objectives of Part~II. Stochastic optimal control matching turns the control problem into a least-squares regression onto a target computed along simulated trajectories, using a reparameterisation that reduces variance \citep{domingoenrich2024stochastic}. Adjoint matching regresses the control onto a ``lean'' adjoint process computed backwards along trajectories \citep{domingoenrich2025adjoint}, and rests on an observation about the noise schedule. Let $r$ be the reward, $\P^u$ and $\P^{\rm base}$ the laws of the controlled and the uncontrolled base process, and $p^{\rm base}$ the density of the base law at time one. Maximising $\E[r(X_1)]-\KL(\P^u\,\|\,\P^{\rm base})$ over controls produces the tilted terminal law $p^{\rm base}(x)e^{r(x)}/Z$, with $Z$ the normalising constant, for every reward only if the base process makes $X_0$ and $X_1$ independent (here, as for flows, $X_0$ is noise and $X_1$ the sample), a property the authors call memoryless. Otherwise the optimal control is in general biased by the initial value function. Adjoint sampling specialises the idea to unnormalised densities, reusing trajectories through a replay buffer and taking many gradient steps per energy evaluation, and scales to molecular conformer generation \citep{havens2025adjoint}.

Trust regions organise the reuse of trajectories. \citet{blessing2025trust} replace a single minimisation of $\KL(\P^u\,\|\,\P^{u^*})$ by a sequence of problems, each constrained to stay within relative entropy $\varepsilon$ of the previous iterate. Each step is again a control problem, with a shifted quadratic cost, and its solution satisfies $\dd\P^{u_{i+1}}/\dd\P^{u_i}\propto(\dd\P^{u^*}/\dd\P^{u_i})^{1/(1+\lambda_i)}$, where $\lambda_i$ is the Lagrange multiplier. Started from the uncontrolled process, the iterates thus form a geometric annealing in path space, $\dd\P^{u_i}/\dd\P^0\propto(\dd\P^{u^*}/\dd\P^0)^{\beta_i}$ with $\beta_i=1-\prod_{j<i}\lambda_j/(1+\lambda_j)$, and the constraint chooses the schedule. To leading order in $\varepsilon$ the iterates are equispaced in the Fisher--Rao distance (\cref{sec:transport}). The multiplier solves a one-dimensional dual problem estimated from a buffer of stored trajectories, and each step is learned off-policy from that buffer with log-variance or SOC-matching losses. Without the target's gradient in the drift, these trust-region samplers avoid the mode collapse that unconstrained losses suffer on Gaussian mixtures beyond about ten dimensions, up to about 150 dimensions; in reward fine-tuning of Stable Diffusion~1.5, with a modified version that fixes the annealing factor $\beta_i$ to one, reusing each buffer lets them match adjoint matching with fewer simulated trajectories.

\subsection{Fine-tuning, alignment and mean-field games}

\paragraph{Fine-tuning and alignment as relative-entropy-regularised control.}

Fine-tuning a pretrained diffusion or flow model to increase a reward $r$ while staying close to it is the control problem above. \citet{uehara2024fine} formulate it as entropy-regularised control and handle the bias from the initial distribution by learning the initial law as well. \citet{tang2024fine} develop the control theory including other regularisers, and \citet{uehara2024understanding} review approaches from reinforcement learning. The practical methods use policy gradients over denoising trajectories \citep{black2024training,fan2023dpok}, direct backpropagation of differentiable rewards through the sampler \citep{clark2024directly}, preference optimisation adapted to diffusion likelihoods \citep{wallace2024diffusion}, or amortised posterior sampling through trajectory-balance objectives \citep{venkatraman2024amortizing}. Policy gradients also extend to flow-matching models, whose samplers are deterministic. Flow-GRPO replaces the probability-flow ODE by an SDE with the same marginals, a member of the family of \cref{prop:interpolant-family} (\cref{fig:eta-family}), so that each Euler--Maruyama step is a Gaussian policy with explicit likelihood ratios and relative entropy to the reference policy. Trained with fewer steps than are used at inference, it raises the GenEval accuracy of a text-to-image model from 63\% to 95\% \citep{liu2025flowgrpo}. Read through \eqref{eq:gibbs}, the relative-entropy-regularised objectives among them, such as DPOK and the control formulations above, target the same Gibbs-tilted path measure, while methods that maximise reward without a relative-entropy term, such as DDPO, correspond to its zero-temperature limit.

\paragraph{Mean-field games.}

When the cost depends on the population law, control becomes a mean-field game. \citet{zhang2023mean} show that normalising flows, score-based models and related methods arise as solutions of mean-field games with different costs and dynamics, and neural solvers for high-dimensional mean-field games and control problems \citep{ruthotto2020machine} provide the computational counterpart.

\begin{openprob}[Open problems for Part IV]
\begin{problem}\label{prob:cfg}
Classifier-free guidance changes the law of the reverse process by a density that is in general not a function of the endpoints, so it is not an $h$-transform. Characterise the law it samples beyond the special settings of \cref{sec:conditioning}, bound its distance from the tilted law $p(x)p(y\mid x)^{1+w}$ in terms of the guidance weight and its schedule, and identify the schedules for which the Feynman--Kac correction of \cref{sec:conditioning} is negligible.
\end{problem}
\begin{problem}
On continuous state spaces, extend the geometric convergence of iterative Markovian fitting beyond weakly log-concave marginals and long time horizons, and that of its bidirectional variant beyond the Gaussian cases of \cref{sec:sb}, and prove quantitative convergence in path space, with errors controlled by the training losses, of its neural implementations, in which each projection is learned, and of neural IPF.
\end{problem}
\begin{problem}
Prove guarantees in high dimension for the off-policy objectives used to train control-based samplers, such as the log-variance, SOC-matching and trust-region schemes of \cref{sec:control}, or find objectives that admit them, and connect the objectives to the Malliavin weights of \cref{sec:ibp}, which are themselves path-space sensitivities.
\end{problem}
\end{openprob}

\part{Beyond Euclidean space and Brownian noise}

\section{Riemannian manifolds}\label{sec:riemannian}

Until now the state space has been Euclidean and the noise Brownian. Many data live elsewhere, on spheres and tori, on groups of rotations, on finite sets and in function spaces. This part relaxes both assumptions. A new state space changes what a velocity field, a Brownian motion, a heat kernel and a score mean, and many constructions of the preceding parts carry over once these objects are redefined. The guiding example is a Riemannian manifold, where the redefinition is most transparent.

\paragraph{Geometric preliminaries.}

On a Riemannian manifold $(\M,g)$, velocity fields are sections of the tangent bundle. The gradient and divergence are defined through the metric, the gradient by $\inner{\nabla_\M f}{X}_g=Xf$ for every vector field $X$ and the divergence by $\dd(\iota_X\dd V_g)=(\divg X)\,\dd V_g$, where $\iota_X$ is the interior product and $\dd V_g$ the Riemannian volume form, so that the divergence theorem holds \citep[pp.~342 and~423, Thm.~16.32]{lee2013introduction}. The divergence does not depend on the orientation, so it is also defined on nonorientable manifolds, where the divergence theorem holds with the Riemannian density in place of $\dd V_g$. The continuity equation and the instantaneous change of variables \eqref{eq:inst-cov} then hold verbatim with the Riemannian divergence. Brownian motion is the diffusion with generator $\frac12\Delta_\M$, constructed for instance by rolling without slipping or by stochastic development \citep{elworthy1982stochastic,hsu2002stochastic}. Its transition density, the heat kernel $p_t(x,y)$, plays the role of the Gaussian kernel. On compact manifolds it has the eigenfunction expansion $p_t(x,y)=\sum_ke^{-\lambda_kt/2}\varphi_k(x)\varphi_k(y)$ in terms of Laplace--Beltrami eigenpairs, $-\Delta_\M\varphi_k=\lambda_k\varphi_k$ with $(\varphi_k)$ orthonormal in $L^2(\dd V_g)$ \citep{grigoryan2009heat}, and for small times Varadhan's formula $\lim_{t\to0}-2t\log p_t(x,y)=d(x,y)^2$ relates it to the geodesic distance \citep{hsu2002stochastic}. Both facts underlie the approximations used in practice.

\paragraph{Flows on manifolds.}

Normalising flows were first built on manifolds by composing a Euclidean flow with charts \citep{gemici2016normalizing}, and then intrinsically on tori and spheres \citep{rezende2020normalizing}. Continuous normalising flows on manifolds integrate a tangent vector field and track the log-density through the Riemannian divergence in \eqref{eq:inst-cov}, either directly or in local exponential-map charts with a volume correction for each chart \citep{mathieu2020riemannian,lou2020neural,falorsi2020neural}. Two approaches avoid simulating the ODE during training. Moser flow parameterises the density as a prior minus the divergence of a learned vector field, following Moser's construction of a diffeomorphism that carries the prior density onto the target density \citep{rozen2021moser}. Matching probability paths gives another simulation-light method \citep{benhamu2022matching}. When the manifold is unknown, flows can learn it jointly with the density \citep{brehmer2020flows}.

\paragraph{Riemannian score-based models.}

Riemannian score-based generative models run Brownian motion, or a Langevin diffusion towards a reference measure, on $\M$ and reverse it \citep{debortoli2022riemannian}. The time-reversal theorem extends to manifolds. If the forward process has drift $b$ and is driven by Brownian motion on $\M$, its path law has finite relative entropy with respect to Brownian motion, and its marginals have smooth positive densities $p_t$, then the reversed process has drift $-b+\nabla_\M\log p_{T-t}$, which is just the score term for Brownian noising. The theorem is stated for complete, connected, orientable and boundaryless $\M$, but its rigorous proof is given for compact $\M$. The scores can be learned by implicit score matching with Riemannian divergences or by denoising score matching with heat-kernel scores, approximated through the eigenfunction expansion or Varadhan's formula. Sampling uses geodesic random walks. Riemannian diffusion models derive a variational bound for the same processes \citep{huang2022riemannian}. The heat kernel is the main computational obstacle, and it has been met in several ways. On symmetric spaces it can be computed accurately enough to scale \citep{lou2023scaling}, mixtures of Riemannian bridges avoid it altogether \citep{jo2024generative}, and on general manifolds it can be learned with a physics-informed network solving the heat equation \citep{ko2026riemannian}. For submanifolds given as level sets, projection schemes yield Riemannian DDPMs \citep{liu2025riemannian}, and Schr\"odinger bridges also extend to compact manifolds \citep{thornton2022riemannian}.

\paragraph{Riemannian flow matching.}

\citet{chen2024flow} extend flow matching to general geometries by prescribing conditional flows through a premetric $\mathsf d(x,x_1)$. A premetric is a nonnegative function on $\M\times\M$ that vanishes exactly on the diagonal. The conditional flow from $x_0$ towards $x_1$ is required to decrease $\mathsf d$ at a prescribed rate, $\mathsf d(x_t,x_1)=\kappa(t)\,\mathsf d(x_0,x_1)$ with $\kappa(0)=1$ and $\kappa(1)=0$, and the conditional velocity of least norm achieving this is
\begin{align*}
u_t(x\mid x_1)=\frac{\dot\kappa(t)}{\kappa(t)}\,\mathsf d(x,x_1)\,\frac{\nabla\mathsf d(x,x_1)}{\norm{\nabla\mathsf d(x,x_1)}_g^2}.
\end{align*}
The gradient is the Riemannian gradient in $x$, and the formula holds wherever $\mathsf d(\cdot,x_1)$ is differentiable with nonzero gradient; for the geodesic distance this excludes only $x_1$ and its cut locus, a set of measure zero. With the geodesic distance on simple manifolds, such as spheres, tori, hyperbolic space and matrix groups, training is simulation-free; on general manifolds such as meshes, spectral distances built from Laplace--Beltrami eigenfunctions replace the geodesic distance. Total-variation rates for Riemannian flow matching on compact and Hadamard manifolds separate discretisation and learning errors \citep{guan2026total}, and few-step flow maps have also been generalised to manifolds \citep{davis2025generalised,zhong2026riemannian}. Viewing flow matching as maximum likelihood for a Gaussian model of the endpoint given the current point \citep{eijkelboom2024variational} suggests other likelihoods on manifolds. Such an objective is compatible if its population minimiser is the marginal velocity. Riemannian flow matching is compatible, whereas a Riemannian Gaussian variant is not; for a distribution on the circle constructed by \citet{ng2026compatible}, its global optimum pushes samples away from the data. Riemannian flow matching in turn places likelihood on velocities that cannot be realised on compact manifolds. An exponential family over endpoints with the conditional velocity as sufficient statistic is both compatible and supported on feasible velocities; it reweights the flow-matching loss so that, on constant-curvature manifolds, errors in speed and direction are treated differently, and it outperforms prior methods on high-dimensional product tori and improves likelihoods on two of four geospatial benchmarks on the sphere \citep{ng2026compatible}.

\paragraph{Constrained domains.}

Convex or polyhedral constraints arise for physical quantities, for images in $[0,1]^d$ and for simplices. Log-barrier metrics and reflected Brownian motion give diffusion models on constrained domains \citep{fishman2023diffusion}, and Metropolis-adjusted random walks provide a cheap alternative with the same continuum limit \citep{fishman2023metropolis}. Reflected diffusion models learn scores satisfying Neumann boundary conditions \citep{lou2023reflected}. Mirror diffusion models, for convex constraint sets, run an unconstrained diffusion in the dual space of a mirror map \citep{liu2023mirror}, and reflected Schr\"odinger bridges can be trained by partially simulation-free bridge matching \citep{haggbom2026reflected}.

\paragraph{Data-dependent metrics and hyperbolic geometry.}

Instead of the ambient metric, one may use a metric that makes paths follow the data. Metric flow matching learns geodesic interpolants under a data-dependent Riemannian metric and produces conditional paths that stay near the data manifold \citep{kapusniak2024metric}, and geometric neural diffusion processes model tensor fields on manifolds with the right equivariance \citep{mathieu2023geometric}. Hyperbolic space, whose exponential volume growth suits hierarchical data, carries wrapped normal distributions \citep{nagano2019wrapped}, normalising flows \citep{bose2020latent} and diffusion models for graphs \citep{fu2024hyperbolic,wen2024hyperbolic}. Hyperbolic space is also the Euclidean section of anti-de Sitter space, which is where \cref{sec:holography} begins.

\section{Lie groups, symmetry and equivariance}\label{sec:lie}

When the data have a symmetry, the state space is better thought of as a space acted on by a group, and two questions arise, namely how a generative model can be made to respect the symmetry, and what generative dynamics live on the group itself. The first is answered by equivariance, the second by diffusion and flow matching on Lie groups. Molecules and proteins are the motivating examples throughout.

\paragraph{Invariance through equivariant dynamics.}

Let a group $G$ act on $\M$ by isometries. If the reference law is $G$-invariant and the velocity field is $G$-equivariant, $v_t(g\cdot x)=\dd g\,v_t(x)$ for $g\in G$, where $\dd g$ is the differential of the isometry $x\mapsto g\cdot x$, then the flow commutes with the action and every $\rho_t$ is $G$-invariant (\citealp{kohler2020equivariant} for linear actions on $\R^n$; \citealp{katsman2021equivariant} for isometries of a Riemannian manifold); the same argument applies to a diffusion with invariant noise and equivariant drift, and differentiating $\rho(g\cdot x)=\rho(x)$ shows that the score of an invariant density is equivariant. Translations of $\R^3$ have no invariant probability measure, so molecular models work on the subspace of zero centre of mass. Equivariance therefore converts a symmetry of the data into an architectural constraint, and it also reduces sample complexity. The sample complexity of estimating the Wasserstein-1 metric and Lipschitz-regularised $\alpha$-divergences between invariant distributions drops in proportion to the group size for finite groups, and is governed by the dimension of a fundamental domain for compact Lie groups; for the MMD the gain also depends on the kernel \citep{chen2023sample}. A Wasserstein path-space analysis gives error bounds for score-based and flow-matching models of symmetric distributions that compare data augmentation quantitatively with equivariant score approximations \citep{chen2024equivariant}.

\paragraph{Diffusion and flows on $SO(3)$, $SE(3)$ and tori.}

On compact Lie groups with a bi-invariant metric the heat kernel is a class function with an explicit character expansion. On $SO(3)$, the isotropic Gaussian $\mathrm{IG}_{SO(3)}$ has a density with respect to Haar measure that depends only on the rotation angle $\omega$, proportional to $\sum_{\ell\ge0}(2\ell+1)e^{-\ell(\ell+1)t/2}\sin((\ell+\frac12)\omega)/\sin(\omega/2)$, where $t$ is the time of the Brownian motion with generator $\frac12\Delta$ that it describes (\cref{tab:kernels}), while the angle itself carries the additional factor $1-\cos\omega$; this makes diffusion models on rotations computable \citep{leach2022denoising,jagvaral2024unified}. Distributions on Lie groups can be reparameterised through the exponential map for variational inference \citep{falorsi2019reparameterizing}. A central application is to proteins, whose backbones are sequences of $N$ rigid frames, one per residue, that is, points of $SE(3)^N$. Diffusion on frames \citep{yim2023se3}, flow matching on frames \citep{yim2023fast}, stochastic and transport-based flow matching on $SE(3)$ \citep{bose2024se3} and sequence-conditioned variants \citep{huguet2024sequence} are now standard, and diffusion over frames drives large-scale protein design \citep{watson2023de}. A recent structure-prediction system, by contrast, generates structures with a diffusion module acting directly on atom coordinates rather than on frames \citep{abramson2024accurate}. Other groups appear elsewhere. Diffusion on tori of torsion angles generates molecular conformers \citep{jing2022torsional}, docking diffuses jointly over translations, rotations and torsions \citep{corso2023diffdock}, and $SE(3)$ diffusion fields learn grasp and motion costs in robotics \citep{urain2023se3}. Flow matching on Lie groups can use exponential curves $s\mapsto g\exp(sX)$, which are the integral curves through $g$ of left-invariant vector fields \citep[Thm.~20.1 and Prop.~20.8(h)]{lee2013introduction}, in place of geodesics; this gives an intrinsic, simulation-free method on Lie groups with surjective exponential map, which for many matrix groups needs only matrix products, inverses, exponentials and logarithms \citep{sherry2025flow}. Kinetic diffusions with trivialised momentum keep the momentum in the fixed Lie algebra rather than in position-dependent tangent spaces, so that Euclidean diffusion machinery transfers to Lie groups \citep{zhu2024trivialized}.

\paragraph{Equivariant flows and diffusions for point clouds.}

For molecules and particle systems the symmetry group is $E(3)$, or more generally $E(n)$, times permutations. The toolkit was built by equivariant normalising flows \citep{satorras2021en}, equivariant diffusion for molecule generation \citep{hoogeboom2022equivariant}, roto-translation-invariant conformer diffusion \citep{xu2022geodiff} and equivariant latent diffusion \citep{xu2023geometric}. Equivariant flow matching aligns samples by the symmetry group before solving the minibatch optimal-transport problem, so that couplings respect the symmetry and paths become shorter \citep{klein2023equivariant}. Hybrid probability transport \citep{song2023equivariant} and mixed continuous--categorical flow matching \citep{dunn2024mixed} handle atom types and coordinates together. Equivariance on manifolds \citep{katsman2021equivariant} and permutation-equivariant flows for point processes and sets \citep{bilos2021equivariant} extend the principle.

\section{Discrete and statistical-manifold state spaces}\label{sec:discrete}

Text, graphs and biological sequences are discrete. On a finite state space there is no gradient and no Brownian motion, yet the main ingredients of the continuous theory have counterparts. A continuous-time Markov chain replaces the forward SDE, ratios of probabilities replace the score, and generator matching replaces flow matching. A second route relaxes discrete data to the probability simplex, where the Fisher--Rao geometry of \cref{sec:transport} becomes the geometry of the state space itself.

\paragraph{Continuous-time Markov chains and their reversal.}

On a finite state space, a forward process is a continuous-time Markov chain with rate matrix $Q_t$, whose off-diagonal entries $Q_t(x,y)\ge0$, $x\ne y$, are jump rates and whose diagonal is $Q_t(x,x)=-\sum_{y\ne x}Q_t(x,y)$. The marginals of the chain solve $\partial_tp_t(y)=\sum_xp_t(x)Q_t(x,y)$.

\begin{proposition}[Time reversal of a continuous-time Markov chain]\label{prop:ctmc-reversal}
Let $t\mapsto Q_t$ be continuous on $[0,T]$. The time reversal of the chain on $[0,T]$ is a Markov chain whose rate of jumping from $y$ to $x\ne y$ at forward time $t\in(0,T]$, for $y$ with $p_t(y)>0$, is
\begin{align*}
\hat Q_t(y,x)=Q_t(x,y)\,\frac{p_t(x)}{p_t(y)},
\end{align*}
with diagonal $\hat Q_t(y,y)=-\sum_{x\ne y}\hat Q_t(y,x)=Q_t(y,y)-\partial_t\log p_t(y)$.
\end{proposition}
The ratios $p_t(x)/p_t(y)$ between neighbouring states are the discrete analogue of the score. The concrete score of \citet{meng2022concrete} is the vector of $p(y)/p(x)-1$ over the neighbours $y$ of $x$, and learning such ratios is the discrete analogue of score matching \citep{campbell2022continuous,sun2023score}. Score entropy is a loss that keeps the ratios positive and recovers them at its minimum \citep{lou2024discrete}. A unified view treats denoising for a broad class of Markov processes, discrete or continuous, under mild regularity conditions \citep{benton2024denoising}.

\paragraph{Discrete diffusion models.}

Structured transition matrices, whether uniform, absorbing or banded, define discrete diffusions in discrete time \citep{hoogeboom2021argmax,austin2021structured} and in continuous time with $\tau$-leaping samplers \citep{campbell2022continuous}. Masked diffusion, in which the absorbing state is a mask token, has emerged as the most effective variant for text, and its structure has been clarified. Its continuous-time ELBO is a weighted cross-entropy \citep{sahoo2024simple,shi2024simplified}. The optimal model is time-independent and models conditional distributions of clean data \citep{ou2025your}, and masked diffusion models are equivalent to time-agnostic masked models, with numerical pitfalls in categorical sampling \citep{zheng2025masked}. Sampling exposes a factorisation error. Each reverse step draws all masked tokens from a product of per-token posteriors, which cannot represent the correlations of the exact joint posterior, so errors accumulate and grow as the number of steps falls. Energy-based diffusion language models multiply the factorised denoising distribution at each step by $e^{-E}$, for an energy $E$ of the whole proposed sequence given the current noisy one, and sample the product by resampling a few proposals from the denoiser. The energy is either learned by noise-contrastive estimation or taken from a pretrained autoregressive model, in which case each step becomes importance sampling from that model with the denoiser as proposal. The correction matters mostly in the early, heavily masked steps; at the scale of GPT-2 it lowers the perplexity bound of the masked diffusion baseline on OpenWebText from 23.8 to 21.5 with a learned energy and to 20.5 with an autoregressive one, and to 17.6, matching the autoregressive model itself, when the autoregressive energy is combined with carrying over the unmasked tokens \citep{xu2025energy}.

Masked diffusion now also works at the scale of large language models. LLaDA trains an 8B-parameter model from scratch \citep{nie2025llada}, and Dream~7B, initialised from the weights of an autoregressive model of the same size, reports performance competitive with that model on standard benchmarks and clear gains on planning tasks such as Sudoku and Countdown \citep{ye2025dream}. Generation need not involve the whole sequence at once. Block diffusion is autoregressive over blocks of tokens and runs masked diffusion within each block, which allows outputs of any length and reuse of cached computation; it traces part of the remaining perplexity gap to autoregression to the variance of the diffusion gradient, apparent even for blocks of a single token, and narrows it with noise schedules that avoid extreme masking rates \citep{arriola2025block}. The order in which tokens are revealed is a further degree of freedom, and its effect depends on the task. Masked diffusion must learn every infilling problem, and for data generated from latent variables some of these are computationally hard, under a conjecture from statistical physics, while the autoregressive problems taken in the right order are easy \citep{kim2025train}. At inference, one sidesteps the hard problems by revealing first the tokens for which the probabilities of the two most likely values differ most. On Sudoku this raises the accuracy of a pretrained masked model from under 7\% to about 90\%, above that of autoregressive models with seven times as many parameters trained with the solving order. On open-ended reasoning the same freedom can hurt. Confidence-based orders defer high-entropy tokens at logical forks, such as ``therefore'', until their uncertainty has collapsed, so that arbitrary-order decoding finds fewer distinct solutions than left-to-right decoding as the number of samples grows; reinforcement learning with an exact left-to-right policy, which keeps parallel decoding at inference, reaches 89.1\% on GSM8K \citep{ni2026flexibility}. Autoregression and diffusion can also share one network. Trained jointly with a next-token loss and a block-wise masked-diffusion loss, with attention bidirectional within blocks and causal across them, a single model can decode left to right, by diffusion block by block, or by self-speculation, in which diffusion drafts several tokens in parallel and the autoregressive pass verifies them against its own next-token predictions. The two losses reinforce each other, and at 8B parameters the model decodes about six times as many tokens per forward pass as an autoregressive model of comparable accuracy \citep{fu2026nemotron}.

Not all of these models use masks. DiffusionGemma, obtained by fine-tuning the mixture-of-experts Gemma~4 model, uses uniform noise and denoises blocks of 256 tokens in parallel, one block after another; it generates about 20 tokens per forward pass while trading some benchmark performance relative to its autoregressive initialisation for much faster decoding \citep{diffusiongemma2026}. Uniform-state diffusion is closely tied to Gaussian diffusion. The argmax of a Gaussian diffusion of one-hot vectors has, at each time, exactly the marginals of a uniform-state discrete diffusion, with an explicit map between the two noise schedules, although the discretised Gaussian path need not itself be a discrete Markov chain \citep{sahoo2025diffusion}. Techniques from the Gaussian side then transfer to this uniform-state model, which its authors call Duo. They include a curriculum that doubles training speed by lowering the variance of the gradient, and a discrete consistency distillation that cuts sampling from 1024 to 8 steps with little loss of quality.

A continuous alternative questions whether discrete noise is needed at all. Flows over one-hot token embeddings, whose denoiser is a tokenwise posterior on the simplex trained by cross-entropy, match state-of-the-art discrete diffusion baselines on LM1B and OpenWebText. Because the factorised posterior determines the velocity exactly, whereas ancestral sampling of a discrete chain incurs factorisation errors, such a flow can be distilled into a flow map (\cref{sec:flowmaps}) whose one-step samples exceed the eight-step quality of recent few-step discrete diffusion language models \citep{lee2026flowmap}. Discrete flow maps take the complementary route of reparameterising any-step flow maps through the mean denoiser, which lives on the probability simplex, so that the Euclidean regression losses of flow-map training become exact cross-entropy and relative-entropy losses; text is then generated in one or a few steps with little loss of quality, and guidance carries over \citep{potaptchik2026discrete}. At larger scale, Sigma denoises Gaussian-corrupted 16-dimensional token embeddings, learned jointly with the model, block by block conditioned on the preceding text; warm-started from autoregressive weights at 3B and 8B parameters, it is competitive with masked diffusion and autoregressive baselines on mathematical reasoning and coding benchmarks, and its continuous trajectories degrade gracefully with few network evaluations and distil efficiently \citep{yang2026sigma}. Trained by likelihood, continuous models also scale \citep{yang2026continuous}. A continuous diffusion language model aligned with the architecture of modern discrete ones has a compute gap of about $20\times$ relative to autoregressive models, outperforms Duo with fewer parameters and masked diffusion when over-trained, and reaches a perplexity bound of 22.1 on OpenWebText, the best among continuous models. Choosing the noise schedule to minimise the variance of the evidence lower bound makes the cross-entropy, the information lost by time $t$, grow almost linearly in $t$, so that the difficulty of denoising is spread evenly over time, and learning the token embeddings by likelihood accounts for most of the gain.

A pure-death forward process, which decays every state to a single absorbing configuration such as a black image, gives another discrete construction \citep{santos2023blackout}. For graphs, discrete diffusion on edges and nodes \citep{vignac2023digress} and systems of SDEs on continuous relaxations \citep{jo2022score} are used.

\paragraph{Discrete flow matching.}

Generator matching on finite spaces gives discrete flow models, in which conditional rate matrices that realise a chosen probability path are averaged into a marginal rate matrix (\cref{sec:generator-matching}). Discrete flow models enable multimodal co-design of protein sequence and structure \citep{campbell2024generative}, and discrete flow matching uses mixture paths and probability velocities \citep{gat2024discrete}. As in the continuous case, the path can be chosen for its kinetic properties. A kinetic-optimal perspective selects velocities and paths that minimise a discrete kinetic energy, the analogue of Benamou--Brenier \citep{shaul2025flow}.

\paragraph{The simplex as a statistical manifold.}

Relaxing one-hot data to the probability simplex allows continuous dynamics. Categorical SDEs diffuse on the simplex \citep{richemond2022categorical}, Jacobi diffusions give Dirichlet diffusion score models \citep{avdeyev2023dirichlet}, Dirichlet probability paths give Dirichlet flow matching for DNA design \citep{stark2024dirichlet}, and diffusion on the simplex has also been studied directly \citep{floto2023diffusion}. Information geometry gives a canonical metric. The simplex with the Fisher--Rao metric is isometric, via $p\mapsto2\sqrt p$, to a positive orthant of a sphere, on which Riemannian flow matching with geodesic interpolants is simulation-free. Fisher flow matching uses this geometry, can bootstrap its flows with Riemannian optimal transport, and proves that the gradient flow it induces is optimal in reducing the forward relative entropy \citep{davis2024fisher}. Statistical flow matching formulates flow matching on statistical manifolds whose geodesics are available in closed form, instantiates it on categorical distributions through the square-root map onto the positive orthant of the sphere, and computes exact likelihoods at interior points of the simplex, while one-hot data require a variational bound \citep{cheng2024categorical}. Here the state space $\M$ is itself a space of probability measures, and the three geometries of \cref{fig:three} fold onto one another.

\paragraph{Theory.}

Convergence theory has been developed along the lines of the continuous one. Uniformisation gives exact simulation of the reverse chain and relative-entropy bounds for discrete diffusion on the hypercube $\{0,1\}^d$ with independent bit flips, using a number of jumps nearly linear in $d$ \citep{chen2024convergence}. A stochastic-integral formulation in terms of Poisson random measures provides a change-of-measure theorem for jump processes playing the role of Girsanov's, with which discrete and continuous diffusion can be analysed in parallel \citep{ren2025how}; high-order $\tau$-leaping solvers are analysed in the same framework \citep{ren2025fast}, and discrete-time analyses are also available \citep{zhang2025convergence}. Working with observables through adjoint equations yields bounds in any integral probability metric that are independent of the vocabulary size, for both masked and uniform priors \citep{kan2026vocabulary}. For masked diffusions with factorised denoisers that unmask positions in random order, the relative-entropy error caused by factorisation is bounded in terms of the number of tokens generated per step, independently of the sequence length; for worst-case orderings the bound grows linearly in that length. The optimal schedule sizes are determined by an information profile of the data distribution \citep{lavenant2025error}.

\section{Infinite dimensions}\label{sec:infinite}

The furthest departure from Euclidean space is to infinite dimensions, where the data are functions. Here the difficulty is not curvature but the absence of a translation-invariant reference measure. There is no Lebesgue measure, so densities, and with them scores, must be defined afresh. Gaussian measures supply the reference, the Cameron--Martin space supplies the admissible directions, and Malliavin calculus, which was born in infinite dimensions, supplies the formula for the score.

\paragraph{Gaussian measures and Cameron--Martin geometry.}

Functions, fields and solutions of PDEs are naturally elements of a separable Hilbert space $\cH$. A centred Gaussian measure $\cN(0,C)$ on $\cH$ has a self-adjoint, nonnegative, trace-class covariance $C$; its Cameron--Martin space is $\cH_C=C^{1/2}\cH$ with norm $\norm{C^{-1/2}h}$. Translating $\cN(0,C)$ by $h$ produces an equivalent measure if and only if $h\in\cH_C$, and by the Feldman--H\'ajek theorem two Gaussian measures are either equivalent or mutually singular \citep{bogachev1998gaussian,daprato2014stochastic}. \citet{hairer2009introduction} is a self-contained introduction to the SPDE side, including Gaussian measure theory and the Cameron--Martin space. The consequences for generative modelling are immediate. First, the noised data must stay in $\cH$. For an Ornstein--Uhlenbeck forward process $\dd X_t=-\tfrac12X_t\,\dd t+C^{1/2}\dd W_t$, with $W$ a cylindrical Wiener process on $\cH$, this requires a trace-class $C$, whereas a smoothing drift, as in the one-dimensional stochastic heat equation, only needs the stochastic convolution to be $\cH$-valued. Second, the laws at different noise levels are typically mutually singular. By the Feldman--H\'ajek theorem, $\cN(m_1,\sigma^2C)$ and $\cN(m_2,\tau^2C)$ are mutually singular whenever $\sigma\ne\tau$ and $C$ has infinite rank, so no common density exists across noise levels. Third, a score makes sense only along admissible directions, in the following sense.

\begin{definition}[Logarithmic derivative]\label{def:log-derivative}
A probability measure $\mu$ on $\cH$ is Fomin differentiable along $h\in\cH$ if there is $\beta_h\in L^1(\mu)$ with $\int\partial_h\varphi\,\dd\mu=-\int\varphi\,\beta_h\,\dd\mu$ for all $\varphi\in C^1_b(\cH)$. The function $\beta_h$ is the logarithmic derivative of $\mu$ along $h$; in finite dimensions $\beta_h=\inner{h}{\nabla\log p}$.
\end{definition}
For $\mu=\cN(0,C)$ and $h$ in the range of $C$, $\beta_h(y)=-\inner{y}{C^{-1}h}$ \citep{bogachev2010differentiable}.

\paragraph{Function-space diffusion and flow matching.}

Several groups developed diffusion models in function space at the same time. \citet{kerrigan2023diffusion} use Gaussian measures with trace-class covariance on Hilbert spaces, and denoising diffusion operators generalise denoising score matching to function space and sample with an infinite-dimensional Langevin dynamics \citep{lim2023score}. \citet{pidstrigach2024infinite} give a rigorous formulation, show how the choice of noise covariance matters and prove dimension-independent error bounds. Further developments include functional diffusion processes \citep{franzese2023continuous}, multilevel training \citep{hagemann2023multilevel}, conditional models for Bayesian linear inverse problems in infinite dimensions \citep{baldassari2023conditional}, and forward processes given by stochastic evolution equations \citep{lim2023evolution}. Practical architectures exploit resolution independence \citep{bondtaylor2024infinity}, spectral representations \citep{phillips2022spectral} and diffusion over function values \citep{dutordoir2023neural}.

Flow-based models transfer as well. Functional flow matching connects a Gaussian measure to the data on a Hilbert space \citep{kerrigan2024functional}, and functional rectified flow, built on the superposition principle in infinite dimensions, recovers functional flow matching without its restrictive measure-theoretic assumptions \citep{zhang2025flow}. Few-step flow maps transfer too. Under the assumptions of functional flow matching, the two-time flow map is Fr\'echet differentiable and satisfies the Eulerian equation of \cref{prop:flowmap}, which yields the identity of \cref{prop:meanflow} on a Hilbert space and a one-step generator started from a Gaussian measure; the same work shows that, unlike velocities, conditional two-time flows do not average to the marginal one \citep{li2025functional}. Function-space diffusion for physics-informed generation attains minimax-optimal rates for density estimation under suitable regularity conditions \citep{wang2025fundiff}, and neural-operator flows serve functional regression \citep{shi2024universal}.

\paragraph{Scores of SPDE laws via Malliavin calculus.}\label{sec:spde-scores}

The integration by parts of \cref{thm:malliavin-score} carries over to Hilbert-space-valued random variables with Malliavin derivatives taking values in Hilbert--Schmidt operators. For an $\cH$-valued $F$ with Malliavin derivative $\Phi_r=\D_rF$ and covariance operator $\gammaM=\int\Phi_r\Phi_r^*\dd r$, any covering field $v_h\in\Dom\Sk$ with $\int\Phi_rv_h(r)\dd r=h$ gives
\begin{equation}\label{eq:log-derivative-malliavin}
\begin{aligned}
\beta_h(F)&=-\E\big[\Sk(v_h)\,\big|\,F\big],\\
v_h(r)&=\Phi_r^*\,\gammaM^{\dagger}h\quad\text{(the canonical choice)},
\end{aligned}
\end{equation}
where $\beta_h$ is the logarithmic derivative of $\law(F)$ along $h$ (\cref{def:log-derivative}), $r$ ranges over $[0,T]$, each $\Phi_r$ is a Hilbert--Schmidt operator from the space of the driving noise into $\cH$, and $\gammaM^\dagger$ is the Moore--Penrose pseudoinverse. Two difficulties distinguish infinite dimensions. The operator $\gammaM$ is trace class, so when it has infinite rank $\gammaM^\dagger$ is unbounded and the canonical field is defined only for directions $h$ in its range; extensions to Cameron--Martin directions go through a density or regularisation argument. And for SPDEs the first variation need not be invertible, so the finite-dimensional factorisation through $Y_t^{-1}$ fails. \citet{mirafzali2025infinite} develop score-based diffusion in separable Hilbert spaces for forward processes given by linear SPDEs driven by space--time coloured noise with trace-class covariance, and obtain, through an infinite-dimensional Bismut--Elworthy--Li formula, a closed-form logarithmic derivative of the transition measure along Cameron--Martin directions. \citet{mirafzali2026logarithmic} extend the construction to variational SPDEs on Gelfand triples under fully local monotonicity, smoothness and non-degeneracy conditions, and build the Moore--Penrose covering field through a Tikhonov regularisation controlled by a trace formula. Their hypotheses are verified only in part for the two motivating examples. For the stochastic $p$-Laplacian, the drift-differentiability hypotheses hold for $p>3$ and the domination estimate for $p\ge4$, while the structural, smoothing, Galerkin-stability and directional non-degeneracy hypotheses remain assumptions. For the two-dimensional Navier--Stokes equation, logarithmic derivatives of finite-dimensional projections follow unconditionally from earlier results, and the full formula holds on a subspace of admissible directions whose denseness in $\cH$ is open. They also obtain integration-by-parts identities for scalar pairings $\inner{u(t)}{\varphi}$ of singular SPDEs within the Bruned--Chandra--Chevyrev--Hairer framework \citep{bruned2021renormalising}. \citet{greco2025malliavin} combine Malliavin and $\Gamma$-calculus with Wiener-chaos expansions and abstract time-reversal formulas to generalise score-based models to random fields. Together these results provide what function-space diffusion models need and the finite-dimensional theory cannot supply, which is a definition of the score intrinsic to the law of a solution and a formula for it.

\section{L\'evy and jump processes}\label{sec:levy}

The last relaxation concerns the noise rather than the space. Brownian motion is the natural driver for data whose fluctuations are light-tailed and continuous; for heavy-tailed data, and for data whose dimension can change, jump processes are more natural. Each piece of the Gaussian theory then needs a counterpart. The jump dynamics need a reversal, the jump part a score, and that score an integration by parts that produces it.

\paragraph{L\'evy and heavy-tailed drivers.}

Replacing Brownian motion by an isotropic $\alpha$-stable L\'evy process changes the generator to a fractional Laplacian. The exact time reversal is then a jump process whose intensity of jumping from $x$ to $x+z$ is reweighted by $p_t(x+z)/p_t(x)$, the continuous counterpart of \cref{prop:ctmc-reversal}. But a reverse-time SDE driven by the same stable noise, with a drift involving a fractional analogue of the score, has the same one-time marginals, which is what sampling requires. The L\'evy--It\^o model derives such a reverse-time SDE and a fractional denoising score matching, and reports faster sampling and better coverage of imbalanced data \citep{yoon2023score}. Denoising L\'evy probabilistic models instead replace the Gaussian noise of DDPM by $\alpha$-stable noise and, using a conditionally Gaussian representation of stable variables, reduce training to a modified DDPM with elementary proofs \citep{shariatian2025denoising}. Heavy tails can also be obtained without stable processes. Student-$t$ noise gives heavy-tailed diffusion and flow models for extreme events \citep{pandey2025heavy}. Memory, rather than heavy tails, is the point of fractional Brownian motion, and Markovian approximations of it give generative fractional diffusion models \citep{nobis2024generative}. Background on L\'evy processes is in \citet{applebaum2009levy}.

\paragraph{Jumps and changing dimension.}

Besides adding heavy tails, jumps can change the dimension of the state. Jump diffusions add or delete components, which enables generation across dimensions \citep{campbell2023trans}, and the time reversal of Markov processes with jumps under a finite-entropy condition provides the corresponding theory \citep{conforti2022time}. Generator matching (\cref{sec:generator-matching}) treats flows, diffusions and jumps in one regression framework.

\paragraph{Malliavin calculus with jumps.}

Malliavin calculus extends to L\'evy processes through chaos expansions with respect to Poisson random measures and a difference-operator derivative on the jump part \citep{dinunno2009malliavin}. An integration-by-parts representation of scores for jump-driven forward processes, extending \cref{thm:malliavin-score} and \cref{sec:spde-scores}, would give L\'evy-driven and SPDE-driven generative models the same footing that denoising score matching gives Gaussian ones, and it is among the open problems below.

\begin{openprob}[Open problems for Part V]
\begin{problem}
Learn heat kernels or their scores on general manifolds with guarantees, and quantify the effect of heat-kernel approximation errors on Riemannian diffusion. Stability bounds that convert the initialisation error and equation residual of a learned log heat kernel into log-density errors, and, for exact evolution, into score errors, are a first step \citep{ko2026riemannian}; end-to-end bounds for the resulting generative model are missing.
\end{problem}
\begin{problem}
Develop a convergence theory for discrete and jump-process models that matches the sharpness of the continuous theory. Specifically, determine the optimal dependence on dimension, vocabulary size and number of steps for the samplers used in practice, beyond the upper bounds of \cref{sec:discrete}, and extend such guarantees to L\'evy-driven models.
\end{problem}
\begin{problem}
Characterise the admissible directions and noise covariances for which function-space scores exist and can be estimated at discretisation-independent rates; combine \eqref{eq:log-derivative-malliavin} with neural operators.
\end{problem}
\begin{problem}
Derive Malliavin score formulas for L\'evy-driven SPDE forward processes, and for singular SPDEs beyond the scalar pairings of \cref{sec:spde-scores}, and test whether such forward processes improve generation of heavy-tailed or rough fields.
\end{problem}
\end{openprob}

\part{Physics of generative dynamics}

\section{Statistical physics of generation}\label{sec:statphys}

Diffusion models were born from nonequilibrium thermodynamics \citep{sohldickstein2015deep}, and physics has continued to shape how they are understood. A generative model is a stochastic system with a very large number of interacting degrees of freedom, and statistical physics is the discipline that predicts the collective behaviour of such systems from a few order parameters, without solving them exactly. Its typical answer is a phase diagram. As a control parameter varies, the behaviour of the system can change abruptly, and these changes are the events that matter for generation, such as the moment at which a sample commits to a class, or to a single training point. In the other direction, physicists have long worked with Gibbs measures on fields in very high dimension, and several of their constructions, such as Langevin dynamics on field space and the coarse-graining flow of the renormalisation group, are noising or denoising dynamics in all but name.

This part follows two threads, statistical physics and field theory. Statistical physics uses high-dimensional limits, phase transitions and fluctuation theorems to explain when generation works, when it memorises and what it costs. Field theory finds generative dynamics already present in stochastic quantisation, in the renormalisation group and in holography. One geometric idea runs through both threads. The time of a generative process behaves like a scale, and generation restores structure one scale at a time, coarse features before fine ones. The first thread and the links with stochastic quantisation and renormalisation rest on established results; the holographic correspondence of \cref{sec:holography} is more speculative, and we mark it as such.

\subsection{Phases of the reverse dynamics}

A sample is not created uniformly in time. Its features are decided in different, often narrow, windows of the reverse dynamics, and locating these windows tells us where the score must be accurate and when a model begins to reproduce its training data. In high dimension the windows become sharp, and they are best described as phase transitions.

\begin{figure}[t]
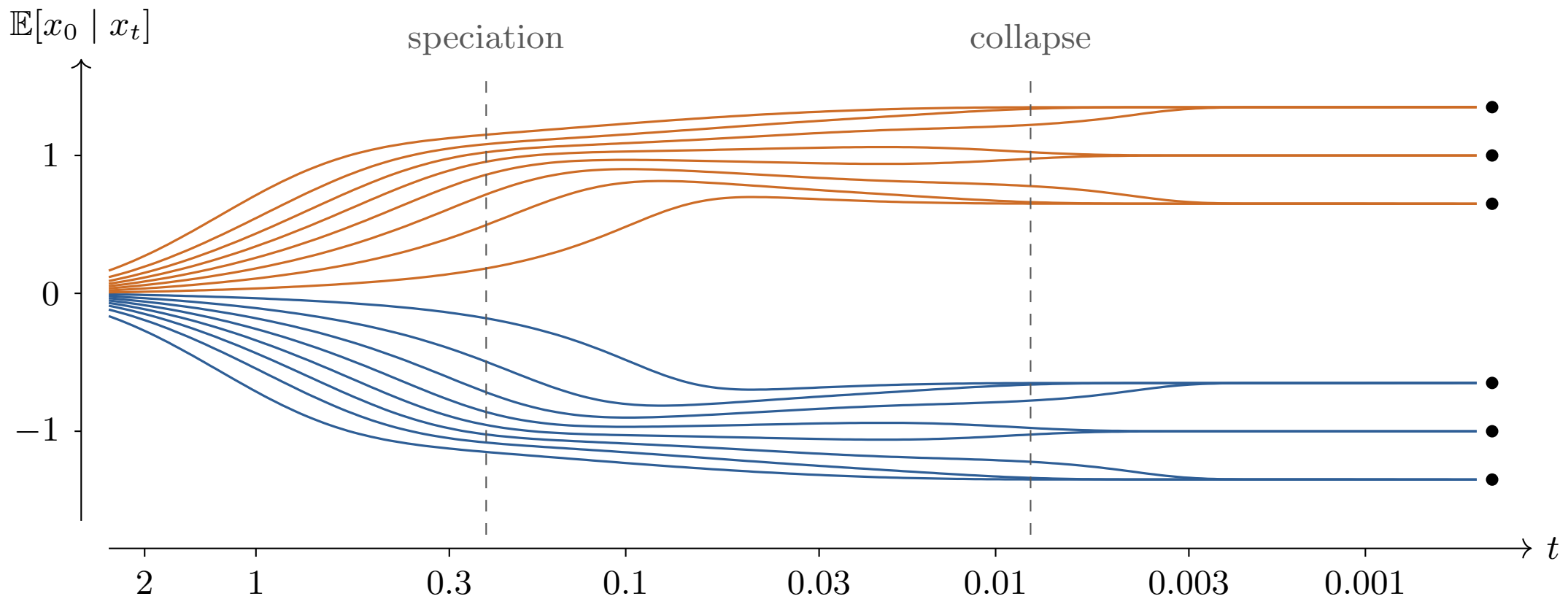

\centering

\caption{Speciation and collapse in one dimension, for a variance-preserving diffusion with $\alpha_t=e^{-t}$, $\sigma_t^2=1-e^{-2t}$ and the exact empirical score of six training points in two clusters (dots). Along probability-flow trajectories started at $t=2.5$ from evenly spaced quantiles of $\cN(0,1)$, the curves show the denoiser $\E[x_0\mid x_t]$ against decreasing time on a logarithmic axis, coloured by the cluster each trajectory reaches. The denoiser starts near the mean of the data, splits into the two cluster means, and only much later into the individual training points. The dashed lines mark the medians, over the sixteen trajectories, of the times at which the posterior probability of each trajectory's eventual cluster, and then of its eventual training point, first exceeds $0.95$, empirical counterparts of the speciation and collapse times of \citet{biroli2024dynamical}.}
\label{fig:speciation}
\end{figure}

Consider the exact empirical score of $n$ training points in dimension $d$, both large. \citet{biroli2024dynamical} show that the reverse dynamics then passes through three regimes, building on the symmetry-breaking transition that \citet{biroli2023generative} found in the Curie--Weiss model. Starting from pure noise, trajectories first undergo a transition that the authors call speciation, at which the gross structure of the data, such as the class of the eventual sample, is decided through a mechanism analogous to symmetry breaking; the speciation time can be read off a spectral analysis of the correlation matrix of the data. Later comes a collapse transition, at which trajectories are attracted to individual training points through a mechanism analogous to condensation in a glass. The collapse time is computed from an excess-entropy quantity, and its dependence on $n$ and $d$ characterises the curse of dimensionality. Avoiding collapse with the exact empirical score requires a number of samples growing exponentially with dimension. The same three regimes, with speciation as a second-order transition and collapse as the condensation of a random energy model, appear in discrete diffusion \citep{takahashi2026dynamical}. \Cref{fig:speciation} shows both transitions in one dimension.

The first transition can be computed by hand in the simplest case, and the computation shows why the language of symmetry breaking applies.

\begin{example}[Speciation as a pitchfork]\label{ex:pitchfork}
Idealise two tight clusters by the points $\pm m$, each with probability $\frac12$, and run the variance-preserving process of \cref{fig:speciation}, with $\alpha_t=e^{-t}$ and $\sigma_t^2=1-e^{-2t}$. Then $p_t$ is the equal mixture of $\cN(\alpha_tm,\sigma_t^2)$ and $\cN(-\alpha_tm,\sigma_t^2)$, and its score is
\begin{equation}\label{eq:pitchfork}
\partial_x\log p_t(x)=\frac{1}{\sigma_t^2}\Big(\alpha_tm\tanh\frac{\alpha_tmx}{\sigma_t^2}-x\Big).
\end{equation}
The critical points of $p_t$ solve $x=\alpha_tm\tanh(\alpha_tmx/\sigma_t^2)$. In the variable $y=x/(\alpha_tm)$ this is $y=\tanh(r_ty)$ with $r_t=\alpha_t^2m^2/\sigma_t^2$, which is $m^2$ times the signal-to-noise ratio $\alpha_t^2/\sigma_t^2$. This is the mean-field equation of the Curie--Weiss ferromagnet, with $r_t$ in the role of the inverse temperature. For $t>t_*=\frac12\log(1+m^2)$ we have $r_t<1$, the only solution is $y=0$, and $-\log p_t$ is a single well. For $t<t_*$ the origin becomes a local minimum of $p_t$ and two maxima appear, one for each cluster. In reverse time, $t_*$ is the moment at which the landscape splits and the trajectories are divided between the two wells. It depends on the data only through their variance $m^2$, a one-dimensional instance of reading speciation off the spectrum of the data covariance; the speciation time $t_S=\frac12\log\Lambda$ of \citet{biroli2024dynamical}, with $\Lambda$ the largest eigenvalue of the data covariance, is its leading order for $\Lambda\gg1$. The bifurcation is a pitchfork only because the clusters have equal weights. With weights $q$ on $m$ and $1-q$ on $-m$ the argument of the hyperbolic tangent is shifted by $\frac12\log\frac{q}{1-q}$, and the pitchfork unfolds into a fold, in line with the generic picture described below. For the clusters of \cref{fig:speciation}, whose means are $\pm1$, this gives $t_*=\frac12\log2\approx0.35$, a little before the dashed line, which marks the time at which the posterior probability of the eventual cluster exceeds $0.95$ rather than the bifurcation itself.
\end{example}

The second transition has an equally concrete description. By \cref{prop:empirical}, the exact empirical denoiser averages the training points with weights $w_i(x)\propto\exp(-\norm{x-\alpha_tx^{(i)}}^2/(2\sigma_t^2))$. These weights form a Boltzmann distribution over the $n$ training points, with energies $\frac12\norm{x-\alpha_tx^{(i)}}^2$ at temperature $\sigma_t^2$, and since the training points are random, so are the energies. As $t$ decreases the temperature falls, and collapse is the condensation of this Boltzmann distribution onto a single training point. The random energy model, the simplest model of condensation in spin-glass theory, is therefore the natural reference. In \cref{fig:speciation} the two transitions are well separated because the clusters are far apart compared with the spread within each of them.

Speciation is not tied to separated means. With classes formalised through Bayes classification, two classes merge when the gap between their average free entropies becomes comparable with its fluctuations, a criterion that involves only those two classes; it recovers the Gaussian-mixture results, applies to classes that differ only in higher-order or collective features, such as Ising chains at different temperatures, and predicts successive speciation times for nested classes \citep{achilli2026theory}. On compact Riemannian manifolds, speciation can be studied without taking the dimension large, as a bifurcation of the critical points of the evolving density \citep{marta2026theory}. Poincar\'e--Hopf and Morse theory then constrain the number and type of the equilibria of the score, and can force critical points of topological origin. For instance, when the critical points of a density on a closed surface are nondegenerate, the numbers of maxima, saddles and minima satisfy $\#\{\text{maxima}\}-\#\{\text{saddles}\}+\#\{\text{minima}\}=\chi$, the Euler characteristic, so on the torus, where $\chi=0$, every such density has at least two saddles, since it has at least one maximum and one minimum. For mixtures of heat kernels a generic speciation is a fold with a one-dimensional critical kernel, and the symmetric pitchfork is non-generic. To first order, an error in the score shifts the speciation time only through its component along the critical direction.

Other analyses reach the same picture by different routes. \citet{raya2023spontaneous} describe it from the side of dynamical systems, as an early phase of linear dynamics around a central fixed point, followed by a bifurcation into attractor dynamics towards the data, which suggests starting generation late from a Gaussian approximation and improves fast samplers. \citet{ambrogioni2023statistical} develops the statistical thermodynamics of these transitions, including critical instabilities. With hierarchical data, forward--backward experiments, which noise a sample to time $t$ and denoise it again, reveal a phase transition in the probability that the class changes, while low-level features vary smoothly, a signature of the hierarchical structure of the data \citep{sclocchi2025phase}. Finally, diffusion models with exact empirical scores behave as dense associative memories whose stored patterns are the training points \citep{hoover2023memory,ambrogioni2024search}. In trained models, the transition from retrieval to generalisation as the training set grows mirrors the behaviour of dense associative memories beyond their critical storage capacity, where new minima that differ from the stored patterns emerge \citep{pham2025memorization}.

The collapse analyses \citep{biroli2024dynamical,takahashi2026dynamical} and the associative-memory analyses of \citet{hoover2023memory} and \citet{ambrogioni2024search} assume the exact empirical score and therefore describe memorisation as the default, whereas speciation is already a property of the population score. They complement \cref{sec:generalisation}, where the trained score deviates from the empirical one, so that collapse is avoided as long as training stops before a memorisation time that grows linearly with $n$, and altogether once $n$ exceeds a model-dependent threshold \citep{bonnaire2025why}.

\subsection{Learning and hardness in high dimension}

Two questions decide whether diffusion models can succeed in high dimension. The statistical one is whether the score can be learned from the available data, and the computational one is whether the denoiser can be computed efficiently along the whole path. Exactly solvable models give sharp learning curves. For Gaussian mixtures, the velocity field of a flow-based model parameterised by a shallow denoiser can be analysed in the high-dimensional limit, and the analysis gives the generated distribution as a function of sample size \citep{cui2024analysis}. Deep networks can represent denoisers of high-dimensional graphical models efficiently by unrolling variational-inference denoising algorithms, which yields sample-efficient learning guarantees \citep{mei2025deep}.

Computational hardness is where phase transitions re-enter. Along the denoising path the posterior law of $x_0$ given $x_t$ sharpens as the signal-to-noise ratio grows, and the overlap of the denoiser $\E[x_0\mid x_t]$ with $x_0$ plays the role of an order parameter. In the statistical-physics description, the overlap reached by the Bayes-optimal denoiser is the global maximiser of a free-entropy potential, while an iterative algorithm started from an uninformative state, such as approximate message passing, reaches the first local maximiser it meets. When the transition is continuous, the two coincide. At a first-order transition the global maximiser jumps to another branch while the algorithm stays behind, and a sampler that needs the Bayes-optimal denoiser at that point receives a different one. For Gibbs measures of spin-glass type, first-order phase transitions encountered along the denoising path obstruct flow-, diffusion- and autoregressive-based samplers, and the tools of disordered systems map both the regions of hardness and the regions where these samplers succeed while Monte Carlo and Langevin dynamics fail \citep{ghio2024sampling}. Algorithmic stochastic localisation (\cref{sec:localisation}) turns the same picture into provable samplers when the denoiser can be computed by approximate message passing. Such guarantees hold for the Sherrington--Kirkpatrick model with inverse temperature $\beta<1/2$, in normalised Wasserstein distance, whereas no stable algorithm can succeed for $\beta>1$ \citep{elalaoui2022sampling}. \citet{celentano2024sudakov} extended the guarantees to the whole high-temperature regime $\beta<1$ by proving a local convexity property that \citet{elalaoui2022sampling} had conjectured. For spherical spin glasses, guarantees in total variation hold under a condition ensuring replica symmetry along the localisation path \citep{huang2024sampling}.

\subsection{Nonequilibrium thermodynamics}

A diffusion model runs an irreversible process backwards. Its forward process erases the structure of the data as it relaxes towards equilibrium, in the sense that the relative entropy of its law to the equilibrium law decreases to zero. The reverse process must rebuild that structure in finite time. Thermodynamics quantifies the cost of such operations through work, free energy and entropy production, and it supplies two exact identities that hold arbitrarily far from equilibrium.

Consider a system with energy $U_\lambda$, depending on a control parameter $\lambda$ unrelated to the log-SNR, in contact with a heat bath at inverse temperature $\beta$ and with dynamics satisfying detailed balance with respect to $e^{-\beta U_\lambda}$ at each fixed $\lambda$. Let it start in equilibrium at $\lambda_0$ and be driven to $\lambda_\tau$ by a protocol $(\lambda_t)_{t\in[0,\tau]}$ of finite duration $\tau$. Along a trajectory $(X_t)$ the protocol performs the work $\mathcal W=\int_0^\tau\partial_\lambda U_{\lambda_t}(X_t)\,\dot\lambda_t\,\dd t$. Jarzynski's equality and Crooks' fluctuation theorem state
\begin{equation}\label{eq:jarzynski}
\begin{aligned}
\E\big[e^{-\beta\mathcal W}\big]&=e^{-\beta\Delta F},\\
\frac{P_{\rm F}(\mathcal W)}{P_{\rm R}(-\mathcal W)}&=e^{\beta(\mathcal W-\Delta F)},
\end{aligned}
\end{equation}
where $\Delta F$ is the free-energy difference and $P_{\rm F},P_{\rm R}$ are the work distributions of the forward protocol and of the reversed protocol, the latter started in equilibrium at $\lambda_\tau$ \citep{jarzynski1997nonequilibrium,crooks1999entropy}. The expectation is over trajectories of the forward protocol, and $\Delta F=F_{\lambda_\tau}-F_{\lambda_0}$ with $F_\lambda=-\beta^{-1}\log\int e^{-\beta U_\lambda}$. Both identities rest on one relation between path measures. Let $\P_{\rm F}$ be the law of the trajectory under the forward protocol and $\tilde\P_{\rm R}$ the law of the trajectory under the reversed protocol, read backwards in time, so that both are laws on the same space of paths. Under the assumptions above, $\dd\tilde\P_{\rm R}/\dd\P_{\rm F}=e^{-\beta(\mathcal W-\Delta F)}$ path by path, where $\mathcal W$ is the work done along the path. Taking the $\P_{\rm F}$-expectation gives Jarzynski's equality, and taking the expectation of the logarithm gives
\begin{equation}\label{eq:dissipated-work}
\beta\big(\E[\mathcal W]-\Delta F\big)=\KL\big(\P_{\rm F}\,\|\,\tilde\P_{\rm R}\big)\ge0,
\end{equation}
which is the second law. The mean dissipated work is a relative entropy on path space, and it measures how distinguishable the forward process is from the reversed protocol read backwards in time. In stochastic thermodynamics the logarithm of the Radon--Nikodym derivative between a path measure and its time reversal is the entropy produced along the trajectory \citep{seifert2012stochastic}. Seen this way, the path-space relative entropies that train diffusion models (\cref{prop:elbo}) and the importance weights of annealed and controlled samplers (\cref{sec:control}) are dissipation functionals, and \eqref{eq:jarzynski} is the identity that makes annealed importance sampling, NETS and stochastic normalising flows unbiased \citep{neal2001annealed,albergo2025nets,caselle2022stochastic,wu2020stochastic}.

Two results make the analogy quantitative. \citet{ikeda2025speed} derive speed--accuracy relations for diffusion models. These inequalities bound the sensitivity of generation by the thermodynamic dissipation of the forward process, that is, the time integral of the temperature times the entropy production rate. The sensitivity is defined as the squared change in the $\W_1$ error divided by the $\chi^2$-divergence between the true and the used initial law of the reverse process and by the duration of the process. The main bound assumes an exactly learned velocity field, and a training-loss term is added otherwise. Minimising the bound selects forward processes along a Wasserstein-2 geodesic, and the authors compare cosine, conditional optimal-transport and optimal-transport schedules numerically. The argument behind this selection is short. Take an overdamped process $\dd X_t=f(t,X_t)\,\dd t+\sqrt{2T_t}\,\dd W_t$ at temperature $T_t$, with unit mobility and Boltzmann's constant set to one, whose law $\rho_t$ is transported, in the sense of \eqref{eq:continuity}, by its current velocity $v_t$, the probability current divided by $\rho_t$. Its rate of total entropy production is $\dot S_t=T_t^{-1}\int\abs{v_t}^2\,\dd\rho_t$. The dissipation rate $T_t\dot S_t$ is therefore the kinetic energy of $v_t$, which equals that of the curve $(\rho_t)$ in the Otto geometry of \cref{sec:transport} when $v_t$ is a gradient, as it is for the linear drifts of diffusion models, and exceeds it otherwise. \Cref{thm:bb}, after time is rescaled to $[0,\tau]$, gives
\begin{equation}\label{eq:dissipation-bound}
\int_0^\tau T_t\dot S_t\,\dd t=\int_0^\tau\!\!\int\abs{v_t}^2\,\dd\rho_t\,\dd t\ \ge\ \frac{\W_2^2(\rho_0,\rho_\tau)}{\tau},
\end{equation}
with equality exactly when $v_t$ is the gradient velocity of a constant-speed Wasserstein geodesic. The dissipation that limits accuracy is thus the action of the forward path in the geometry of \cref{sec:transport}, and the optimal-transport schedule is the forward process that dissipates least.

The second result turns entropy production into a measure of information. \citet{premkumar2024neural} define the neural entropy of a diffusion model, related to the total entropy produced by the forward diffusion and interpreted as the information the network must store to reverse it; measured on image models, it grows logarithmically with the number of training samples, which indicates very efficient compression. Information theory provides a third view. By the I-MMSE relation \citep{guo2005mutual}, the mutual information between the data and its noisy observation $\sqrt{\mathsf s}\,X_0+\epsilon$ at signal-to-noise ratio $\mathsf{s}$ has derivative $\frac12\,\mathrm{mmse}(\mathsf{s})$ in $\mathsf s$, where $\mathrm{mmse}(\mathsf s)=\E\norm{X_0-\E[X_0\mid\sqrt{\mathsf s}\,X_0+\epsilon]}^2$ is the error of the optimal denoiser. The entropy of discrete data is therefore $\frac12\int_0^\infty\mathrm{mmse}(\mathsf s)\dd\mathsf s$, and for data with a density the differential entropy is $\frac d2\log(2\pi e)-\frac12\int_0^\infty\big(\frac d{1+\mathsf s}-\mathrm{mmse}(\mathsf s)\big)\dd\mathsf s$, a comparison with the standard Gaussian, whose error is $d/(1+\mathsf s)$. \citet{kong2023information} give a pointwise version, which expresses the log-likelihood of each data point exactly through the pointwise error of the optimal denoiser and, with a learned denoiser in its place, refines the usual variational likelihood bounds of diffusion models. Wasserstein gradient flows also find their way back into physics as solvers for quantum many-body ground states \citep{neklyudov2023wasserstein}.

\section{Field theory, renormalisation and holography}\label{sec:fieldtheory}

Statistical physics studies generative models from outside. Field theory approaches them differently, through generative dynamics that were already present inside physics before diffusion models existed. A Euclidean lattice field theory is a Gibbs measure $Z^{-1}e^{-S[\phi]}\,\dd\phi$ on the configurations of a field, one real number per lattice site or one group element per link, with action $S$, partition function $Z$ and $\dd\phi$ the product of Lebesgue or Haar measures. Its correlation functions are expectations under this measure. Computing them is a sampling problem of the kind generative models address, with three features that make it both hard and instructive. The dimension grows with the volume of the lattice; near a critical point correlations extend over many lattice spacings, so that local updates mix slowly; and the measures carry exact symmetries, global or local. In the continuum, the quadratic part of the action defines a Gaussian measure on a space of distributions, so the Gaussian geometry of \cref{sec:infinite} reappears in field theory.

Stochastic quantisation samples a field theory with a Langevin equation, the renormalisation group coarse-grains a field by a heat flow, and holography organises a field theory by a radial coordinate that behaves like a scale. We take each correspondence in turn and close the section by placing the three side by side.

\subsection{Stochastic quantisation and lattice field theory}

\paragraph{Stochastic quantisation.}

Parisi and Wu observed that a Euclidean field theory with action $S[\phi]$ can be quantised by introducing a fictitious time $\tau$ and running the Langevin equation
\begin{equation}\label{eq:parisi-wu}
\begin{aligned}
\partial_\tau\phi(x,\tau)&=-\frac{\delta S}{\delta\phi(x,\tau)}+\eta(x,\tau),\\
\E[\eta(x,\tau)\eta(x',\tau')]&=2\delta(x-x')\delta(\tau-\tau'),
\end{aligned}
\end{equation}
whose stationary distribution is $e^{-S[\phi]}$ \citep{parisi1981perturbation,damgaard1987stochastic}. Here $x$ ranges over Euclidean space-time, $\eta$ is a centred Gaussian white noise, and $\delta S/\delta\phi(x,\tau)$ is the functional derivative of $S$ at $\phi(\cdot,\tau)$. This is the Langevin diffusion of \eqref{eq:fp-gradflow} in field space. For the free field of mass $m>0$, with $S[\phi]=\frac12\int(\abs{\nabla\phi}^2+m^2\phi^2)\,\dd x$, the drift is linear and each Fourier mode is an independent Ornstein--Uhlenbeck process,
\begin{equation}\label{eq:free-sq}
\partial_\tau\hat\phi(k,\tau)=-(k^2+m^2)\,\hat\phi(k,\tau)+\hat\eta(k,\tau),
\end{equation}
where $\hat\phi(k,\tau)$ and $\hat\eta(k,\tau)$ are Fourier transforms in $x$ at wave vector $k$ and $k^2=\abs{k}^2$. Each mode relaxes at rate $k^2+m^2$ to its stationary variance $(k^2+m^2)^{-1}$, the free propagator. The slowest mode relaxes at rate $m^2=\xi^{-2}$, where $\xi=1/m$ is the correlation length, so the relaxation time of this local dynamics grows like $\xi^2$ as the critical point $m\to0$ is approached. This is critical slowing down in its simplest form, and it explains the appeal of proposals that update all scales at once.

A diffusion model trained on field configurations replaces the action gradient by a learned score of noised configurations. \citet{wang2024quantization} make the correspondence explicit and show that a diffusion model trained on hybrid Monte Carlo configurations can serve as a global Metropolis--Hastings proposal for two-dimensional $\phi^4$ theory, with likelihoods from the probability-flow ODE, and that this proposal substantially reduces autocorrelation times near criticality, where local Markov chain Monte Carlo suffers critical slowing down. Follow-up work studies how diffusion models learn higher-order cumulants \citep{aarts2024learning} and distributions generated by complex Langevin dynamics for theories with complex actions \citep{habibi2024diffusion}.

\paragraph{Flows for lattice field theory.}

Normalising flows entered lattice field theory as proposal distributions for Markov chain Monte Carlo. If a model with tractable density $q$ proposes a configuration $\phi'$ independently of the current one $\phi$, accepting it with probability $\min\{1,e^{-S[\phi']}q(\phi)/(e^{-S[\phi]}q(\phi'))\}$ leaves the Gibbs measure invariant for any $q$. Provided $q$ is positive wherever $e^{-S}$ is, the quality of $q$ affects only the acceptance rate, and with it the autocorrelation time. An independence Metropolis sampler with a flow proposal is asymptotically exact, and a flow trained to a fixed acceptance rate shows no critical slowing down in the sampling stage, the difficulty being shifted to training \citep{albergo2019flow,albergo2021introduction}. 

Gauge fields bring geometry of their own. On a lattice, a gauge field is a discrete connection. It assigns to each oriented link from $x$ to $x+\hat\mu$, with $\hat\mu$ the unit lattice vector in direction $\mu$, an element $U_\mu(x)$ of a compact Lie group $G$, the parallel transport along that link. A gauge transformation, an arbitrary map $\Omega$ from sites to $G$, acts by $U_\mu(x)\mapsto\Omega(x)U_\mu(x)\Omega(x+\hat\mu)^{-1}$. The action is built from the traces of the holonomies around the smallest closed loops, the plaquettes, and these traces are gauge invariant. The symmetry group is thus a copy of $G$ at every site. Gauge theories require flows equivariant under local gauge transformations, and these have been built for $U(1)$ \citep{kanwar2020equivariant}, for $SU(N)$ \citep{boyda2021sampling} and in arbitrary space-time dimension \citep{abbott2023normalizing}, while continuous equivariant flows scale further for scalar field theories \citep{dehaan2021scaling,gerdes2023learning}. Reweighting with the exact sampling density of a generative neural sampler, such as a flow or an autoregressive network, gives asymptotically unbiased estimates of observables, including thermodynamic ones such as free energies and entropies \citep{nicoli2020asymptotically}. Stochastic normalising flows are naturally read as nonequilibrium transformations obeying Jarzynski's equality \citep{caselle2022stochastic}. Lattice field theory is thus a demanding test bed in which the symmetries of \cref{sec:lie} are local, the target is known up to normalisation as in \cref{sec:control}, and exactness matters.

\subsection{Renormalisation as transport and coarse-graining as noising}

Physical fields near a critical point, like natural images, have structure on every scale, and the renormalisation group is the systematic way of treating such systems one scale at a time. Integrating out fluctuations at short distances leaves an effective description of the long-distance behaviour that depends on the microscopic details only through a few relevant parameters, which is why systems as different as a fluid at its critical point and a uniaxial magnet share the same critical exponents. For generative modelling the point is that, under coarse-graining, the Boltzmann factor of the interaction evolves by a heat flow, that is, by a noising process.

Wilson's renormalisation group integrates out short-distance fluctuations \citep{wilson1974renormalization}. \citet{polchinski1984renormalization} writes the exact renormalisation group, which he credits to Wilson, as a differential equation for the interaction, chosen so that the partition function does not change as the cutoff is lowered, and uses it to prove perturbative renormalisability. Its solution is a Gaussian convolution. With $C_\Lambda$ the regularised propagator,
\begin{equation}\label{eq:polchinski}
e^{-V_{\Lambda}(\phi)}=\int e^{-V_{\Lambda_0}(\phi+\chi)}\,\dd\mu_{C_{\Lambda_0}-C_{\Lambda}}(\chi),
\end{equation}
where $\mu_C$ is the centred Gaussian measure with covariance $C$, $V_{\Lambda_0}$ is the bare interaction at the initial cutoff $\Lambda_0$, and $\Lambda\le\Lambda_0$. Here $C_\Lambda$ is the free covariance with the modes above $\Lambda$ suppressed, for instance, in Fourier variables, $C_\Lambda(k)=K(k^2/\Lambda^2)/(k^2+m^2)$ for a smooth nonincreasing cutoff function $K$ equal to $1$ near $0$ and decaying to $0$ at infinity, so that $C_{\Lambda_0}-C_\Lambda$ is nonnegative and carries the modes between $\Lambda$ and $\Lambda_0$. The formula has a transparent probabilistic meaning. Split the free field into independent Gaussian fields, $\phi=\phi_\Lambda+\chi$, with covariances $C_\Lambda$ and $C_{\Lambda_0}-C_\Lambda$, the slow and the fast modes. Integrating out $\chi$ from $e^{-V_{\Lambda_0}(\phi)}\mu_{C_{\Lambda_0}}(\dd\phi)$ leaves for $\phi_\Lambda$ the law proportional to $e^{-V_\Lambda(\phi_\Lambda)}\mu_{C_\Lambda}(\dd\phi_\Lambda)$, so $V_\Lambda$ is the effective interaction of the slow modes. As the cutoff is lowered, \eqref{eq:polchinski} convolves $e^{-V_{\Lambda_0}}$ with Gaussian measures of growing covariance $C_{\Lambda_0}-C_\Lambda$. Acting on $e^{-V_{\Lambda_0}}$, viewed as an unnormalised density on field space, this is a heat semigroup, in other words a variance-exploding forward process, and the Polchinski equation is the PDE satisfied by minus its logarithm. Explicitly, let the cutoff $\Lambda(t)$ decrease from $\Lambda_0$, and set $\Gamma_t=C_{\Lambda_0}-C_{\Lambda(t)}$ and $V_t=V_{\Lambda(t)}$. Then $Z_t=e^{-V_t}$, a function on field space rather than a normalising constant, solves the heat equation $\partial_tZ_t=\frac12\tr(\dot\Gamma_t\nabla^2Z_t)$ on field space, with diffusion matrix $\dot\Gamma_t$, and
\begin{equation}\label{eq:polchinski-pde}
\partial_tV_t=\tfrac12\tr\big(\dot\Gamma_t\nabla^2V_t\big)-\tfrac12\inner{\nabla V_t}{\dot\Gamma_t\nabla V_t},
\end{equation}
where $\nabla$ denotes the functional derivative. The structure is the same as for any heat flow. If $p_t=e^{u_t}$ solves $\partial_tp=\Delta p$, so that $u_t=\log p_t$ is a log-density rather than a velocity, then
\begin{equation}\label{eq:viscous-hj}
\partial_tu_t=\Delta u_t+\abs{\nabla u_t}^2,
\end{equation}
a viscous Hamilton--Jacobi equation for the log-density whose gradient is the score. Equation \eqref{eq:polchinski-pde} is \eqref{eq:viscous-hj} for $u_t=-V_t$, with the Euclidean diffusion matrix $2\Id$ replaced by $\dot\Gamma_t$, so that in this dictionary the effective force $-\nabla V_t$ of the running interaction plays the role of the score. For the law of the field itself, the same coarse-graining is variance-preserving. Suppose that $K>0$, so that $\alpha_t^2=C_{\Lambda(t)}/C_{\Lambda_0}$ is a Fourier multiplier with values in $(0,1]$, and let $\sigma_t^2=1-\alpha_t^2$. Then the rescaled slow field $\alpha_t^{-1}\phi_{\Lambda(t)}$ has the law of $\alpha_t\phi+\sigma_t\zeta$, where $\phi$ follows the interacting measure and $\zeta\sim\mu_{C_{\Lambda_0}}$ is an independent free field. This is a variance-preserving forward process whose schedule depends on the Fourier mode, with the free field in the role of the noise.

\citet{cotler2023renormalization} show, by formal functional calculus, that Polchinski's flow (more generally, Wegner--Morris flow with a quadratic seed action) is the optimal-transport (Wasserstein) gradient flow of the relative entropy to a Gaussian reference functional. A regularised version of that relative entropy is then an RG monotone, and a JKO-type implicit discretisation turns RG into a sequence of variational problems. \citet{cotler2023renormalizing} then train diffusion models to learn inverse RG flows, interpret multiscale diffusion physically and build adaptive bridge samplers for lattice theories. Renormalisation-group diffusion models generate data coarse-to-fine by reversing an RG flow, with applications to protein structures and images \citep{masuki2025generative}. Using an exact Polchinski flow as the probability path of flow matching makes the velocity a functional derivative of the running interaction \citep{masuki2026renormalization}. For local, or conditionally local, data, and under the assumption that the running interaction stays local on the scale $\Lambda^{-1}$ along the flow, a strengthened form of the quasi-locality generally assumed in exact RG analyses, this velocity is quasi-local. Restricted to a patch with a buffer of width $\ell$, the velocity at cutoff scale $\Lambda$ has a root-mean-square error that decays like $e^{-c\Lambda\ell}$, up to a factor growing at most polynomially in the system size. If the local velocity is moreover Lipschitz uniformly in the scale, an assumption the authors do not prove, patches with buffers of width of order $\Lambda^{-1}(\log L_{\rm sys}+\log(1/\varepsilon))$, where $L_{\rm sys}$ is the linear size of the system, give $\W_2$ accuracy $\varepsilon$. Generation then costs time nearly linear in the volume and, in one-dimensional tests, still reproduces correlations far beyond a single patch.

Multiscale generative models predate this vocabulary. Several of them rest on a factorisation over scales. If, starting from a field $x_0$, an orthogonal wavelet transform maps $x_{j-1}$ to a coarser field $x_j$ and the details $\bar x_j$ lost in between, for $j=1,\dots,J$, then
\begin{equation}\label{eq:scale-factorisation}
p(x_0)=p(x_J)\prod_{j=1}^Jp(\bar x_j\mid x_j),
\end{equation}
and a field can be generated from the coarsest scale down, one conditional at a time. Each conditional describes fluctuations at a single scale given the coarser scales. Neural-network renormalisation group builds a hierarchical bijection that disentangles scales \citep{li2018neural}, and RG-Flow combines it with sparse priors for interpretable latent variables \citep{hu2022rg}. Wavelet conditional renormalisation models a field through conditional probabilities of wavelet coefficients scale by scale \citep{marchand2023wavelet}. Wavelet score-based generative models explain the acceleration of multiscale generation by the ill-conditioning of the score; by factorising across scales they use the same number of time steps at every scale, so that cost grows only linearly with image size. The authors prove this for stationary Gaussian processes with power-law spectra, where normalised wavelet coefficients have a dimension-independent condition number, and support it empirically for the critical $\varphi^4$ field and natural images \citep{guth2023wavelet}. Learned local Markov models of wavelet coefficients conditioned on the next coarser scale, with a global model only for the coarsest band, capture face images using conditioning neighbourhoods of fixed size independent of the image size \citep{kadkhodaie2023learning}. The correspondence between neural networks and field theories, in which network ensembles define field theories with computable correlators, gives a complementary dictionary \citep{halverson2021neural}.

\subsection{Holography}\label{sec:holography}

This subsection is more speculative than the rest of the survey. The correspondence it develops between radial depth in anti-de Sitter space and diffusion time rests on formal structure and on small-scale experiments, and we present it as a direction for research rather than as an established result.

The AdS/CFT correspondence equates a quantum gravity theory on asymptotically anti-de Sitter space $\mathrm{AdS}_{d+1}$ with a conformal field theory on its $d$-dimensional boundary \citep{maldacena1997large,gubser1998gauge,witten1998anti}; \citet{hubeny2015ads} gives a short review. In Euclidean Poincar\'e coordinates the metric is
\begin{align*}
\dd s^2=\frac{L^2}{z^2}\big(\dd z^2+\dd x^2\big),\qquad z>0,\ x\in\R^d,
\end{align*}
which is the hyperbolic space $\mathbb H^{d+1}$ of radius $L$, the AdS radius; in Lorentzian signature the constant-time slices of AdS are hyperbolic spaces. This geometry already encodes scale. The dilation $(z,x)\mapsto(\lambda z,\lambda x)$ is an isometry, so moving into the bulk by a factor $\lambda$ is the same as viewing the boundary at a scale larger by $\lambda$; more generally, the conformal transformations of the boundary are realised by isometries of the bulk. A geodesic joining two boundary points at distance $\ell$ is a semicircle that reaches depth $\ell/2$, so structure of size $\ell$ on the boundary is represented at depth of order $\ell$, and the volume element $(L/z)^{d+1}\dd z\,\dd x$ places most of the volume near the boundary, where the fine scales live.

\begin{figure}[t]
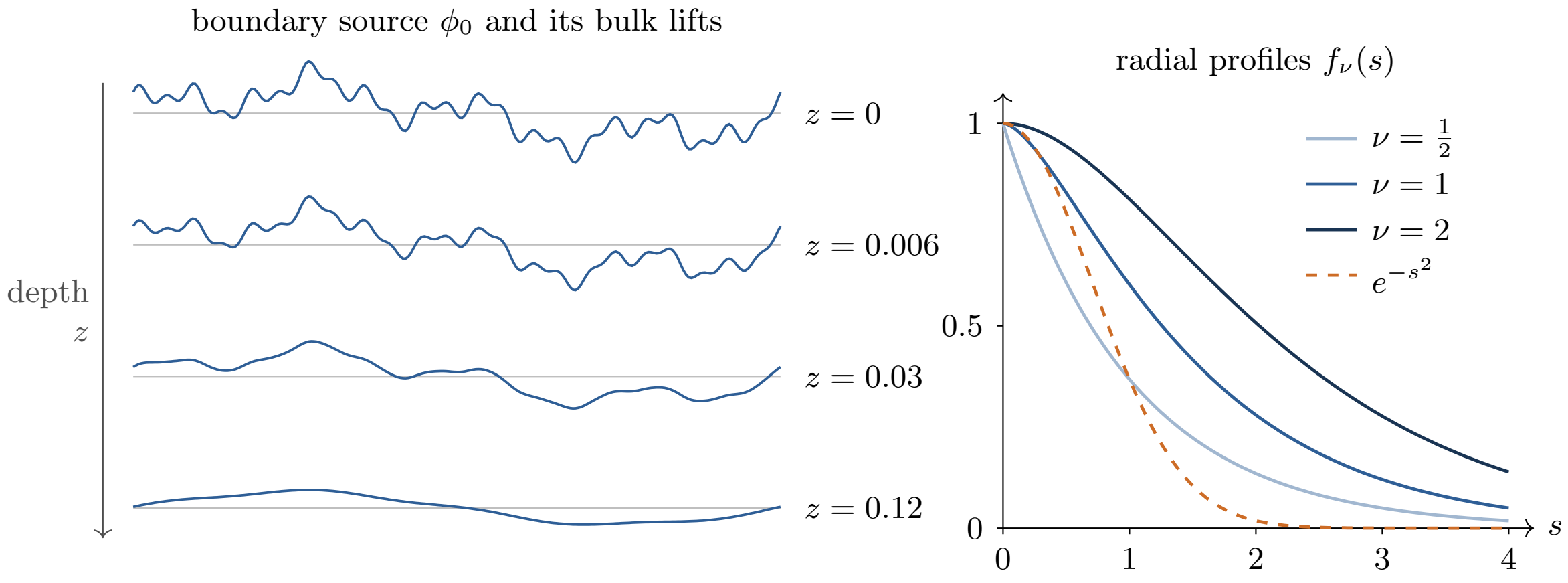

\centering

\caption{Depth in the bulk acts as a coarse-graining scale. In the left panel, a boundary source $\phi_0$ made of four Fourier modes is carried into the bulk by \eqref{eq:bulk-lift} with $\nu=1$, drawn with the common factor $z^{d-\Delta}$ removed. The horizontal axis of this panel is the boundary coordinate $x\in[0,1]$, and the modes have momenta $k=2\pi n$ with $n=1,4,11,23$. Each mode of momentum $k$ is multiplied by $f_\nu(kz)$, so the finest modes disappear first as the depth $z$ grows, and at $z=0.12$ little beyond the coarsest mode survives. In the right panel, the radial profiles $f_\nu(s)$ for three values of $\nu$ equal $1$ at $s=0$ and decay once $s\gtrsim1$, like the heat-kernel factor $e^{-k^2t}$ written in the variable $s=k\sqrt t$. For $\nu=\tfrac12$ the profile is exactly $e^{-s}$.}
\label{fig:bulk-lift}
\end{figure}

A free bulk scalar of mass $m$ satisfies the Klein--Gordon equation $(-\Box_{\rm AdS}+m^2)\phi=0$ and behaves near the boundary as $\phi(z,x)\simeq z^{d-\Delta}\phi_0(x)+\dots+z^{\Delta}A(x)$ with $\Delta(\Delta-d)=m^2L^2$. Here $\Box_{\rm AdS}$ is the Laplace--Beltrami operator of the metric above, $\phi_0$ and $A$ are functions on the boundary, and $\Delta>\frac d2$ is the scaling dimension, a number rather than the Laplacian. For a Fourier mode $e^{\mathrm ik\cdot x}$ the equation reads $z^2\partial_z^2\phi-(d-1)z\partial_z\phi=(\abs{k}^2z^2+m^2L^2)\phi$. Near $z=0$ the term $\abs{k}^2z^2$ is negligible, and substituting $\phi=z^\alpha$ gives $\alpha(\alpha-d)=m^2L^2$, with roots $\Delta$ and $d-\Delta$. Since $\Delta>\frac d2$, the power $z^{d-\Delta}$ dominates as $z\to0$, which is why its coefficient $\phi_0$ serves as the boundary value. The boundary value $\phi_0$ is the source of a CFT operator $\mathcal O$ of dimension $\Delta$, and the dictionary $\big\langle\exp\int\phi_0\mathcal O\big\rangle_{\rm CFT}=Z_{\rm bulk}[\phi\to\phi_0]$ relates boundary correlators to bulk dynamics, with the integral over the boundary, the expectation taken in the boundary theory and $Z_{\rm bulk}[\phi\to\phi_0]$ the bulk partition function with boundary value $\phi_0$ \citep{gubser1998gauge,witten1998anti}. In Fourier space the regular bulk solution with source $\phi_0$ is
\begin{equation}\label{eq:bulk-lift}
\begin{aligned}
\phi(z,k)&=\phi_0(k)\,z^{d-\Delta}\,f_\nu(kz),\\
f_\nu(s)&=\frac{2^{1-\nu}}{\Gamma(\nu)}\,s^{\nu}K_\nu(s),\qquad \nu=\Delta-\tfrac d2>0,
\end{aligned}
\end{equation}
with $K_\nu$ the modified Bessel function of the second kind and $k=\abs{k}$ in $f_\nu(kz)$ and below; $f_\nu(s)\to1$ as $s\to0$ and decays like $s^{\nu-1/2}e^{-s}$ as $s\to\infty$. The lift therefore suppresses boundary modes of momentum $k$ once $z\gtrsim1/k$, so depth in the bulk acts as a coarse-graining scale (\cref{fig:bulk-lift}), much as the heat-kernel factor $e^{-k^2t}$ suppresses high frequencies in the heat-dissipation and blurring models of \cref{sec:diffusion}. In a standard diffusion model, additive white noise has a similar effect on natural images, whose weak high-frequency components are the first to be drowned. In the geometric language above, a boundary mode of momentum $k$, of wavelength $2\pi/k$, penetrates the bulk to a depth of order $1/k$, just as a geodesic anchored at boundary points a distance $\ell$ apart reaches depth $\ell/2$.

This radial-scale correspondence is the holographic renormalisation group. The radial evolution of the on-shell action obeys a Hamilton--Jacobi equation that encodes the RG flow of the boundary theory \citep{deboer2000holographic}. With a cutoff in the bulk, the radial evolution of boundary conditions parallels Wilsonian renormalisation, although which field-theory cutoff corresponds to the radial cutoff is left open \citep{heemskerk2011holographic}. Holographic renormalisation supplies the counterterms that make the dictionary finite \citep{skenderis2002lecture}. Tensor networks give a discrete version. The multiscale entanglement renormalisation ansatz \citep{vidal2007entanglement,vidal2008class} has the geometry of a discretised hyperbolic space, which suggested that holography is a form of entanglement renormalisation \citep{swingle2012entanglement}. On the hyperbolic slices themselves, heat kernels are explicit. They are elementary in odd dimensions, as \cref{tab:kernels} records for $\mathbb H^3$, and a one-dimensional integral in even dimensions \citep{grigoryan1998heat}, so hyperbolic forward processes come with closed-form or quadrature-computable kernels and scores (\cref{sec:riemannian}).

Machine learning met holography well before GenAdS \citep{mirafzali2026holographic}, the holographic generative flow described below. \citet{hashimoto2018deep} identified the radial evolution of a bulk field on a discretised AdS geometry with a deep network whose weights encode the metric, and learned the metric from boundary data; the method was applied to holographic QCD \citep{hashimoto2018qcd}, and AdS/CFT was cast as a deep Boltzmann machine \citep{hashimoto2019ads}. \citet{you2018machine} learned spatial geometry from entanglement features with random tensor networks. A flow-based hierarchical generative model trained as an exact holographic mapping already produces an emergent bulk geometry matching $\mathbb H^3$ \citep{hu2020machine}.

\paragraph{Holographic generative flows.}

GenAdS \citep{mirafzali2026holographic} makes the bulk direction the flow parameter of a generative model. Boundary samples are encoded in Fourier space as scalar sources $\phi_0$ and carried into the bulk by Klein--Gordon dynamics, with flow matching staged in a normalised radial phase space and a neural network learning a residual correction to the Klein--Gordon backbone; physical parameters such as the bulk mass, which sets $\Delta$ and hence the radial profile in \eqref{eq:bulk-lift}, become design choices. In the reported experiments, most of the gain on a two-dimensional checkerboard benchmark comes from the Fourier representation and the convolutional architecture, while on MNIST the holographic encoding trades a little fidelity for coverage, its linear-path models reaching significantly higher recall than a convolutional baseline at comparable precision and slightly worse Fr\'echet statistics, whereas a Hermite path collapses recall. For the theory of Part~II, the construction raises the question of which properties of a radial, scale-ordered transport make the velocity field easier to learn.

\subsection{Scale as an extra dimension}

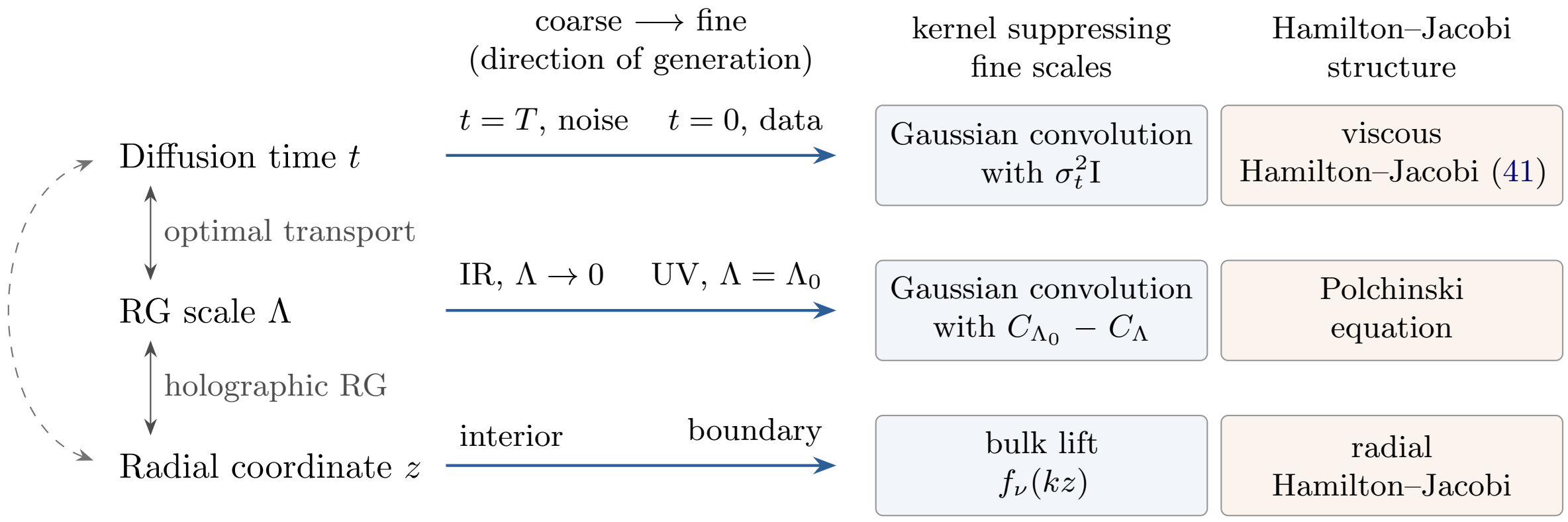
\begin{figure}[t]
\centering
\begin{tikzpicture}[>=Stealth,
  lab/.style={font=\footnotesize},
  head/.style={font=\footnotesize,align=center,anchor=south},
  par/.style={font=\small,anchor=west},
  cell/.style={draw=black!40,rounded corners=2pt,font=\footnotesize,align=center,text width=#1,minimum height=1.0cm,inner sep=3pt}]
\foreach \y/\name/\coarse/\fine/\ker/\hj in {
  3.1/{Diffusion time $t$}/{$t=T$, noise}/{$t=0$, data}/{Gaussian convolution\\ with $\sigma_t^2\Id$}/{viscous\\ Hamilton--Jacobi \eqref{eq:viscous-hj}},
  1.55/{RG scale $\Lambda$}/{IR, $\Lambda\to0$}/{UV, $\Lambda=\Lambda_0$}/{Gaussian convolution\\ with $C_{\Lambda_0}-C_\Lambda$}/{Polchinski\\ equation},
  0.0/{Radial coordinate $z$}/{interior}/{boundary}/{bulk lift\\ $f_\nu(kz)$}/{radial\\ Hamilton--Jacobi}}{
  \node[par] at (0,\y) {\name};
  \draw[->,thick,figblue] (3.4,\y) -- (7.3,\y);
  \node[lab,anchor=south west] at (3.4,\y+0.06) {\coarse};
  \node[lab,anchor=south east] at (7.3,\y+0.06) {\fine};
  \node[cell=3.1cm,fill=figblue!6] at (9.35,\y) {\ker};
  \node[cell=3.2cm,fill=figorange!8] at (12.85,\y) {\hj};
}
\node[head] at (5.35,3.75) {coarse $\longrightarrow$ fine\\ (direction of generation)};
\node[head] at (9.35,3.75) {kernel suppressing\\ fine scales};
\node[head] at (12.85,3.75) {Hamilton--Jacobi\\ structure};
\draw[<->,black!70] (0.45,2.8) -- (0.45,1.85) node[midway,right,lab] {optimal transport};
\draw[<->,black!70] (0.45,1.25) -- (0.45,0.3) node[midway,right,lab] {holographic RG};
\draw[<->,black!55,dashed] (-0.12,3.05) to[out=200,in=160] (-0.12,0.05);
\end{tikzpicture}
\caption{Three parameters that order a system by scale. Each row runs from coarse to fine, which is the direction of generation, and in each a Gaussian or Bessel kernel suppresses fine scales while a Hamilton--Jacobi equation governs a log-density or an on-shell action. The optimal-transport formulation of RG makes the identification of diffusion time with the RG scale precise at the level of formal functional calculus \citep{cotler2023renormalization}, and holographic renormalisation relates the radial coordinate to the RG scale \citep{deboer2000holographic,heemskerk2011holographic}. A direct identification of radial depth with diffusion time, drawn dashed, remains open.}
\label{fig:scales}
\end{figure}

Three parameters now play analogous roles (\cref{fig:scales}). Diffusion time orders the generative process from coarse to fine, the renormalisation scale orders a field theory from infrared to ultraviolet, and the holographic radial coordinate orders the bulk from deep interior to boundary. In each case a Gaussian or Bessel kernel suppresses fine scales first, and a Hamilton--Jacobi structure governs the evolution of a log-density or on-shell action, through \eqref{eq:viscous-hj} for diffusion, Polchinski's equation for exact RG and the radial Hamilton--Jacobi equation of holographic RG. The optimal-transport formulation of RG \citep{cotler2023renormalization} makes the identification of the first two precise at the level of formal functional calculus. Making the third precise remains open. It would mean relating learned flows in AdS to holographic RG, learning bulk metrics as forward processes in the spirit of \citet{hashimoto2018deep}, and exploiting the explicit hyperbolic heat kernel.

For generative modelling, the three rows of \cref{fig:scales} offer a way of reading the time variable. A forward process is a choice of the order in which the structure of the data is forgotten. All three constructions forget fine structure first and restore it last (for additive noise, on data whose spectra decay as those of natural images do), so that generation proceeds from coarse to fine. The phases of \cref{sec:statphys} fit the same reading. Speciation fixes the coarse structure of a sample early, and collapse onto a training point, when it happens, is the last and finest decision.

\begin{openprob}[Open problems for Part VI]
\begin{problem}
Extend the speciation and collapse analysis from exact empirical scores to trained networks with the inductive biases of \cref{sec:generalisation}, beyond random-features models, and describe how the collapse time of the reverse dynamics moves as training passes the memorisation time $\tau_{\rm mem}$.
\end{problem}
\begin{problem}
Turn the speed--accuracy relations of stochastic thermodynamics into schedule-design principles with guarantees. Specifically, bound the time-discretisation error of \cref{sec:theory} by the dissipation of the forward process, or show that no such bound holds, and decide whether minimum-dissipation schedules minimise the number of steps for a given accuracy.
\end{problem}
\begin{problem}
Derive the radial flow of GenAdS from holographic renormalisation, and test whether bulk geometries other than pure AdS (for example with a horizon) induce useful forward processes; the hyperscaling-violating geometries interpolating between AdS and flat space tested so far, on a two-dimensional checkerboard and without matching the bulk mass, underperform planar AdS \citep{mirafzali2026holographic}.
\end{problem}
\begin{problem}
Use RG-as-transport to design forward processes adapted to the correlation structure of physical fields, with provable step-count savings beyond the Gaussian case.
\end{problem}
\end{openprob}

\part{Synthesis}

\section{A unified view}\label{sec:zoo}

Seen from a distance, the methods of this survey are assembled from a small number of choices. Each method fixes a state space, chooses dynamics that realise a curve of laws on it, chooses a path or a coupling between source and data, trains a network by some signal, and generates by integrating or simulating the learned dynamics. \Cref{tab:zoo} records these choices for the main methods discussed, grouped as in the text and with the section in which each is treated. Its rows show how few ingredients the field uses, and its columns show where the open choices lie. Almost every training signal in the table is a conditional-expectation regression or a relative entropy, on path space or, for likelihood-trained flows, on the state space. The exceptions minimise a discrepancy such as the MMD or a sliced Wasserstein distance (MMD and sliced-Wasserstein flows, inductive moment matching), a log-variance divergence or a PDE residual (the physics-informed, or PINN, loss of NETS), or maximise a reward (policy-gradient fine-tuning).

\clearpage
{\footnotesize
\setlength{\LTcapwidth}{\textwidth}
\setlength{\tabcolsep}{3pt}
\renewcommand{\arraystretch}{1.12}
\newcommand{\grp}[2]{\addlinespace[5pt]\multicolumn{6}{@{}l}{\normalsize\scshape #1\quad\footnotesize\normalfont\S\ref{#2}}\\*\addlinespace[2pt]}
\newcommand{\mth}[2]{#1\newline{\scriptsize\citealp{#2}}}
\begin{longtable}{@{}P{3.0cm}P{1.7cm}P{2.2cm}P{2.45cm}P{2.85cm}P{2.25cm}@{}}
\caption{The methods of the survey. Each row records the state space, the dynamics that realise the curve of laws, the path or coupling between source and data, the signal by which the network is trained, and the procedure that generates samples. In the state-space column, $a$ and $D$ are numbers of added coordinates, $N$ is the number of residues and $n$ the number of atoms, and $L$ and $K$ are the sequence length and the vocabulary size; other symbols are as in the sections indicated.}\label{tab:zoo}\\
\toprule
Method & Space & Dynamics & Path or coupling & Training signal & Generation \\
\midrule
\endfirsthead
\multicolumn{6}{@{}l}{\small\Cref{tab:zoo}, continued}\\
\addlinespace[2pt]
\toprule
Method & Space & Dynamics & Path or coupling & Training signal & Generation \\
\midrule
\endhead
\midrule
\multicolumn{6}{r@{}}{\small Continued on the next page}\\
\endfoot
\bottomrule
\endlastfoot
\grp{Normalising flows}{sec:nf}
\mth{Coupling flows: NICE, RealNVP, Glow}{dinh2014nice,dinh2016density,kingma2018glow} & $\R^d$ & composed coupling layers & none; a single map & maximum likelihood & one pass \\
\mth{Neural spline flows}{durkan2019neural} & $\R^d$ & rational-quadratic spline couplings or autoregressive layers & none & maximum likelihood & one pass (coupling form) \\
\mth{Autoregressive flows: IAF, MAF, NAF}{kingma2016improved,papamakarios2017masked,huang2018neural} & $\R^d$ & triangular (Knothe--Rosenblatt) layers & none & maximum likelihood (MAF, NAF); variational bound (IAF; NAF also in VAEs) & one pass (IAF); $d$ sequential passes (MAF); $d$ sequential numerical inversions (NAF) \\
\mth{Invertible residual flows}{behrmann2019invertible,chen2019residual} & $\R^d$ & residual layers, $\Lip(g)<1$ & none & likelihood with power-series log-determinant & fixed-point inversion \\
\mth{Convex potential flows}{huang2021convex} & $\R^d$ & Brenier maps $\nabla\varphi$ & none & maximum likelihood & inversion of $\nabla\varphi$ by convex optimisation \\
\mth{Boltzmann generators}{noe2019boltzmann} & $\R^d$ & composed diffeomorphisms & none & reverse relative entropy and likelihood & one pass with reweighting \\
\mth{Stochastic normalising flows}{wu2020stochastic} & $\R^d$ & flow layers and Langevin steps & forward and backward chains & path-space relative entropy & Markov chain \\
\grp{Continuous normalising flows}{sec:cnf}
\mth{Neural ODEs, FFJORD}{chen2018neural,grathwohl2019ffjord} & $\R^d$ & ODE & none & likelihood with trace estimation & ODE solve \\
\mth{Augmented neural ODEs}{dupont2019augmented} & $\R^{d+a}$ & augmented ODE & none & studied for representation, not density & ODE solve \\
\mth{Kinetic regularisation, OT-Flow}{finlay2020train,onken2021otflow} & $\R^d$ & ODE; potential velocity in OT-Flow & least kinetic energy & likelihood plus transport penalty & ODE solve \\
\mth{Hamiltonian generative networks}{toth2020hamiltonian} & $\R^{2d}$ & Hamiltonian flow & none & likelihood bound; no log-determinant & ODE solve \\
\grp{Simulation-free flows}{sec:simfree}
\mth{Flow matching}{lipman2023flow} & $\R^d$ & ODE & Gaussian conditional paths, independent coupling & conditional velocity regression & ODE solve \\
\mth{Rectified flow and reflow}{liu2023flow} & $\R^d$ & ODE & straight interpolation; induced coupling & velocity regression, repeated & ODE solve, few steps \\
\mth{$c$-rectified flow}{liu2022rectified} & $\R^d$ & ODE, velocity $\nabla c^*(\nabla f_t)$ (a gradient for quadratic cost) & straight interpolation & Bregman-type matching (least squares for quadratic cost) & ODE solve \\
\mth{Stochastic interpolants}{albergo2023building,albergo2023stochastic} & $\R^d$ & ODE or SDE family & any coupling plus latent noise & velocity and score regression & ODE or SDE, tunable noise \\
\mth{Minibatch optimal-transport flow matching}{tong2024improving,pooladian2023multisample} & $\R^d$ & ODE & minibatch optimal coupling & velocity regression & ODE solve \\
\mth{Optimal flow matching}{kornilov2024optimal} & $\R^d$ & ODE, convex potential & straight optimal displacement & velocity regression & one step, $x_1=\nabla\Psi(x_0)$ \\
\mth{Action matching}{neklyudov2023action} & $\R^d$ & ODE, gradient velocity & marginal samples only & variational action & ODE solve \\
\mth{Generator matching}{holderrieth2025generator} & any & Markov process & conditional generators & generator regression & Markov simulation \\
\mth{Poisson flows: PFGM, PFGM++}{xu2022poisson,xu2023pfgm} & augmented $\R^{d+D}$ & electric field lines & physical field & field regression & ODE solve \\
\mth{Holographic generative flows}{mirafzali2026holographic} & Fourier phase space of a bulk scalar in AdS & Klein--Gordon backbone and residual & linear or Hermite interpolation from Gaussian bulk modes to lifted data & flow matching in radial phase space & ODE solve \\
\grp{Flow maps and few-step generation}{sec:flowmaps}
\mth{Consistency models}{song2023consistency} & $\R^d$ & one-time map & probability-flow trajectories & consistency loss & one or few steps \\
\mth{Consistency trajectory models}{kim2024consistency} & $\R^d$ & two-time map & probability-flow trajectories & consistency across pairs & any-step \\
\mth{Truncated consistency models}{lee2025truncated} & $\R^d$ & one-time map on $[t_{\rm tr},T]$ & probability-flow trajectories & consistency loss with a boundary model & one or two steps \\
\mth{Shortcut models}{frans2025one} & $\R^d$ & two-time map & dyadic steps & semigroup loss & one or few steps \\
\mth{MeanFlow, improved and pixel MeanFlow}{geng2025mean,geng2025improved,lu2026one} & $\R^d$ & average velocity & conditional paths & MeanFlow identity & one step \\
\mth{Flow-map self-distillation}{boffi2025build} & $\R^d$ & two-time map & learned velocity & Lagrangian, Eulerian or progressive loss & one or few steps \\
\mth{Inductive moment matching}{zhou2025inductive} & $\R^d$ & two-time map & pairs of times & kernel discrepancy & few steps \\
\mth{Progressive distillation}{salimans2022progressive} & $\R^d$ & student sampler & teacher trajectories & regression, halving steps & few steps \\
\mth{Distribution matching, score identity distillation}{yin2024one,zhou2024score} & $\R^d$ & one-step generator & pretrained diffusion & score difference; score identities & one step \\
\mth{Transition matching}{shaul2025transition} & $\R^d$ & Markov chain with learned kernels & linear interpolation, difference latent & flow-matching head per step & 16--32 stochastic steps, fewer than flow matching \\
\mth{Strong stochastic flow maps}{mccallum2026strong} & $\R^d$ & It\^o map of an additive-noise SDE & Brownian path as input & simulation-free matching to the SDE & few stochastic steps \\
\mth{Meta flow maps, Diamond Maps}{potaptchik2026meta,holderrieth2026diamond} & $\R^d$ & stochastic flow map to endpoint posteriors & posterior-targeting flows indexed by $(t,x)$ & amortised flow-map training or distillation & one-step posterior samples for reward alignment \\
\grp{Gradient-flow generative models}{sec:gradflows}
\mth{JKO-trained flows, JKO-iFlow}{mokrov2021large,xu2023normalizing} & $\R^d$ & Wasserstein gradient flow & JKO steps & proximal relative entropy & block-wise transport \\
\mth{MMD flows, sliced MMD flows}{arbel2019maximum,hertrich2024generative} & $\R^d$ & particle flow & discrepancy to data & MMD witness gradient & particle evolution \\
\mth{Sliced-Wasserstein flows}{liutkus2019sliced} & $\R^d$ & interacting-particle SDE & one-dimensional transports & sliced Wasserstein plus entropy & particle evolution \\
\mth{Stein variational gradient descent}{liu2016stein} & $\R^d$ & kernelised gradient flow & Stein geometry & relative entropy & particle evolution \\
\mth{Wasserstein--Fisher--Rao dynamics}{lu2019accelerating,yan2024learning} & $\R^d$ & transport with birth and death & Wasserstein--Fisher--Rao & relative entropy (sampling); negative log-likelihood (mixture estimation) & interacting particles \\
\mth{Drifting models; long--short flow-map derivation}{deng2026generative,li2026longshort} & $\R^d$ & one-step generator drifted during training & attraction to data, repulsion from samples & regression onto drifted samples & one pass \\
\mth{Equilibrium matching}{wang2025equilibrium} & $\R^d$ & time-independent gradient field & interpolation, target vanishing at the data & regression onto $(\epsilon-x)c(\gamma)$ & gradient descent \\
\mth{W-Flow}{han2026one} & $\R^d$ & one-step generator trained along a gradient flow & Wasserstein gradient flow of a Sinkhorn divergence & regression onto Euler-step targets & one step \\
\grp{Diffusion models}{sec:diffusion}
\mth{DDPM}{ho2020denoising} & $\R^d$ & variance-preserving chain & Gaussian noising & noise prediction (reweighted ELBO) & ancestral sampling \\
\mth{Noise-conditional score networks}{song2019generative} & $\R^d$ & variance-exploding noise levels & Gaussian noising & denoising score matching & annealed Langevin \\
\mth{Score SDEs}{song2021score} & $\R^d$ & VP or VE SDE & Gaussian noising & denoising score matching & reverse SDE or probability-flow ODE \\
\mth{Variational diffusion models}{kingma2021variational} & $\R^d$ & SDE with learned schedule & Gaussian noising & continuous-time ELBO & ancestral sampling \\
\mth{EDM}{karras2022elucidating} & $\R^d$ & probability-flow ODE, $\sigma(t)=t$ & Gaussian noising & weighted denoising loss, preconditioned denoiser & Heun ODE solver \\
\mth{Critically damped Langevin diffusion}{dockhorn2022score} & $\R^{2d}$ & hypoelliptic SDE & noise in the velocity & hybrid score matching & splitting sampler \\
\mth{Malliavin score targets}{mirafzali2025malliavin,mirafzali2025approach} & $\R^d$ & nonlinear SDE & nonlinear drift; state-independent noise, state-dependent in theory & regression onto Skorokhod-integral targets & reverse SDE \\
\mth{Divergence--kernel estimators}{ni2025scores,ni2025divergence} & $\R^d$ & multiplicative-noise SDE & forward simulation from a Gaussian & scores and density gradients from a forward covector process; relative entropy to the data & forward simulation of the fitted SDE \\
\mth{Latent diffusion}{rombach2022high} & latent space & linear SDE & noising in the latent & denoising score matching & decode after sampling \\
\mth{Heat dissipation, cold diffusion}{rissanen2023generative,bansal2023cold} & $\R^d$ & heat equation with small additive noise; arbitrary deterministic degradations & coarse-graining; blur, masking and other transforms & restoration regression & iterative restoration, with noise injected for heat dissipation \\
\mth{Blurring diffusion}{hoogeboom2023blurring} & $\R^d$ & Gaussian diffusion in frequency space & blurring plus Gaussian noise, with a frequency-dependent schedule & noise prediction & ancestral sampling \\
\grp{Score estimation}{sec:score-estimation}
\mth{Target score matching}{debortoli2024target} & $\R^d$ & linear SDE & Gaussian noising & data-score regression & reverse SDE \\
\grp{Samplers}{sec:samplers}
\mth{DDIM}{song2021denoising} & $\R^d$ & probability-flow ODE & non-Markovian forward & unchanged & first-order deterministic \\
\mth{DPM-Solver, DPM-Solver++, DEIS, UniPC}{lu2022dpm,lu2022dpmpp,zhang2023fast,zhao2023unipc} & $\R^d$ & semilinear ODE & half log-SNR & unchanged & exponential integrator \\
\mth{Predictor--corrector}{song2021score} & $\R^d$ & reverse SDE & Langevin corrections & unchanged & noise-corrected SDE \\
\mth{Restart sampling}{xu2023restart} & $\R^d$ & probability-flow ODE & large forward noise injections & unchanged & ODE with restart cycles \\
\mth{Parallel Picard sampling}{shih2023parallel} & $\R^d$ & sliding window of time steps & Picard iteration & unchanged & parallel in time \\
\grp{Conditioning}{sec:conditioning}
\mth{Classifier guidance}{dhariwal2021diffusion} & $\R^d$ & guided reverse SDE & $h$-transform by a classifier, exact only at unit guidance scale & noise-conditional classifier & guided sampling \\
\mth{Classifier-free guidance}{ho2022classifier} & $\R^d$ & guided reverse SDE & extrapolated scores & joint conditional and unconditional training & guided sampling \\
\mth{Diffusion posterior sampling}{chung2023diffusion} & $\R^d$ & guided reverse SDE & likelihood at the denoiser & pretrained prior & guided sampling \\
\mth{Pseudoinverse guidance, DDRM}{song2023pseudoinverse,kawar2022denoising} & $\R^d$ & guided reverse dynamics & Gaussian posterior; singular basis & pretrained prior & guided sampling \\
\mth{Twisted diffusion sampler}{wu2023practical} & $\R^d$ & reverse SDE with particles & twisting potentials & pretrained prior & sequential Monte Carlo \\
\mth{Feynman--Kac steering}{singhal2025general} & $\R^d$ or finite set & reverse process with particles & reward-tilted law & pretrained model, intermediate rewards & resampling with potentials \\
\mth{DEFT, Doob's matching}{denker2024deft,chang2026inference} & $\R^d$ & $h$-transformed SDE & learned $h$-transform & fine-tuned $h$ network; gradient-regularised regression & guided sampling \\
\mth{Feynman--Kac correctors}{skreta2025feynman} & $\R^d$ & weighted guided SDE & annealed, product or geometric-average targets & pretrained scores & sequential Monte Carlo \\
\mth{Variational flow maps}{mammadov2026variational} & $\R^d$ & one-step flow map & conditional noise from a learned adapter & joint variational bound with a mean-flow loss & one or few steps \\
\grp{Schr\"odinger bridges}{sec:sb}
\mth{Diffusion Schr\"odinger bridges}{debortoli2021diffusion} & $\R^d$ & SDE & iterative proportional fitting & half-bridge regressions & SDE \\
\mth{Bridge matching, iterative Markovian fitting}{shi2023diffusion,peluchetti2023diffusion} & $\R^d$ & SDE & Markovian and reciprocal projections & bridge matching & SDE \\
\mth{Twisted Schr\"odinger bridge matching}{noble2026twisted} & $\R^d$ & SDE & Feynman--Kac reference with a potential & bridge matching with the potential gradient & SDE \\
\mth{Categorical Schr\"odinger bridge matching}{ksenofontov2025categorical} & finite spaces & discrete-time Markov chain & reciprocal and Markovian projections & relative-entropy bridge matching & Markov chain \\
\mth{Simulation-free score and flow matching}{tong2024simulation} & $\R^d$ & SDE & entropic minibatch coupling, Brownian bridges & score and drift regression & SDE \\
\mth{I2SB, denoising diffusion bridges}{liu2023i2sb,zhou2024denoising} & $\R^d$ & SDE & paired data & bridge score matching & SDE \\
\mth{Light Schr\"odinger bridges}{korotin2024light} & $\R^d$ & SDE with closed-form drift & Gaussian-mixture potentials & relative entropy & direct sampling of the Gaussian-mixture conditional; SDE for paths \\
\grp{Stochastic optimal control}{sec:control}
\mth{Path integral sampler}{zhang2022path} & $\R^d$ & controlled SDE & Dirac source to target & path-space relative entropy & SDE \\
\mth{Denoising diffusion samplers}{vargas2023denoising} & $\R^d$ & controlled SDE & Ornstein--Uhlenbeck reference & path-space relative entropy & SDE \\
\mth{Log-variance training}{richter2024improved} & $\R^d$ & controlled SDE & off-policy trajectories & log-variance divergence & SDE \\
\mth{Iterated denoising energy matching}{akhoundsadegh2024iterated} & $\R^d$ & reverse SDE & noised target & Monte Carlo score from the energy & SDE \\
\mth{Stochastic localisation via iterative posterior sampling}{grenioux2024stochastic} & $\R^d$ & stochastic localisation SDE & observation process & MCMC estimate of the denoiser, no training & SDE \\
\mth{Controlled Monte Carlo diffusions, NETS}{vargas2024transport,albergo2025nets} & $\R^d$ & annealed SDE & prescribed annealing path & path-space relative entropy (CMCD); PINN or action-matching loss (NETS) & SDE with Crooks--Jarzynski weights \\
\mth{Adjoint matching, adjoint sampling}{domingoenrich2025adjoint,havens2025adjoint} & $\R^d$ & controlled SDE & memoryless reference & lean adjoint regression & SDE \\
\mth{Trust-region stochastic optimal control}{blessing2025trust} & $\R^d$ & controlled SDE & geometric annealing in path space & trust-region log-variance or SOC matching & SDE \\
\mth{Underdamped diffusion bridges}{blessing2025underdamped} & $\R^{2d}$ & controlled kinetic SDE & learned forward and backward processes & path-space relative entropy & splitting integrator \\
\mth{Policy-gradient fine-tuning: DDPO, DPOK}{black2024training,fan2023dpok} & $\R^d$ & fine-tuned sampler & pretrained model & reward, with or without relative entropy & stochastic sampler with Gaussian steps \\
\mth{Flow-GRPO}{liu2025flowgrpo} & $\R^d$ & marginal-preserving SDE & pretrained flow model & group-relative policy gradient & default sampler and step count \\
\grp{Riemannian manifolds}{sec:riemannian}
\mth{Riemannian CNFs}{mathieu2020riemannian,lou2020neural} & manifold & ODE on $\M$ & none & likelihood, Riemannian divergence & manifold ODE \\
\mth{Moser flow}{rozen2021moser} & manifold & ODE on $\M$ & prior minus divergence & likelihood without ODE simulation & manifold ODE \\
\mth{Riemannian score-based models}{debortoli2022riemannian} & manifold & Brownian motion on $\M$ & heat-kernel noising & implicit or heat-kernel score matching & geodesic random walk \\
\mth{Riemannian flow matching}{chen2024flow} & manifold & ODE on $\M$ & premetric conditional flows & velocity regression & manifold ODE \\
\mth{Reflected, mirror and barrier diffusions}{lou2023reflected,liu2023mirror,fishman2023diffusion} & constrained domain & reflected or dual SDE & boundary-respecting noise & score matching & constrained SDE \\
\mth{Metric flow matching}{kapusniak2024metric} & $\R^d$ & ODE & data-dependent geodesics & velocity regression & ODE solve \\
\grp{Lie groups and equivariance}{sec:lie}
\mth{Denoising on $SO(3)$}{leach2022denoising} & $SO(3)$ & discrete-time chain with $\mathrm{IG}_{SO(3)}$ transitions & isotropic Gaussian $\mathrm{IG}_{SO(3)}$ & noise prediction in $\mathfrak{so}(3)$ & ancestral sampling from the uniform law \\
\mth{Frame diffusion and flow matching}{yim2023se3,yim2023fast,bose2024se3} & $SE(3)^N$ & SDE or ODE on frames & frame noising or flows & score or velocity regression & reverse dynamics \\
\mth{Torsional diffusion}{jing2022torsional} & torus & SDE on torsions & wrapped noising & denoising score matching & reverse SDE \\
\mth{Equivariant diffusion for molecules}{hoogeboom2022equivariant} & zero-mean $\R^{3n}$ coordinates with atom features & discrete-time Gaussian chain with $E(3)$-equivariant denoiser & Gaussian noising & noise-prediction loss & ancestral sampling \\
\mth{Equivariant flow matching}{klein2023equivariant} & $\R^{3n}$ & equivariant ODE & aligned minibatch coupling & velocity regression & ODE solve \\
\grp{Discrete and statistical-manifold spaces}{sec:discrete}
\mth{Structured discrete diffusion}{austin2021structured} & finite set & discrete-time chain & uniform, absorbing or banded & variational bound & ancestral sampling \\
\mth{Continuous-time discrete diffusion}{campbell2022continuous} & finite set & CTMC & rate matrices & continuous-time ELBO & $\tau$-leaping \\
\mth{Score entropy}{lou2024discrete} & finite set & CTMC & rate matrices & score entropy & reverse CTMC \\
\mth{Masked diffusion}{sahoo2024simple,shi2024simplified} & finite set & absorbing CTMC & masking & weighted cross-entropy & unmasking \\
\mth{Masked diffusion language models: LLaDA, Dream}{nie2025llada,ye2025dream} & token sequences & absorbing CTMC & masking & weighted cross-entropy & parallel unmasking \\
\mth{Block diffusion}{arriola2025block} & token sequences & absorbing CTMC within blocks & masking & block-wise weighted cross-entropy, clipped schedules & block-autoregressive with caching \\
\mth{Nemotron-Labs-Diffusion}{fu2026nemotron} & token sequences & absorbing CTMC within blocks, with next-token prediction & masking within blocks & joint autoregressive and block-diffusion losses & autoregressive, block diffusion or self-speculation \\
\mth{Energy-based diffusion language models}{xu2025energy} & token sequences & absorbing CTMC & masking & residual energy, by NCE or from an autoregressive model & importance resampling per step \\
\mth{DiffusionGemma}{diffusiongemma2026} & token sequences & uniform-noise CTMC & blockwise, conditioned on the prefix & fine-tuning, then reinforcement learning and sampler distillation & parallel denoising of 256-token blocks \\
\mth{Duo}{sahoo2025diffusion} & token sequences & uniform-state CTMC & argmax of Gaussian diffusion & curriculum; discrete consistency distillation & few-step sampling \\
\mth{Flow map language models}{lee2026flowmap} & one-hot token embeddings & ODE and its flow map & interpolant to one-hot data & cross-entropy on the simplex & one or few steps \\
\mth{Discrete flow maps}{potaptchik2026discrete} & one-hot token sequences in $\R^{L\times K}$ & flow map through the mean denoiser, which lies on the simplex & linear interpolant from noise in $\R^K$ to one-hot vertices & cross-entropy and relative-entropy flow-map losses & one or few steps \\
\mth{Sigma}{yang2026sigma} & learned 16-dimensional token embeddings & Gaussian diffusion & blockwise, conditioned on the prefix & likelihood bound, autoregressive warm start & ODE or SDE per block \\
\mth{Discrete flow matching}{gat2024discrete,campbell2024generative} & finite set & CTMC & mixture paths & posterior cross-entropy & CTMC simulation from source to data \\
\mth{Dirichlet flow matching}{stark2024dirichlet} & simplex & ODE & Dirichlet paths & posterior cross-entropy & ODE solve \\
\mth{Fisher flow matching}{davis2024fisher} & simplex as sphere & ODE & spherical geodesics & Riemannian velocity regression & ODE solve \\
\grp{Infinite dimensions}{sec:infinite}
\mth{Denoising diffusion operators}{lim2023score} & Hilbert space & Gaussian-random-field noise levels & trace-class noise & operator score matching & annealed Langevin \\
\mth{Infinite-dimensional diffusion}{kerrigan2023diffusion,pidstrigach2024infinite} & Hilbert space & discrete-time Gaussian chain or SDE in $\cH$ & Gaussian measures & reweighted ELBO; score matching in function space & reverse process \\
\mth{Functional flow matching}{kerrigan2024functional} & Hilbert space & ODE & Gaussian measure to data & velocity regression & ODE solve \\
\mth{Functional mean flow}{li2025functional} & Hilbert space & two-time flow map & Gaussian measure to data & mean-flow identity & one step \\
\mth{Malliavin scores of SPDE laws}{mirafzali2025infinite,mirafzali2026logarithmic} & Hilbert space & linear or variational SPDE & coloured additive noise, or state-dependent noise for variational SPDEs & logarithmic-derivative formulas, in closed form for linear SPDEs; no training & none; score formulas only \\
\grp{L\'evy and jump processes}{sec:levy}
\mth{L\'evy--It\^o models}{yoon2023score} & $\R^d$ & $\alpha$-stable SDE & heavy-tailed noising & fractional score matching & jump SDE \\
\mth{Denoising L\'evy probabilistic models}{shariatian2025denoising} & $\R^d$ & $\alpha$-stable chain & heavy-tailed noising & modified DDPM loss & ancestral sampling \\
\mth{Trans-dimensional jump diffusions}{campbell2023trans} & varying dimension & jump diffusion & component deletion & score and jump-rate learning & jump simulation \\
\grp{Field theory}{sec:fieldtheory}
\mth{Flow-based lattice sampling}{albergo2019flow,kanwar2020equivariant,boyda2021sampling} & lattice fields & flows, gauge-equivariant for gauge theories & none & reverse relative entropy & independence Metropolis \\
\mth{Diffusion as stochastic quantisation}{wang2024quantization} & lattice fields & reverse SDE & Gaussian noising & denoising score matching & global sampler \\
\mth{Inverse renormalisation diffusion}{cotler2023renormalizing,masuki2025generative} & fields, images & reverse RG flow & coarse-graining & denoising score matching & coarse to fine \\
\mth{Renormalisation group flow matching}{masuki2026renormalization} & fields, images & ODE along an exact RG flow & Polchinski flow with heat-kernel cutoff & flow matching with local velocities & patchwise local generation \\
\mth{Wavelet score-based models}{guth2023wavelet} & $\R^d$ & SDE per scale & wavelet factorisation & conditional score matching & scale by scale \\
\end{longtable}
}

\section{Open problems}\label{sec:outlook}

The open problems collected at the end of Parts~II to VI fall under five themes, which are also the organising ideas of the survey, namely integration by parts, the geometry of paths, changes of measure, the state space, and physics.

\paragraph{Variance-optimal integration by parts.}
The score targets used to train diffusion and bridge models are unbiased conditional estimators obtained by integration by parts (\cref{tab:ibp}). The choice among them is a free design variable, whether it be the covering field on Wiener space, the free parameter of divergence--kernel estimators, or a combination of noise-side and data-side targets with weights summing to one. A theory that selects these to minimise variance, possibly by learning the covering field jointly with the score, would sharpen sample complexity where diffusion models are weakest, at small noise. Making Malliavin targets practical for nonlinear, state-dependent, kinetic and jump-driven forward processes is the computational side of the same problem.

\paragraph{The geometry of paths in $\cP(\M)$.}
Which curve of laws should a model follow? Optimal transport gives straight but hard-to-estimate paths, diffusion gives smooth but curved ones, and the probability flow is not optimal transport \citep{lavenant2022flow}. Thermodynamic speed limits \citep{ikeda2025speed}, straightness theory for rectified flow \citep{bansal2025wasserstein} and the multiscale analysis of wavelet scores \citep{guth2023wavelet} each address part of the question. A unified criterion trading learnability of the velocity against integration cost is missing, as are intrinsic-dimension guarantees for one-step maps trained by consistency or self-distillation losses.

\paragraph{Changes of measure on path space.}
Guidance, bridges, control, fine-tuning and stochastic samplers all reweight path measures (\cref{sec:doob,sec:follmer}), at least when the new path law is absolutely continuous with respect to the original. Classifier-free guidance is a reweighting that is in general not an $h$-transform, and the law it samples is understood only in special settings; on continuous state spaces, iterative Markovian fitting has geometric rates only for weakly log-concave marginals and long horizons, and its neural implementations and neural IPF lack guarantees in path space with errors controlled by the training losses; and the off-policy objectives that train control-based samplers lack guarantees in high dimension. Malliavin calculus, as the calculus of path-space sensitivities, is a natural source of low-variance gradient estimators for all of them.

\paragraph{Beyond Euclidean state spaces.}
Heat kernels on general manifolds, admissible directions and noise covariances in infinite dimensions, a sharp convergence theory for discrete and jump processes, and intrinsic definitions of scores for SPDE laws \citep{mirafzali2026logarithmic} are each partly solved. Connecting function-space scores to neural operators with discretisation-independent guarantees would make generative modelling of physical fields as principled as image generation.

\paragraph{Physics as a source of structure.}
Phase-transition analyses explain memorisation for exact empirical scores; for trained networks there are so far training-time accounts \citep{bonnaire2025why} and an associative-memory account of the transition from memorisation to generalisation \citep{pham2025memorization}, but there is no theory of speciation and collapse for networks with the inductive biases of \cref{sec:generalisation}. Renormalisation suggests forward processes adapted to correlation structure, and holography suggests that a radial, scale-ordered transport with an explicit propagator can organise generation \citep{mirafzali2026holographic}. Turning these analogies into design principles with guarantees is the most open of the five themes.

\section{Conclusion}\label{sec:conclusion}

The rapid progress of flow and diffusion models has been driven by simple regressions run at scale. The mathematics reviewed here explains why those regressions work. Their targets are random variables whose conditional expectations are the quantities of interest, made computable by marginalisation and integration by parts. The same mathematics indicates where the next gains may come from, namely better choices of paths in the space of measures, better estimators on path space, better geometry for the state space, and the structural insight of physics. We hope that the three geometries, of the state space, of the space of probability measures and of path space, help readers move between these levels, and between communities that have often studied the same objects independently.

\appendix
\section{A dictionary of parameterisations}\label{app:dictionary}

A diffusion network may be asked to predict the score, the noise, the clean data or a velocity. Given $x_t$ these are equivalent. This appendix records the conversions between them and lists the paths most often used (\cref{tab:schedules}). Throughout, $x_t=\alpha_tx_0+\sigma_t\epsilon$ with $\epsilon\sim\cN(0,\Id)$ independent of $x_0$, and hats denote network predictions of conditional expectations given $x_t$, so that $\hat x_0\approx\E[x_0\mid x_t]$, $\hat\epsilon\approx\E[\epsilon\mid x_t]$, $\hat s\approx\nabla\log p_t(x_t)$. By Tweedie's formula (\cref{prop:tweedie}) the predictions are affinely related, as listed in \cref{tab:dictionary}. The velocity of the probability-flow ODE in the same variables is
\begin{align*}
\dot x_t&=\dot\alpha_t\,\hat x_0+\dot\sigma_t\,\hat\epsilon\\
&=\frac{\dot\alpha_t}{\alpha_t}\,x_t+\Big(\dot\sigma_t-\frac{\dot\alpha_t}{\alpha_t}\sigma_t\Big)\hat\epsilon,
\end{align*}
which agrees with $f(t)x-\frac12g(t)^2\hat s$ for $f=\dot\alpha_t/\alpha_t$ and $g^2=2\sigma_t\dot\sigma_t-2(\dot\alpha_t/\alpha_t)\sigma_t^2$; this is also the marginal flow-matching velocity generated by the conditional Gaussian paths (\cref{prop:marginalisation}).

\begin{table}[ht]
\centering\small
\begin{tabular}{@{}llll@{}}
\toprule
Target & from $\hat\epsilon$ & from $\hat x_0$ & from $\hat s$ \\
\midrule
score $\hat s$ & $-\hat\epsilon/\sigma_t$ & $(\alpha_t\hat x_0-x_t)/\sigma_t^2$ & $\hat s$ \\
noise $\hat\epsilon$ & $\hat\epsilon$ & $(x_t-\alpha_t\hat x_0)/\sigma_t$ & $-\sigma_t\hat s$ \\
data $\hat x_0$ & $(x_t-\sigma_t\hat\epsilon)/\alpha_t$ & $\hat x_0$ & $(x_t+\sigma_t^2\hat s)/\alpha_t$ \\
\bottomrule
\end{tabular}
\caption{Conversions between prediction targets for Gaussian paths. Squared losses transform accordingly, with $\norm{\hat x_0-x_0}^2=\frac{\sigma_t^2}{\alpha_t^2}\norm{\hat\epsilon-\epsilon}^2$ and $\norm{\hat s+\epsilon/\sigma_t}^2=\sigma_t^{-2}\norm{\hat\epsilon-\epsilon}^2$, so the choice of target is a choice of loss weighting in terms of $\snr(t)=\alpha_t^2/\sigma_t^2$. The $v$-prediction $\hat v=\alpha_t\hat\epsilon-\sigma_t\hat x_0$ is given in \eqref{eq:vpred}.}
\label{tab:dictionary}
\end{table}

The $v$-prediction target of \citet{salimans2022progressive} converts as
\begin{equation}\label{eq:vpred}
\begin{aligned}
\hat v&=\frac{\alpha_t^2+\sigma_t^2}{\alpha_t}\,\hat\epsilon-\frac{\sigma_t}{\alpha_t}\,x_t\\
&=\frac{\alpha_tx_t-(\alpha_t^2+\sigma_t^2)\,\hat x_0}{\sigma_t}\\
&=-\frac{(\alpha_t^2+\sigma_t^2)\,\sigma_t}{\alpha_t}\,\hat s-\frac{\sigma_t}{\alpha_t}\,x_t.
\end{aligned}
\end{equation}

\begin{table}[ht]
\centering\small
\begin{tabularx}{\textwidth}{@{}P{3.6cm}P{5.2cm}L@{}}
\toprule
Path & $(\alpha_t,\sigma_t)$ & Notes \\
\midrule
Variance preserving (DDPM) & $\alpha_t=e^{-\frac12\int_0^t\beta}$, $\sigma_t^2=1-\alpha_t^2$ & \citet{ho2020denoising,song2021score} \\
Variance exploding (NCSN) & $\alpha_t=1$, $\sigma_t$ increasing & \citet{song2019generative,song2021score} \\
EDM & $\alpha_t=1$, $\sigma_t=t$ & \citet{karras2022elucidating} \\
Cosine & $\alpha_t=\cos(\frac{\pi t}{2})$, $\sigma_t=\sin(\frac{\pi t}{2})$, $t\in[0,1]$ & \citet{nichol2021improved} (discrete version with a small offset, $\alpha_t^2=f(t)/f(0)$, $f(t)=\cos^2\big(\frac{t+s}{1+s}\frac{\pi}{2}\big)$, $s=0.008$) \\
Straight (FM-OT, rectified flow) & $\alpha_t=1-t$, $\sigma_t=t$ & \citet{lipman2023flow,liu2023flow}; time reversed \\
\bottomrule
\end{tabularx}
\caption{Common Gaussian paths in the diffusion orientation (data at $t=0$). Flow-matching papers usually reverse time so that data sit at $t=1$. In the first row $\beta$ is the noise schedule of the variance-preserving SDE (\cref{sec:diffusion}); in the cosine row $f$ is a schedule function and $s$ an offset, not a drift and a time.}
\label{tab:schedules}
\end{table}

\section{Forward kernels on standard spaces}\label{app:kernels}

\Cref{tab:kernels} lists forward processes whose transition kernels (or good approximations) are available in closed form, which is what denoising score matching and bridge matching require. Where no closed form exists, the integration-by-parts estimators of \cref{sec:ibp} and learned heat kernels \citep{ko2026riemannian} are the alternatives.

\begin{table}[ht]
\centering\small
\begin{tabularx}{\textwidth}{@{}P{2.4cm}P{3.4cm}L@{}}
\toprule
Space & Forward process & Kernel or score target \\
\midrule
$\R^d$ & OU / VP, VE & $\cN(\alpha_tx_0,\sigma_t^2\Id)$; score $-(x-\alpha_tx_0)/\sigma_t^2$ \\
$\R^{2d}$ (kinetic) & critically damped Langevin & joint Gaussian in position and velocity with explicit mean and covariance \citep{dockhorn2022score} \\
Torus $\mathbb T^d$ & Brownian motion & wrapped normal $\sum_{k\in\mathbb Z^d}\cN(x+2\pi k;\,x_0,\sigma_t^2\Id)$ \citep{jing2022torsional} \\
Sphere $S^d$ & Brownian motion & spherical-harmonic (Gegenbauer) series; Varadhan approximation for small $t$ \citep{debortoli2022riemannian} \\
$SO(3)$ & Brownian motion & $\mathrm{IG}_{SO(3)}$, whose density with respect to Haar measure is $f(\omega,t)=\sum_\ell(2\ell+1)e^{-\ell(\ell+1)t/2}\frac{\sin((\ell+1/2)\omega)}{\sin(\omega/2)}$ in the rotation angle $\omega$ at time $t$ (generator $\frac12\Delta$) \citep{yim2023se3}; \citet{leach2022denoising} write $\epsilon^2$ for $t/2$ (the angle itself has density $\propto(1-\cos\omega)f$) \\
$SE(3)$ & product of $SO(3)$ and $\R^3$ diffusions & $\mathrm{IG}_{SO(3)}\times$ Gaussian on translations \citep{yim2023se3} \\
Hyperbolic $\mathbb H^3$ & Brownian motion (generator $\frac12\Delta$) & $(2\pi t)^{-3/2}\frac{r}{\sinh r}e^{-t/2-r^2/(2t)}$ with $r$ the geodesic distance \citep{grigoryan1998heat}, who state it for $\partial_tu=\Delta u$ \\
Simplex $\Delta^{K-1}$ & independent univariate Jacobi diffusions in stick-breaking coordinates & Dirichlet stationary law; each univariate kernel a Jacobi-polynomial series \citep{avdeyev2023dirichlet} \\
Finite set & CTMC, uniform or masking & $e^{tQ}$ in closed form, a mixture of identity and uniform, or of identity and mask \citep{austin2021structured,campbell2022continuous} \\
Hilbert space $\cH$ & $\dd X=-\frac12X\dd t+\sqrt C\,\dd W$, $C$ trace class & $\cN\big(e^{-t/2}x_0,(1-e^{-t})C\big)$; target the preconditioned score $s(t,x)=-(x-e^{-t/2}\E[x_0\mid x_t=x])/(1-e^{-t})$, equal to $C\nabla\log p_t$ in finite dimensions \citep{pidstrigach2024infinite} \\
\bottomrule
\end{tabularx}
\caption{Forward kernels on standard state spaces. The hyperbolic kernel is for curvature $-1$. On $\mathbb H^3$ of radius $L$, as in \cref{sec:holography}, replace $r$ by $r/L$ and $t$ by $t/L^2$, and divide by $L^3$. In the last row, $W$ is a cylindrical Wiener process on $\cH$.}
\label{tab:kernels}
\end{table}

\clearpage
\bibliographystyle{plainnat}
\bibliography{References}

\end{document}